%% file: supnorm.tex
\documentclass[11pt,twoside]{article}
\newcommand{\ifims}[2]{#1} 
\newcommand{\ifAMS}[2]{#1}   
\newcommand{\ifau}[4]{#1}  
\newcommand{\ifbook}[2]{#1}   
\newcommand{\ifNL}[2]{#2}  
\newcommand{\ifapp}[2]{#2}  
\newcommand{\ifsqnorm}[2]{#2}  

\usepackage{tocloft}

\advance\cftsecnumwidth 0.3em\relax
\advance\cftsubsecnumwidth 0.3em\relax

\input myfrontmatter

\input mydef
\input statdef

\def\thetitle{Semiparametric plug-in estimation, sup-norm risk bounds, marginal optimization, and
inference in BTL model}
\def\theruntitle{Semiparametric plug-in estimation, sup-norm risk, and marginal optimization}

\def\theabstract{
The recent paper \cite{GSZ2023} on estimation and inference for top-ranking problem in Bradley-Terry-Lice (BTL) model
presented a surprising result: component-wise estimation and inference can be done under much weaker conditions 
on the number of comparison then it is required for the full dimensional estimation. 
The present paper revisits this finding from completely different viewpoint.
Namely, we show how a theoretical study of \emph{estimation in sup-norm} can be reduced to the analysis of \emph{plug-in semiparametric estimation}. 
For the latter, we adopt and extend the general approach from \cite{Sp2024} to high-dimensional estimation and inference.
The main tool of the analysis is a theory of \emph{perturbed marginal optimization} 
when an objective function depends on a low-dimensional target parameter along with a high-dimensional nuisance parameter.
A particular focus of the study is the critical dimension condition.
Full-dimensional estimation requires in general the condition \( \neff \gg \dimA \) between the effective parameter dimension \( \dimA \)
and the effective sample size \( \neff \) corresponding to the smallest eigenvalue of 
the Fisher information matrix \( \IF \).
Inference on the estimated parameter is even more demanding: the condition \( \neff \gg \dimA^{2} \) cannot be generally avoided;
see \cite{Sp2024}.
However, for the sup-norm estimation, 
the critical dimension condition can be reduced to
\( \neff \geq \CONST \log(\dimp) \). 
}

\def\kwdp{62F10,62E17}
\def\kwds{62J12}

\def\thekeywords{Fisher and Wilks expansions, risk bounds, critical dimension, BTL model}

\def\thankstitle{}

\aua
{Vladimir Spokoiny}
{Spokoiny, V. }
{
    Weierstrass Institute Berlin,  
    HSE 
    and IITP RAS,
    \\
    Mohrenstr. 39, 10117 Berlin, Germany
    }
{spokoiny@wias-berlin.de}
{Weierstrass-Institute Berlin}
{
  Financial support by the German Research Foundation (DFG) through the Collaborative Research Center 1294 ``Data assimilation'' is gratefully acknowledged.
}

\begin{document}
\thispagestyle{empty}
{
\title{\thetitle}
\theauthors

\maketitle
\begin{abstract}
{\footnotesize \theabstract}
\end{abstract}

\ifAMS
    {\par\noindent\emph{AMS 2010 Subject Classification:} Primary \kwdp. Secondary \kwds}
    {\par\noindent\emph{JEL codes}: \kwdp}

\par\noindent\emph{Keywords}: \thekeywords
} 

\tableofcontents

\Chapter{Introduction}
\label{Sgenintr}
Recent paper \cite{Sp2024} states some finite sample accurate results for the (quasi) maximum 
likelihood estimator under the so-called \emph{critical dimension} condition.
This condition means ``effective parameter dimension is smaller than effective sample size'',
where the \emph{effective dimension} \( \dimA \) can be defined as the trace of the covariance matrix of the standardized score,
while the effective sample size \( \neff \) corresponds to the smallest eigenvalue of the Fisher information matrix.
Parametric inference is even more demanding and requires \( \dimA^{2} \ll \neff \). 
%
In view of recent results from \cite{katsevich2023tight}, condition \( \dimA^{2} \ll \neff \) is essential and cannot be avoided.
Surprisingly, \cite{GSZ2023} established very strong inference results for the Bradley-Terry-Luce (BTL) problem of ranking of \( \dimp \) players 
from pairwise comparisons under much weaker condition \( \neff \gg \log(\dimp) \).
One of the aims of this paper is to revisit this problem and to understand the nature of this phenomen.
Compared to \cite{GSZ2023}, a different view on the considered problem is offered; the whole analysis is reduced 
to estimation in sup-norm which, in turn, is treated as a special case of semiparametric plug-in estimation.

We mention two possible motivations for the study of plug-in estimation.
Let a model be described by a high dimensional parameter \( \prmtv \) consisting of a low dimensional target \( \tarpv \) and a high dimensional 
nuisance parameter \( \nupv \).
Given a full dimensional log-likelihood function \( \LL(\prmtv) = \LL(\tarpv,\nupv) \), the full dimensional MLE is defined by a joint optimization
\begin{EQA}
	\tilde{\prmtv}
	& \eqdef &
	\argmax_{\prmtv} \LL(\prmtv)
	=
	\argmax_{(\tarpv,\nupv)} \LL(\tarpv,\nupv) .
\label{u8yda8yiiken5555t}
\end{EQA}
The \emph{profile MLE} \( \tilde{\tarpv} \) is the \( \tarpv \)-component of \( \tilde{\prmtv} \).
This estimator possesses a number of nice theoretical properties like asymptotic efficiency and normality; see e.g.
\cite{MR1245941}, \cite{robinson1988root}, \cite{Do1994}, \cite{van2000asymptotic}, \cite{ChCh2018}.
\cite{Sp2024} explained how a problem of high-dimensional nonlinear regression can be reduced to a semiparametric estimation
for a special class of \emph{stochastically linear smooth} (SLS) models and presented some finite sample properties 
of this estimator including Fisher and Wilks expansions.
Unfortunately, the full dimensional optimization is often hard to implement and study.
At the same time, partial optimization of \( \LL(\tarpv,\nupv) \) w.r.t. \( \tarpv \) for \( \nupv \) fixed can be a much simpler problem.
This suggests a \emph{plug-in} approach: given a preliminary/pilot estimate \( \hat{\nupv} \), define \( \hat{\tarpv} \) by maximization 
of \( \LL(\tarpv,\hat{\nupv}) \):
\begin{EQA}
	\hat{\tarpv}
	& \eqdef &
	\argmax_{\tarpv} \LL(\tarpv,\hat{\nupv}) .
\label{dhyiu3iwextq4e3esarg}
\end{EQA}
Furthermore, one can continue by using \( \hat{\tarpv} \) as a preliminary estimate of \( \tarpv \) and reestimate \( \hat{\nupv} \).
This leads to \emph{alternated optimization} and \emph{EM}-type methods: starting from \( \hat{\tarpv}_{1} \) and \( k=1 \), estimate
\begin{EQA}[rcccl]
	\hat{\nupv}_{k}
	&=&
	\argmax_{\nupv} \LL(\hat{\tarpv}_{k},\nupv) \, ,
	\qquad
	\hat{\tarpv}_{k+1}
	&=&
	\argmax_{\tarpv} \LL(\tarpv,\hat{\nupv}_{k}) \, ,
\label{dc7iu27b7w23uyikgvbb76}
\end{EQA}
and increase \( k \) by one; see e.g. 
\cite{Balakrishnan2014}, 
\cite{YiCa2015}, 
\cite{AASP2015}, %
\cite{JiWa2016}, 
\cite{Dwivedi2018SingularityMA}, 
\cite{Ho2020InstabilityCE}, 
\cite{WaDwi2020}, 
\cite{lin2023semi}, 
 and references therein.
Analysis of convergence for such procedures requires to evaluate the accuracy of each step.
This leads to a semiparametric plug-in study.

Another situation naturally leading to a plug-in semiparametric estimation arises when the MLE \( \tilde{\prmtv} \)
of the full dimensional parameter \( \prmtv \in \R^{\dimp} \)
is studied component-wise (in sup-norm);
see e.g.
\cite{LiRi2022} and
\cite{GSZ2023}
in context of estimation for a BTL model.
For each \( j \leq \dimp \), one can consider the \( j \)-th entry \( \prmts_{j} \) of \( \prmtvs \) as the target component 
while all the remaining components are treated as a nuisance. 
This leads to the sequential (linewise) optimization procedure: at each step, the objective function is optimized
w.r.t. only one parameter, while the remaining ones are kept.
Later in this paper, we study this approach in general and then apply to the BTL model.

\subsection*{This paper contribution}
For effective estimation and inference in the BTL model,
we bring together three very different topics:
semiparametric plugin estimation, 
coordinate-wise estimation with sup-norm loss function,
and marginal perturbed optimization. 
 
\paragraph{Semiparametric plug-in estimation.}
We explore a plug-in estimator \( \hat{\tarpv} \) from \eqref{dhyiu3iwextq4e3esarg}
without specifying the nature of the pilot \( \hat{\nupv} \) and structure of the 
objective function \( \LL(\tarpv,\nupv) \).
The only condition on this pilot is a local concentration around the true value \( \nupvs \).
In the parametric case with \( \nupvs \) known, one can use the following expansion for 
the estimator \( \tilde{\tarpv} = \argmax_{\tarpv} \LL(\tarpv,\nupvs) \):
\begin{EQA}[c]
	\tilde{\tarpv} - \tarpvs 
	=
	\IF^{-1} \nabla_{\tarpv} \zeta + \operatorname{remander} 
	\; ;
\end{EQA}
see \cite{Sp2024}.
Here \( \IF = - \nabla_{\tarpv\tarpv} \E \LL(\tarpvs,\nupvs) \) 
and \( \nabla_{\tarpv} \zeta \eqdef \nabla_{\tarpv} \LL(\tarpvs,\nupvs) \).
The use of a plug-in estimate \( \hat{\nupv} \) in place of \( \nupvs \) leads to an extra term in the expansion:
\begin{EQA}[c]
	\hat{\tarpv} - \tarpvs
	=
	\IF^{-1} \nabla_{\tarpv} \zeta
	- \IF^{-1} \IFT_{\tarpv\nuiv} (\hat{\nuiv} - \nuivs)
	+ \operatorname{remander} 
	\, ,
\end{EQA}
where \( \IFT_{\tarpv\nuiv} = - \nabla_{\tarpv\nupv} \E \LL(\tarpvs,\nupvs) \); see Theorem~\ref{PFisemi} in Section~\ref{Csemiplug}.
A closed form representation of the remainder for a finite sample setup allows us 
to explore the convergence and accuracy of many iterative EM-type algorithms in
a high-dimensional situation.

\paragraph{Sup-norm estimation} is studied as a special setup in
the semiparametric plug-in approach:
one component of the vector \( \prmtv \) is treated as the target while the rest is the nuisance. 
For a special case with a stochastically dominant information matrix \( \IF  = - \nabla^{2} \E L(\upsvs) \),
we prove a coordinate-wise expansion of the MLE \( \tilde{\upsv} \)
\begin{EQA}
	\bigl\| \DPN (\tilde{\upsv} - \upsvs) - \DPN^{-1} \nabla \zeta) \bigr\|_{\infty}
	& \leq &
	\frac{\dltwunss}{1 - \crosssup} \, \| \DPN^{-1} \nabla \zeta \|_{\infty}^{2} + \frac{\crosssup}{1 - \crosssup} \, \| \DPN^{-1} \nabla \zeta \|_{\infty}
	\, ,
\label{usdhyw6hikhurnetrspp2ei}
\end{EQA}
where \( \DPN_{j}^{2} \eqdef \nabla_{\ups_{j} \ups_{j}} \E \LL(\upsvs) \),
\( \DPN^{2} \eqdef \diag\{ \DPN_{1}^{2}, \ldots, \DPN_{\dimp}^{2} \} \),
\( \nabla \zeta = \nabla \LL(\upsvs) \),
and \( \dltwunss,\crosssup \) are explicitly given small constants; see Section~\ref{Ssupsemi}.
In particular,
\begin{EQA}[c]
	\crosssup^{2}
	= 
	\max_{j = 1,\ldots,\dimp} \frac{1}{\DPN_{j}^{2}}
	\sum_{m \neq j} \frac{\IFT_{jm}^{2}}{\DPN_{m}^{2}} 
	\, .
\label{dujekfyxrrveivut3c}
\end{EQA}
The result \eqref{usdhyw6hikhurnetrspp2ei} means a uniform in \( j \in \{ 1, \ldots, \dimp \} \) approximation
\begin{EQA}[c]
	\DPN_{j} \bigl( \tilde{\ups}_{j} - \upss_{j} \bigr)
	\approx
	\frac{1}{\DPN_{j}} \nabla_{j} \zeta
	\, .
\end{EQA}
Moreover, the remainder is given in a closed form
and the established results allow to relax the critical dimension condition \( \neff \gg \dimp^{2} \)
from \cite{Sp2024} to \( \neff \gg \log(\dimp) \).
Section~\ref{ScoBTL} illustrates the obtained results on the case 
of a BTL model.
The derived first-order expansions for the MLE are more accurate than 
ones from \cite{GSZ2023} and established under milder model assumptions.
For the experimental setup from \cite{GSZ2023},
we provide an empirical evidence that \( \crosssup \) from \eqref{dujekfyxrrveivut3c} 
is small and the expansion in \eqref{usdhyw6hikhurnetrspp2ei} is very accurate.

\paragraph{Marginal perturbed optimization.}
The main technical tool of the study is the theory of \emph{perturbed optimization}. 
The results from \cite{Sp2024} explain how the solution and value of an optimization problem change 
after a linear, quadratic, or smooth perturbation of the objective function.
This paper focuses on the case of partial and marginal optimization with a low-dimensional target variable 
and a possible large-dimensional nuisance variable; see Section~\ref{Ssemiopt}. 
This setup requires to evaluate a special \emph{semiparametric bias} caused by using an inexact value of a nuisance variable.
It appears that the value of this bias strongly depends on the norm in which smoothness of the objective function 
w.r.t. the nuisance variable is measured. 
An important special case is given by the sup-norm; see Section~\ref{Ssupnormco}. 


\paragraph{Organization of the paper.}
General problem of semiparametric plug-in estimation is discussed in Section~\ref{Csemiplug}.
Section~\ref{Ssupsemi} discusses the estimation problem in sup-norm.
Section~\ref{Slocalsmooth} collects the results for perturbed optimization including the cases of partial and marginal optimization.
Section~\ref{ScoBTL} illustrates the obtained results for a BTL model. 
Appendix~\ref{SgenBounds} and Appendix~\ref{Squadnquad} overview the general SLS theory from \cite{Sp2024}.
The proofs are gathered in Appendix~\ref{Ssemiopttools}.

{
\renewcommand{\Section}[1]{\section{#1}}
\renewcommand{\Subsection}[1]{\subsection{#1}}

}

%

\input semi_plugin

{
\renewcommand{\Section}[1]{\section{#1}}
\input semi_sup

}

\Chapter{Marginal optimization and sup-norm estimation}
\label{Slocalsmooth}
This section discusses a general problem of a marginal optimization.
Then the results are applied to optimization in sup-norm.
All the proofs are gathered in Section~\ref{Ssemiopttools}.

\input semi_opt_short

\input BTL_short

\bibliography{exp_ts,listpubm-with-url}

\newpage
\appendix

\appendix
\input genMLE_short

\input localbounds_short

\Chapter{Technical proofs}
\label{Ssemiopttools}
This section collects the proofs of the general statements from Section~\ref{Slocalsmooth}.
\input semi_opt_tools


\end{document}

%% file: myfrontmatter.tex
\usepackage{ifthen}
\usepackage[ampersand]{easylist}
\usepackage{amsmath,amssymb,amsthm}
\usepackage[bb=stixtwo]{mathalpha}
\usepackage{mathtools}
\usepackage{natbib}
\usepackage{epsfig,graphicx}
\usepackage{comment}
\usepackage{color}
\usepackage{srcltx}
\usepackage[mathscr]{eucal}
\usepackage[math]{easyeqn}
\usepackage{etoolbox}
\usepackage{hyperref}
\hypersetup{
            colorlinks,
            linkcolor=hookersgreen,
            linktoc=blue,
            citecolor=ultramarine,
            urlcolor=black,
            filecolor=black
            }
\usepackage{nameref}
\usepackage{multirow}

\usepackage{accents}
\usepackage{mathrsfs}
\usepackage{bbm}
\usepackage{algorithm}
\usepackage{algpseudocode}

\ifims{
\textheight=22cm
\textwidth=14.8cm
\topmargin=0pt
\oddsidemargin=1.0cm
\evensidemargin=1.0cm
\linespread{1.3}

}{ 
}

\numberwithin{equation}{section}
\numberwithin{figure}{section}
\newcounter{example}[section]
\numberwithin{example}{section}
\newcounter{remark}[section]
\numberwithin{remark}{section}
\newtheorem{theorem}{Theorem}[section]
\newtheorem{proposition}[theorem]{Proposition}
\newtheorem{lemma}[theorem]{Lemma}

\newtheorem{exmp}[example]{Example}
\newtheorem{rmrk}[remark]{Remark}
\newenvironment{example}{\begin{exmp}\rm}{\end{exmp}}
\newenvironment{remark}{\begin{rmrk}\rm}{\end{rmrk}}

\ifbook{
    \newcommand{\Chapter}[1]{\section{#1}}
    \newcommand{\Section}[1]{\subsection{#1}}
    \newcommand{\Subsection}[1]{\subsubsection{#1}}
    \def\Chname{Section }
    \def\chname{section }
    
  }
  {
    \newcommand{\Chapter}[1]{\chapter{#1}}
    \newcommand{\Section}[1]{\section{#1}}
    \newcommand{\Subsection}[1]{\subsection{#1}}
    \def\Chname{Chapter}
    
  }

%% file: mydef.tex
\renewcommand{\(}{$\,}
\renewcommand{\)}{\,$}

\def\nquad{\hspace{-1cm}}
\def\eqdef{\stackrel{\operatorname{def}}{=}}

\DeclareMathAlphabet{\mathbbmsl}{U}{bbm}{bx}{sl}

\DeclareMathSymbol{\Alpha}{\mathalpha}{operators}{"41}
\DeclareMathSymbol{\Beta}{\mathalpha}{operators}{"42}
\DeclareMathSymbol{\Epsilon}{\mathalpha}{operators}{"45}
\DeclareMathSymbol{\Zeta}{\mathalpha}{operators}{"5A}
\DeclareMathSymbol{\Eta}{\mathalpha}{operators}{"48}
\DeclareMathSymbol{\Iota}{\mathalpha}{operators}{"49}
\DeclareMathSymbol{\Kappa}{\mathalpha}{operators}{"4B}
\DeclareMathSymbol{\Mu}{\mathalpha}{operators}{"4D}
\DeclareMathSymbol{\Nu}{\mathalpha}{operators}{"4E}
\DeclareMathSymbol{\Omicron}{\mathalpha}{operators}{"4F}
\DeclareMathSymbol{\Rho}{\mathalpha}{operators}{"50}
\DeclareMathSymbol{\Tau}{\mathalpha}{operators}{"54}
\DeclareMathSymbol{\Chi}{\mathalpha}{operators}{"58}
\DeclareMathSymbol{\omicron}{\mathord}{letters}{"6F}

\newcommand{\cc}[1]{\mathscr{#1}}
\newcommand{\bb}[1]{\boldsymbol{#1}}

\DeclareFontFamily{U}{mathx}{\hyphenchar\font45}
\DeclareFontShape{U}{mathx}{m}{n}{
<5><6><7><8><9><10>
<10.95><12><14.4><17.28><20.74><24.88>
mathx10
}{}
\DeclareSymbolFont{mathx}{U}{mathx}{m}{n}
\DeclareFontSubstitution{U}{mathx}{m}{n}
\DeclareMathAccent{\widebar}{0}{mathx}{"73}

\renewcommand{\bar}[1]{\widebar{#1}}
\renewcommand{\hat}[1]{\widehat{#1}}
\renewcommand{\tilde}[1]{\widetilde{#1}}

\makeatletter
\def\mathcenterto#1#2{\mathclap{\phantom{#1}\mathclap{#2}}\phantom{#1}}
\let\old@widetilde\widetilde
\def\widetildeto#1#2{\mathcenterto{#2}{\old@widetilde{\mathcenterto{#1}{#2\,}}}}
\let\old@widehat\widehat
\def\widehatto#1#2{\mathcenterto{#2}{\old@widehat{\mathcenterto{#1}{#2\,}}}}
\makeatother

\newcommand{\thankstitle}[1]{\ifthenelse{\equal{#1}{}}{}{\thanks{#1}}}
\newcommand{\thanksau}[1]{\ifthenelse{\equal{#1}{}}{}{\thanks{#1}}}

\newcommand{\aua}[6]
{\def\authora{#1}
\def\runauthora{#2}
\def\addressa{#3}
\def\emaila{#4}
\def\affiliationa{#5}
\def\thanksa{#6}}

\def\theauthors{
\ifau{ 
  \author{
    \authora
    \thanksau{\thanksa}
    \\[5.pt]
    \addressa \\
    \texttt{ \emaila}
  }
}
{  
  \author{
    \authora
    \thanksau{\thanksa}
    \\[5.pt]
    \addressa \\
    \texttt{ \emaila}
    \and
    \authorb
    \thanksau{\thanksb}
    \\[5.pt]
    \addressb \\
    \texttt{ \emailb}
  }
}
{   
  \author{
    \authora
    \thanksau{\thanksa}
    \\[5.pt]
    \addressa \\
    \texttt{ \emaila}
    \and
    \authorb
    \thanksau{\thanksb}
    \\[5.pt]
    \addressb \\
    \texttt{ \emailb}
    \and
    \authorc
    \thanksau{\thanksc}
    \\[5.pt]
    \addressc \\
    \texttt{ \emailc}
  }
} {   
  \author{
    \authora
    \thanksau{\thanksa}
    \\[5.pt]
    \addressa \\
    \texttt{ \emaila}
    \and
    \authorb
    \thanksau{\thanksb}
    \\[5.pt]
    \addressb \\
    \texttt{ \emailb}
    \and
    \authorc
    \thanksau{\thanksc}
    \\[5.pt]
    \addressc \\
    \texttt{ \emailc}
    \and
    \authord
    \thanksau{\thanksd}
    \\[5.pt]
    \addressd \\
    \texttt{ \emaild}
  }
}
}

\renewcommand{\Gamma}{\varGamma}
\renewcommand{\Pi}{\varPi}
\renewcommand{\Sigma}{\varSigma}
\renewcommand{\Delta}{\varDelta}
\renewcommand{\Lambda}{\varLambda}
\renewcommand{\Psi}{\varPsi}
\renewcommand{\Phi}{\varPhi}
\renewcommand{\Theta}{\varTheta}
\renewcommand{\Omega}{\varOmega}
\renewcommand{\Xi}{\varXi}
\renewcommand{\Upsilon}{\varUpsilon}

\def\argmax{\operatornamewithlimits{argmax}}
\def\argmin{\operatornamewithlimits{argmin}}

\def\av{\bb{a}}

\def\ev{\bb{e}}

\def\uv{\bb{u}}

\def\wv{\bb{w}}

\def\zv{\bb{z}}

\def\Av{\bb{A}}

\def\Uv{\bb{U}}

\def\epsv{\bb{\varepsilon}}

\newenvironment{myitem}{
  \setlength{\parskip}{0pt}
\begin{itemize}
  \setlength{\itemsep}{0pt}
  \setlength{\parskip}{0pt}
  \setlength{\parsep}{0pt}
}{\end{itemize}}

\usepackage{color}

\definecolor{blue(pigment)}{rgb}{0.2, 0.2, 0.6}
\definecolor{ultramarine}{rgb}{0.07, 0.04, 0.56}
\definecolor{darkspringgreen}{rgb}{0.09, 0.45, 0.27}
\definecolor{hookersgreen}{rgb}{0.0, 0.44, 0.0}
\definecolor{hgreen}{rgb}{0.0, 0.44, 0.0}
\definecolor{plum(traditional)}{rgb}{0.56, 0.27, 0.52}
\definecolor{purple(html/css)}{rgb}{0.5, 0.0, 0.5}
\definecolor{magenta(dye)}{rgb}{0.79, 0.08, 0.48}

%% file: statdef.tex
\newcommand{\binomb}[2]{\genfrac{}{}{0pt}{}{#1}{#2}}
\newcommand{\scorem}[1]{\score_{\!#1}}

\def\neff{\mathbbmsl{N}}
\def\hspm{\hspace{1pt}}

\def\AFN{\mathbbmsl{U}}
\def\Avm{\bb{M}}

\def\inftyi{\scalebox{0.66}{$\infty$}}
\def\inftyii{\scalebox{0.44}{$\infty$}}
\def\duali{1}

\def\crossnorm{\rho_{\dual}}

\def\crossgrad{\rho_{2}}
\def\crosssup{\rho_{\duali}}

\def\cmax{\const}

\def\AFN{\mathbb{Z}}

\def\dmax{\kappa}

\def\AAv{{\mathcal{A}}}
\def\CCv{{\mathcal{C}}}

\def\Range{\mathcal{R}}

\def\IFTN{\IFT_{0}}

\def\AFN{\mathbb{U}}
\def\HFN{\mathbb{H}}

\def\DFNb{\breve{\DFN}}

\def\prmt{\ups}
\def\prmts{\prmt^{*}}
\def\prmtv{\bb{\prmt}}

\def\prmtvs{\prmtv^{*}}

\def\targ{x}
\def\targv{\bb{\targ}}
\def\targvs{\targv^{*}}
\def\targvb{\breve{\targv}}

\def\targvn{\targv^{\circ}}

\def\tarp{\theta}
\def\tarpv{\bb{\tarp}}

\def\tarpvs{\tarpv^{*}}

\def\CFT{\cc{C}}

\def\hmax{\mathsf{c}}

\def\hL{h}

\def\Matr{\mathfrak{M}}
\def\Graph{\mathcal{G}}
\def\GrV{\mathcal{V}}
\def\GrE{\mathcal{E}}

\def\nano{\circ}
\def\dual{*}

\def\feta{\phi}

\def\Eta{\mathcal{H}}

\def\nbin{N}

\def\HL{\mathbb{m}}

\def\dltwb{\omega}

\def\dltwu{\tau}

\def\dltwbss{\dltwb^{+}}

\def\dltwun{\dltwu_{12}}
\def\dltwunn{\dltwu_{21}}
\def\dltwuns{\dltwu_{\nano}}
\def\dltwunss{\dltwu_{\inftyi}}

\def\dltwbun{\dltwb}

\def\dblk{D}

\def\DFN{\DVL}

\def\DFblk{\dblk}

\def\R{\mathbbmsl{R}}
\def\E{\mathbbmsl{E}}

\def\P{\mathbbmsl{P}}

\def\kappa{\varkappa}

\def\blk{\operatorname{block}}

\def\diag{\operatorname{diag}}

\def\Var{\operatorname{Var}}

\def\T{\top}
\def\tr{\operatorname{tr}}

\def\ex{\mathrm{e}}

\def\Id{\mathbbmsl{I}}
\def\Ind{\operatorname{1}\hspace{-4.3pt}\operatorname{I}}

\def\alp{\alpha}

\def\BBH{B}
\def\BBN{B}

\def\cdens{\phi}

\def\CONST{\mathtt{C} \hspace{0.1em}}
\def\CONSTi{\mathtt{C}}

\def\DPN{D}

\def\DVL{\mathbb{D}}

\def\dimA{\mathbb{p}}

\def\dimttl{\bar{\dimp}}

\def\dimp{p}

\def\dimG{\dimA_{\GP}}

\def\dimQ{\dimA_{\QP}}

\def\dimq{q}

\def\dimD{\dimA_{\scalebox{0.666}{${\DPN}$}}}

\def\Eta{\cc{H}}

\def\fs{f}
\def\fsb{\breve{\fs}}

\def\fn{g}

\def\gp{g}

\def\GP{G}

\def\IF{\mathbbmsl{F}}

\def\IFN{\IFL}

\def\IFL{{\mathbbmss{F}}}

\def\IFtotal{F}
\def\IFT{\mathscr{\IFtotal}}

\def\IFTb{\Phi}

\def\IFL{\mathbb{F}}

\def\Kappa{\cc{K}}

\def\LT{L}
\def\LGP{\LT_{\GP}}

\def\LL{\cc{L}}

\def\nuov{\bb{a}}
\def\nuovs{\nuov^{*}}

\def\nupv{\bb{\nup}}
\def\nupvs{\nupv^{*}}

\def\nui{s}
\def\Nui{\mathcal{S}}

\def\nuiv{\bb{\nui}}

\def\nuivs{\nuiv^{*}}

\def\nuo{\tau}
\def\nuov{\bb{\nuo}}

\def\QP{Q}

\def\rhoIF{\rho}

\def\riskt{\cc{R}}

\def\rr{\mathtt{r}}

\def\rru{\rr_{\circ}}

\def\rrB{\rr_{\scalebox{0.666}{${\BBH}$}}}
\def\rrD{\rr_{\scalebox{0.666}{${\DPN}$}}}

\def\rrinf{\rr_{\inftyi}}
\def\rrinfi{\rr_{\inftyii}}

\def\rrn{\rr}

\def\score{\nabla}

\def\Tau{T}

\def\ups{\upsilon}
\def\upsv{\bb{\ups}}

\def\upsvd{\upsv^{\circ}}
\def\upsvs{\upsv^{*}}

\def\upsvn{\upsvd}

\def\ups{\upsilon}
\def\upsv{\bb{\ups}}

\def\upsd{\ups^{\circ}}

\def\upss{\ups^{*}}

\def\UVw{\mathcal{W}}

\def\Ups{\varUpsilon}
\def\Upsd{\Ups^{\circ}}

\def\VP{V}

\def\wv{\bb{w}}

\def\WV{\mathcal{W}}

\def\xx{\mathtt{x}}

\def\zq{z}

%% file: semi_plugin.tex

\def\nupv{\nuiv}
\def\nquadr{\hspace{-12pt}}
\def\dltwbss{\dltwb}
\def\HPN{H}

\Chapter{Semiparametric plug-in estimation}
\label{Csemiplug}
This section explains an approach to semiparametric estimation 
based on partial optimization which is not limited to the profile MLE.
It assumes that a reasonable pilot estimator of the nuisance parameter is available.
The aim is to establish some accurate finite sample theoretical bound on the estimation accuracy 
which correspond to the dimension of the target parameter only.
The idea can be explained as follows.
Define a full dimensional truth \( \prmtvs = (\tarpvs,\nupvs) \) by
\begin{EQA}
	\prmtvs 
	&=& 
	(\tarpvs,\nupvs)
	=
	\argmax_{(\tarpv,\nupv)} \E \LL(\tarpv,\nupv) .
\label{ytxtysys8w725ew62}
\end{EQA}
Further, consider a family of partial optimization problems with respect to 
the target parameter \( \tarpv \) for each fixed value of the nuisance parameter \( \nupv \):
\begin{EQ}[rcccl]
	\tilde{\tarpv}(\nupv)
	& \eqdef &
	\argmax_{\tarpv} \LL(\tarpv,\nupv) \, ,
\label{ttn2ttELGettint}
	\qquad
	\tarpvs(\nupv)
	& \eqdef &
	\argmax_{\tarpv} \E \LL(\tarpv,\nupv) \, .
\label{ttn2ttELGettintpo}
\end{EQ}
We also assume that each partial problem can be easily solved.
Further, let a pilot estimate \( \hat{\nupv} \) of the nuisance parameter \( \nupv \) be given.
This leads to the \emph{plug-in estimator} 
\(
	\hat{\tarpv}
	=
	\tilde{\tarpv}(\hat{\nupv}) 
	.
\)
The profile MLE is a special case of the plug-in method with \( \hat{\nupv} = \tilde{\nupv} \), where 
\( \tilde{\prmtv} = (\tilde{\tarpv},\tilde{\nupv}) \) is a full dimenstional maximizer of \( \LL(\tarpv,\nupv) \).
The results for the profile MLE can be derived from the general results on estimation of the full-dimensional 
parameter \( \prmtv \); see \cite{Sp2024}.
Later we do not limit ourselves to any particular choice of \( \hat{\nupv} \). 
This requires to develop different tools and approaches.
The aim is to show that under reasonable conditions on the model and the pilot \( \hat{\nupv} \),
this plug-in estimator behaves nearly as the ``oracle'' estimator \( \tilde{\tarpv}(\nupvs) \).
A very important requirement for this construction is that the pilot estimator \( \hat{\nupv} \) 
concentrates on a small vicinity \( \Eta_{0} \) of the point \( \nupvs \).
Further, the results of \cite{Sp2024} briefly sketched in Section~\ref{SgenBounds} imply a bound on the difference 
\( \tilde{\tarpv}(\nupv) - \tarpvs(\nupv) \) between
the partial estimates \( \tilde{\tarpv}(\nupv) \) and its population counterpart 
\( \tarpvs(\nupv) \) in terms of the effective target dimension.
Moreover, the full dimensional condition \nameref{Eref} enables us to state such bounds 
uniformly over the concentration set \( \Eta_{0} \) of the pilot \( \hat{\nupv} \).
%
Under a so called ``small bias'' condition \( \tarpvs(\nupv) \approx \tarpvs \)
and \( \IF(\nupv) \approx \IF(\nupvs) \),
behavior of \( \tilde{\tarpv}(\nupv) \) only weakly depends on \( \nupv \in \Eta_{0} \)
yielding the desirable properties of \( \hat{\tarpv} \).
An extension to a general situation requires a careful evaluation of the so-called ``semiparametric bias''
caused by using the pilot value \( \hat{\nupv} \) in place of the true value \( \nupvs \).

\Section{Uniform bounds in partial estimation}
This section establishes some bounds for \emph{partial estimation}  of the target parameter \( \tarpv \) that apply uniformly over 
\( \nupv \in \Eta_{0} \) for  some set \( \Eta_{0} \in \Eta \).
Usually this is a concentration set of a pilot estimator \( \hat{\nupv} \). 
The analysis follows \cite{Sp2024}.
In particular, we assume the full-dimensional general condition \nameref{Eref} of linearity in \( \upsv \) of the stochastic component
\( \zeta(\upsv) = \LL(\upsv) - \E \LL(\upsv) \);
see Section~\ref{SgenBounds}.
\cite{Sp2024} explains how calming trick and localization can be used to enforce this condition 
for rather complex models like deep neuronal networks.

For each \( \nupv \in \Eta_{0} \), consider \( \LL(\tarpv,\nupv) \) 
as a function of \( \tarpv \).
We suppose that the function \( \fs_{\nupv}(\tarpv) = \E \LL(\tarpv,\nupv) \) is concave and
satisfies condition \nameref{LLsoT3ref} meaning smoothness of \( \fs_{\nupv}(\tarpv) \) in \( \tarpv \)
for \( \nupv \) fixed  
but with constant \( \dltwu_{3} \), the metric tensor \( \DPN \), and the radius \( \rr \)
which do not depend on \( \nupv \in \Eta_{0} \).
The matrix 
\begin{EQA}
	\IF(\nupv) 
	&=& 
	- \nabla_{\tarpv\tarpv}^{2} \E \LL(\tarpvs(\nupv),\nupv) 
\label{tdiucxu8cucudfuyrdfyh}
\end{EQA}
describes the information about the target parameter \( \tarpv \) in the partial model with \( \nupv \) fixed.
Later we need that 
the matrix \( \IF(\nupv) \) does not vary too much within \( \Eta_{0} \), that is, 
\( \IF(\nupv) \approx \IF = \IF(\nupvs) \).
More precisely, we assume that for some \( \dltwbss \leq 1/3 \) 
\begin{EQA}
	- \dltwbss \DPN^{2} 
	& \leq &
	\IF(\nupv) - \IF
	\leq 
	\dltwbss \DPN^{2} \, ,
	\quad
	\forall \nupv \in \Eta_{0} \, .
\label{kjhgtrewsxsdefxzsedr}
\end{EQA}
This can be effectively checked under condition \nameref{LLpT3ref}; see Lemma~\ref{L3IFNNui} and Lemma~\ref{LIFvari}. 
Due to assumption \nameref{Eref} about linearity of the stochastic term \( \zeta(\prmtv) \),
the \( \tarpv \)-gradient \( \scorem{\tarpv} \zeta \) does not depend on \( \nupv \).
This substantially simplifies our study.
We also assume that \( \scorem{\tarpv} \zeta \) satisfies \nameref{EU2ref}.
Define \( \BBH_{\DPN} = \Var( \DPN \, \IF^{-1} \scorem{\tarpv} \zeta) \), \( \dimD = \tr \BBH_{\DPN} \), and
\begin{EQA}
	\rr_{\DPN}
	& \eqdef &
	\frac{1}{(1 - \dltwbss)^{1/2}} \,\, \zq(\BBH_{\DPN},\xx)
	=
	\frac{1}{(1 - \dltwbss)^{1/2}} \,\, \bigl( \sqrt{\tr \BBH_{\DPN}} + \sqrt{2\xx \, \| \BBH_{\DPN} \|} \bigr) .
	\qquad
\label{2emxGPm12nz122s}
\end{EQA}
It is important that the radius \( \rr_{\DPN} \) corresponds to the dimension of the target parameter \( \tarpv \).
Now we apply the general results from Section~\ref{SgenBounds} to the partial pMLEs \( \tilde{\tarpv}(\nupv) \).

\begin{proposition}
\label{PMLEsemipart}
Assume \nameref{Eref}, and let the \( \tarpv \)-component\( \nabla_{\tarpv} \zeta \) of \( \nabla \zeta \) satisfy \nameref{EU2ref}.
Let also \eqref{kjhgtrewsxsdefxzsedr} hold with \( \dltwbss \leq 1/3 \).
For any \( \nupv \in \Eta_{0} \), let \( \fs_{\nupv}(\tarpv) = \E \LL(\tarpv,\nupv) \) 
satisfy \nameref{LLsoT3ref} at \( \tarpvs(\nupv) \) with \( \rr_{\nupv} \equiv \rr \), \( \DPN_{\nupv} \equiv \DPN \), 
and \( \dltwu_3 \) such that 
\begin{EQA}[c]
	\DPN^{2} \leq \dmax^{2} \, \IF(\nupv) ,
	\quad
	\rr \geq \frac{3}{2} \, \rr_{\DPN} \, ,
	\quad
	\dltwu_{3} \, \dmax^{2} \, \rr_{\DPN} < \frac{4}{9} \, .
\label{8difiyfc54wrboesep}
\end{EQA}
Then on a random set \( \Omega(\xx) \) with \( \P(\Omega(\xx)) \geq 1 - 3 \ex^{-\xx} \), 
it holds for all \( \nupv \in \Eta_{0} \)
\begin{EQA}[rcl]
	&& \nquad
	\bigl\| \QP \bigl\{ \tilde{\tarpv}(\nupv) - \tarpvs(\nupv)
		- \IF^{-1}(\nupv) \, \scorem{\tarpv} \zeta \bigr\}
	\bigr\|
	\leq 
	\| \QP \IF^{-1}(\nupv) \DPN \| \, \frac{3 \dltwu_{3}}{4} \, 
	\| \DPN \, \IF^{-1}(\nupv) \, \scorem{\tarpv} \zeta \|^{2} 
	\\
	& \leq &
	\| \QP \IF^{-1} \DPN \| \, \frac{(3/4) \dltwu_{3}}{(1 - \dltwbss)^{3/2}} \, 
	\| \DPN \, \IF^{-1} \, \scorem{\tarpv} \zeta \|^{2} 
	\, .
	\qquad
\label{3a3G1ma3G2Dm1e}
\end{EQA}
Also, it holds on \( \Omega(\xx) \)
\begin{EQA}
	\| \QP \{ \IF^{-1}(\nupv) \, \scorem{\tarpv} \zeta - \IF^{-1} \, \scorem{\tarpv} \zeta \} \|
	& \leq &
	\| \QP \IF^{-1} \DPN \| \, \frac{\dltwbss}{1 - \dltwbss} \, 
	\| \DPN \, \IF^{-1} \scorem{\tarpv} \zeta \| \, .
\label{vivkj7vjefv98vkehyvui}
\end{EQA}
\end{proposition}

\begin{proof}
Condition \nameref{EU2ref} implies on \( \Omega(\xx) \) 
\begin{EQA}
	\| \DPN \, \IF^{-1} \, \scorem{\tarpv} \zeta \|
	& \leq &
	\zq(\BBH_{\DPN},\xx) \, .
\label{c}
\end{EQA} 
This, \eqref{2emxGPm12nz122s}, and \eqref{kjhgtrewsxsdefxzsedr} yield on \( \Omega(\xx) \) for any \( \nupv \in \Eta_{0} \)
\begin{EQA}
	\bigl\| \DPN \, \IF^{-1}(\nupv) \, \scorem{\tarpv} \zeta \bigr\|
	& \leq &
	\rr_{\DPN} \, .
\label{1d3Gm12rGDGm1e}
\end{EQA}
The first bound in \eqref{3a3G1ma3G2Dm1e}  follows from Theorem~\ref{TFiWititG3} applied to
\( \tilde{\tarpv}(\nupv) = \argmax_{\tarpv} \{ \LL(\tarpv,\nupv) - \| \GP \tarpv \|^{2}/2 \} \). 
Further, for each \( \nupv \in \Eta_{0} \), by \eqref{kjhgtrewsxsdefxzsedr}
\begin{EQA}
	\bigl\| 
		\DPN^{-1} \IF \, \{ \IF^{-1}(\nupv) - \IF^{-1} \} \IF(\nuiv) \DPN^{-1}
	\bigr\| 
	& = &
	\bigl\| \DPN^{-1} \{ \IF - \IF(\nupv) \} \DPN^{-1} \bigr\| 
	\leq 
	\dltwbss \,  
\label{8ju8j8okyt443w3e34e3d}
\end{EQA}
and
\begin{EQA}
	&& \nquad
	\| \QP \{ \IF^{-1}(\nupv) \, \scorem{\tarpv} \zeta - \IF^{-1} \, \scorem{\tarpv} \zeta \} \|
	\leq 
	\| \QP \IF^{-1} \DPN \| \,\,
	\| \DPN^{-1} \IF \{ \IF^{-1} - \IF^{-1}(\nupv) \} \scorem{\tarpv} \zeta \|
	\\
	& \leq &
	\| \QP \IF^{-1} \DPN \| \, \,
	\dltwbss \,\, \| \DPN \, \IF^{-1}(\nuiv) \scorem{\tarpv} \zeta \|
\label{ydjcuvuggugh36hfjrfttyy}
\end{EQA}
With \( \QP = \DPN^{-1} \IF \), this implies
\begin{EQA}
	(1 - \dltwbss) \| \DPN \, \IF^{-1} \scorem{\tarpv} \zeta \|
	\leq 
	\| \DPN \, \IF^{-1}(\nuiv) \scorem{\tarpv} \zeta \|
	& \leq &
	(1 + \dltwbss) \| \DPN \, \IF^{-1} \scorem{\tarpv} \zeta \|
\label{gduheidhiewhdi6erd6edu}
\end{EQA}
and \eqref{vivkj7vjefv98vkehyvui} follows as well as the second bound in \eqref{3a3G1ma3G2Dm1e}.
\end{proof}

A benefit of \eqref{3a3G1ma3G2Dm1e} is that the accuracy of estimation corresponds to the dimension 
of the target component only.
Another benefit of \eqref{3a3G1ma3G2Dm1e} and \eqref{vivkj7vjefv98vkehyvui} is that these bounds hold on \( \Omega(\xx) \)
for all \( \nupv \in \Eta_{0} \) simultaneously. 
This is granted by \nameref{Eref}, \nameref{EU2ref}, and \eqref{kjhgtrewsxsdefxzsedr}.
Note, however, that fixing the nuisance parameter \( \nupv \) 
changes the value \( \tarpvs = \tarpvs(\nupvs) \) to \( \tarpvs(\nupv) \), 
and we have to control the variability of this estimate w.r.t. \( \nupv \in \Eta_{0} \).

\Section{Semiparametric bias under (semi)orthogonality}
\label{SOporthogo}

The results of Proposition~\ref{PMLEsemipart} rely on variability of  
\( \tarpvs(\nupv) = \argmax_{\tarpv} \fs(\tarpv,\nupv) \) w.r.t. the nuisance parameter \( \nupv \),
where \( \fs(\prmtv) = \E \LL(\prmtv) = \E \LL(\tarpv,\nupv) \).
This section studies the \emph{semiparametric bias} \( \tarpvs(\nupv) - \tarpvs \).
It appears that local quadratic approximation of the function \( \fs \) in a vicinity of 
\( \prmtvs \) yields a nearly linear dependence of \( \tarpvs(\nupv) \) in \( \nupv \).
To see this, represent \( \IFT(\prmtv) = - \nabla^{2} \fs(\prmtv) \) in the block form 
\begin{EQA}
	\IFT(\prmtv)
	&=& 
	\begin{pmatrix}
		\IFT_{\tarpv\tarpv}(\prmtv) & \IFT_{\tarpv\nupv}(\prmtv) \\
		\IFT_{\nupv\tarpv}(\prmtv) & \IFT_{\nupv\nupv}(\prmtv)
	\end{pmatrix} .
\label{HAuuuvvuvvTHo}
\end{EQA}
We write \( \IF(\prmtv) = \IFT_{\tarpv\tarpv}(\prmtv) \) and
\( \IFT = \IFT(\prmtvs) \).
Consider first a \emph{quadratic} \( \fs \).
As \( \nabla \fs(\prmtvs) = 0 \), it holds 
\( \fs(\prmtv) = \fs(\prmtvs) - (\prmtv - \prmtvs)^{\T} \IFT \, (\prmtv - \prmtvs)/2 \). 
For \( \nupv \) fixed, the point \( \tarpvs(\nupv) \) satisfies 
\begin{EQA}
	\tarpvs(\nupv)  
	&=& 
	\argmax_{\tarpv} \fs(\tarpv,\nupv)
	=
	\argmin_{\tarpv} \{ 
		(\tarpv - \tarpvs)^{\T} \IF \, (\tarpv - \tarpvs) /2 
		+ (\tarpv - \tarpvs)^{\T} \IFT_{\tarpv\nupv} \, (\nupv - \nupvs) 
	\}
	\\
	&=&
	\tarpvs - \IF^{-1} \, \IFT_{\tarpv\nupv} \, (\nupv - \nupvs) .
\label{9nvyf6eyndfvthwehgvhe}
\end{EQA}
This observation is in fact discouraging because the bias \( \tarpvs(\nupv) - \tarpvs \) has the same (in order) magnitude as 
the nuisance parameter \( \nupv - \nupvs \).
%
However, the \emph{orthogonality} condition 
\begin{EQA}[c]
	\scorem{\tarpv} \nabla_{\nupv} \, \fs(\tarpv,\nupv) \equiv 0 , \qquad
	\forall (\tarpv,\nupv) \in \UVw \, ,
\label{7edjcgt45rtyecjmf76w3}
\end{EQA}
still ensures a vanishing bias, even if \( \fs(\prmtv) \) is not quadratic.
Indeed, it implies the decomposition \( \fs(\tarpv,\nupv) = \fs_{1}(\tarpv) + \fs_{2}(\nupv) \) for some functions 
\( \fs_{1} \) and \( \fs_{2} \).
As a corollary, the maximizer \( \tarpvs(\nupv) \) and the corresponding negative Hessian \( \IF(\nupv) \)
do not depend on \( \nupv \) yielding 
\( \tarpvs(\nupv) \equiv \tarpvs \) and \( \IF(\nupv) \equiv \IF \).
This is a very useful property allowing to obtain accurate results about estimation accuracy of the target parameter \( \tarpv \) 
as if the true value of the nuisance parameter \( \nupvs \) were known.
In practice one may apply the plug-in estimator \( \tilde{\tarpv}(\hat{\nupv}) \),
where \( \hat{\nupv} \) is any reasonable estimate of \( \nupvs \).

Unfortunately, the orthogonality condition \( \scorem{\tarpv} \nabla_{\nupv} \, \fs(\tarpv,\nupv) \equiv 0 \) is too restrictive
and fulfilled only in some special cases.
One of them corresponds to the already mentioned additive case \( \fs(\tarpv,\nupv) = \fs_{1}(\tarpv) + \fs_{2}(\nupv) \).
If \( \fs(\tarpv,\nupv) \) is quadratic, then orthogonality can be achieved by a linear transform of the nuisance parameter \( \nupv \).
For a general function \( \fs \), such a linear transform helps to only ensure 
the \emph{one-point orthogonality} condition \( \scorem{\tarpv} \nabla_{\nupv} \, \fs(\prmtvs) = 0 \).
In some situation, \( \scorem{\tarpv} \nabla_{\nupv} \, \fs(\prmtvs) = 0 \) implies
\( \scorem{\tarpv} \nabla_{\nupv} \, \fs(\tarpvs,\nupv) = 0 \) for all \( \nupv \in \Eta_{0} \).
We refer to this situation as \emph{semi-orthogonality}.
A typical example is given by models for which the cross-derivative 
\( \scorem{\tarpv} \nabla_{\nupv} \, \fs(\tarpvs,\nupv) \)
depends on \( \tarpv \) only%
\ifapp{; see e.g. the case of nonlinear regression in Section~\ref{Scalming}.}{.}
In this situation,  the semiparametric bias vanishes.
More precisely, Lemma~\ref{Lsemiorto} yields 

\begin{proposition}
\label{Psemiortose}
Let \( \scorem{\tarpv} \nabla_{\nupv} \, \fs(\tarpvs,\nupv) = 0 \) for all \( \nupv \in \Eta_{0} \)
and \( \fs(\tarpv,\nupv) \) is concave in \( \tarpv \) for each \( \nupv \).
Then 
\begin{EQA}
	\tarpvs(\nupv)
	& \equiv &
	\tarpvs \, ,
	\qquad
	\IF(\nupv)
	\eqdef
	- \nabla^{2}_{\tarpv\tarpv} \fs(\tarpvs,\nupv)
	\equiv
	\IF \, .
\label{ocuerjecjhrwekjgv9ei}
\end{EQA}
\end{proposition}

As in the orthogonal case, the condition of semi-orthogonality
 allows to ignore the semiparametric bias; see Theorem~\ref{PWilkssemiso} later.


%

\Section{Loss and risk of a plug-in estimator}
This section discusses the properties of the plug-in estimator 
\( \hat{\tarpv} = \tilde{\tarpv}(\hat{\nupv}) \).
The major condition on \( \hat{\nupv} \) is that it belongs 
with high probability to the local set \( \Eta_{0} \).
To simplify our notation, we fix a random set \( \Omega(\xx) \) with 
\( \P(\Omega(\xx)) \geq 1 - 3\ex^{-\xx} \) and assume that \( \hat{\nupv} \in \Eta_{0} \) and condition \eqref{PMLEsemipart} of Proposition~\ref{PMLEsemipart} hold true on this set.

\Subsection{Orthogonal case}
In some special cases like orthogonality or semi-orthogonality considered in Proposition~\ref{Psemiortose},
it holds \( \tarpvs(\nupv) \equiv \tarpvs \),
\( \IF \equiv \IF(\nupv) \).
In particular, this condition meets for the special case when 
\( \IFT_{\tarpv\nupv}(\prmtv) \) 
and \( \IFT_{\nupv\nupv}(\prmtv) \) in \eqref{HAuuuvvuvvTHo} for \( \prmtv = (\tarpv,\nupv) \) depend on \( \tarpv \) only;
see Proposition~\ref{Psemiortose}.
Application of Theorem~\ref{TFiWititG3} to the squared risk of \( \hat{\tarpv} \)
in the partial model with \( \nupv = \nupvs \) yields similarly to Proposition~\ref{PMLEsemipart}
the following very strong result.

\begin{theorem}
\label{PWilkssemiso}
Let \( \| \HPN (\hat{\nupv} - \nupvs) \|_{\nano} \leq \rru \) on \( \Omega(\xx) \), and let \( \tarpvs(\nupv) \equiv \tarpvs \), \( \IF(\nupv) \equiv \IF \).
Under the conditions of Proposition~\ref{PMLEsemipart},
the plug-in estimator \( \hat{\tarpv} = \tilde{\tarpv}(\hat{\nupv}) \) satisfies on \( \Omega(\xx) \)
\begin{EQA}
	\bigl\| \DPN^{-1} \IF \bigl( \hat{\tarpv} - \tarpvs - \IF^{-1} \scorem{\tarpv} \zeta \bigr) \bigr\|
	& \leq &
	\frac{3\dltwu_{3}}{4} \, \| \DPN \, \IF^{-1} \scorem{\tarpv} \zeta \|^{2} \, .
\label{2biGe233GDttGseFso}
\end{EQA}
Moreover, for any linear mapping \( \QP \), define \( \dimQ = \E \, \| \QP \IF^{-1} \scorem{\tarpv} \zeta \|^{2} \) and
\begin{EQA}
	\riskt_{\QP}
	& \eqdef &
	\E \bigl\{ \| \QP \IF^{-1} \scorem{\tarpv} \zeta \|^{2} \Ind_{\Omega(\xx)} \bigr\}
	\leq 
	\dimQ \, .
\label{6edsjcu7d5ygvvjg3wtgse}
\end{EQA} 
With \( \dimD \) and \( \rrD \) from \eqref{2emxGPm12nz122s}, suppose
\begin{EQ}[rcl]
	\alp_{\QP}
	& \eqdef & 
	\frac{\| \QP \, \IF^{-1} \DPN \| (3/4)\dltwu_{3} \, \rrD \sqrt{\dimD}} {\sqrt{\riskt_{\QP}}}
	< 1
	\, . 
\label{c56hgwejufutee3gt5cse}
\end{EQ}
Then
\begin{EQA}
	(1 - \alp_{\QP})^{2} \riskt_{\QP} 
	\leq 
	\E \bigl\{ \| \QP \, (\hat{\tarpv} - \tarpvs) \|^{2} \Ind_{\Omega(\xx)} \bigr\}
	& \leq &
	(1 + \alp_{\QP})^{2} \riskt_{\QP} \, .
\label{EQtuGmstrVEQtGQSEP}
\end{EQA}
\end{theorem}

The standardized score \( \DFblk^{-1} \scorem{\tarpv} \zeta \) does not depend on \( \hat{\nupv} \).
Thus, for the statement of Theorem~\ref{PWilkssemiso}, it suffices
that \( \hat{\nupv} \) concentrates on \( \Eta_{0} \), all partial smoothness conditions 
on \( \fs(\tarpv,\nupv) = \E \LL(\tarpv,\nupv) \) w.r.t. \( \tarpv \) hold uniformly over \( \nupv \in \Eta_{0} \),
and \( \tarpvs(\nupv) \equiv \tarpvs \), \( \IF(\nupv) \equiv \IF \).

\Subsection{Concentration and Fisher expansion for the plug-in estimator}
This section presents our main results for the semiparametric plug-in estimator 
\( \hat{\tarpv} = \tilde{\tarpv}(\hat{\nupv}) \). 
The main issue in the analysis is the semiparametric bias \( \tarpvs(\nupv) - \tarpvs \) caused by using a pilot \( \hat{\nupv} \) in place of the truth \( \nupvs \).
This bias can be bounded under smoothness properties of \( \fs(\prmtv) \) 
using an accurate expansion of the bias term \( \tarpvs(\nupv) - \tarpvs \)
from Section~\ref{Sonepoint} for all \( \nupv \) from a vicinity \( \Eta_{0} \) of \( \nupvs \).
The approach allows to incorporate the anisotropic case when conditions on smoothness of 
\( \fs(\prmtv) = \fs(\tarpv,\nupv) \) are stated in different norm 
for the target parameter \( \tarpv \) and the nuisance parameter \( \nupv \).
A typical example is given by using a \( \ell_{2} \)-norm for \( \tarpv \) and
a sup-norm for \( \nupv \); see Section~\ref{Ssupnormco}.
Given a norm \( \| \cdot \|_{\nano} \) in \( \R^{\dimq} \), 
a metric tensor \( \HPN(\nupv) \) for \( \nupv \in \R^{\dimq} \), and a radius \( \rru \), consider local sets 
\( \Eta_{0} \) of the form
\begin{EQA}
	\Eta_{0}
	&=&
	\{ \nupv \colon \| \HPN (\nupv - \nupvs) \|_{\nano} \leq \rru \} 
	\, .
\label{8xjxcwbd6dejdccyw3du}
\end{EQA}
We assume conditions \nameref{LLsoT3ref}, \nameref{LLpT3ref} from Section~\ref{Ssemiopt}.
The first one requires that \( \fs(\tarpv,\nupv) \) is smooth in \( \tarpv \) for \( \nupv \) fixed, while
the second one describes the smoothness properties of \( \fs(\tarpvs,\nupv) \) w.r.t. \( \nupv \).
These conditions are weaker than the full dimensional smoothness condition \nameref{LLsT3ref}
and only involve partial derivatives in \( \tarpv \) and cross-derivatives of \( \fs(\prmtv) = \fs(\tarpv,\nupv) \);
see Section~\ref{Sonepoint} for a detailed discussion.
Moreover, in some sense, these conditions are nothing but definitions of the important quantities 
\( \dltwu_{3} \),  \( \crossgrad \), \( \dltwun \), and \( \dltwunn \).
We do not require \( \IFT_{\tarpv\nuiv} = 0 \), however, it is implicitly assumed that this operator 
is close to zero.
To quantify this condition, introduce the dual norm of an operator
\( \BBH \colon \R^{\dimq} \to \R^{\dimp} \):
\begin{EQA}
	\| \BBH \|_{\dual}
	&=&
	\sup_{\zv \colon \| \zv \|_{\nano} \leq 1} \| \BBH \zv \| .
\label{ygtuesdfhiwdsoif9rutaG}
\end{EQA}
If \( p=1 \) and \( \| \cdot \|_{\nano} \) is the sup-norm \( \| \cdot \|_{\infty} \) then 
\( \| \BBH \|_{\dual} = \| \BBH \|_{1} \).
Define
\begin{EQ}[rcccl]
	\crossnorm 
	& \eqdef &
	\| \DPN^{-1} \IFT_{\targv\nuiv} \, \HPN^{-1} \|_{\dual} \, ,
	\qquad
	\crossgrad
	& \eqdef &
	\crossnorm + \dltwun \, \rru /2
	\, .
\label{f8vuehery6gv65ftehweeG}
\end{EQ}
Now we apply Proposition~\ref{Pconcsupp} yielding the following \emph{concentration bound} and
\emph{semiparametric Fisher expansion}.

\begin{theorem}
\label{PFisemi}
Let \( \hat{\nupv} \in \Eta_{0} \) on \( \Omega(\xx) \), where \( \Eta_{0} \) is the local set 
from \eqref{8xjxcwbd6dejdccyw3du} with a norm \( \| \cdot \|_{\nano} \),
a metric tensor \( \HPN \), and a radius \( \rru \).
Let also \nameref{LLpT3ref} and 
\nameref{LLsoT3ref} with \( \DFN_{\nuiv} \equiv \DPN \) hold for a metric tensor
\( \DPN \) and some constants \( \dltwu_{3} \), \( \dltwun \), \( \dltwunn \),
and \( \rr \) such that 
\begin{EQA}[c]
	\DPN^{2} \leq \dmax^{2} \IF \, ,
	\quad
	\dltwbun \eqdef \dmax^{2} \dltwunn \rru \leq 1/4 \, ,
	\quad
	\rr 
	\geq 
	\frac{3 \dmax^{2} \crossgrad}{2(1 - \dltwbun)} \, \rru \, ,
	\quad
	\frac{\dmax^{4} \crossgrad \, \dltwu_{3}}{1 - \dltwbun} \, \rru
	\leq 
	\frac{4}{9} 
	\, ,
	\qquad
\label{6cuydjd5eg3jfggu8eyt4hG}
\end{EQA}
for \( \crossnorm \) and \( \crossgrad \) from \eqref{f8vuehery6gv65ftehweeG}.
Then it holds on \( \Omega(\xx) \) 
\begin{EQA}
	\| \DPN (\hat{\tarpv} - \tarpvs) \|
	& \leq &
	\frac{3 \dmax \, \crossgrad}{2(1 - \dltwbun)} \, 
	\bigl( \dmax^{2} \crossgrad \, \| \HPN (\hat{\nuiv} - \nuivs) \|_{\nano} 
		+ \| \DPN \, \IF^{-1} \scorem{\tarpv} \zeta \|
	\bigr) 
	\, ,
\label{rhDGtuGmusGU0a2sptG}
\end{EQA}
and, moreover, 
\begin{EQA}
	&& \nquad
	\| \QP \{ \hat{\tarpv} - \tarpvs 
		+ \IF^{-1} \IFT_{\tarpv\nuiv} (\hat{\nuiv} - \nuivs) 
		- \IF^{-1} \scorem{\tarpv} \zeta
	\} 
	\|
	\\
	& \leq &
	\| \QP \, \IF^{-1} \DPN \| \, 
	\, \Bigl\{ (\dltwuns + \dmax \, \dltwunn) \, \| \HPN (\hat{\nuiv} - \nuivs) \|_{\nano}^{2}
		+ (2 \dltwu_{3} + \dmax^{-1} \dltwunn/2) \, \| \DPN \, \IF^{-1} \scorem{\tarpv} \zeta \|^{2}
	\Bigr\}
	\, ,
	\quad
\label{usdhyw6hikhurnsem}
\end{EQA}
where 
\begin{EQA}
	\dltwuns
	& \eqdef &
	\frac{1}{1 - \dltwbun} \, 
	\biggl( \dmax^{2} \crossnorm \, \dltwunn + \frac{\dltwun}{2}
		+ \frac{3\dmax^{4} \crossgrad^{2} \, \dltwu_{3}}{4(1 - \dltwbun)^{2}} 
	\biggr) \, .
\label{36gfijh94ejdvtwekoisedG}
\end{EQA}
\end{theorem}

\Subsection{Risk of the plug-in estimator}
This section describes concentration sets of the semiparametric plug-in estimator 
\( \hat{\tarpv} = \tilde{\tarpv}(\hat{\nupv}) \)
and provides some bounds on the squared localized risk 
\( \E \, \{ \| \hat{\tarpv} - \tarpvs \|^{2} \Ind_{\Omega(\xx)} \} \) 
using the Fisher expansion \eqref{usdhyw6hikhurnsem}.
Informally it can be written as 
\begin{EQA}
	\hat{\tarpv} - \tarpvs
	& \approx &
	- \IF^{-1} \IFT_{\tarpv\nupv} (\hat{\nupv} - \nupvs)
	+ \IF^{-1} \nabla \zeta
	 \, .
\label{uye94oyuhrt8ytgiutr}
\end{EQA}
This decomposition describes two sources of the estimation loss:
the \emph{variance} (stochastic) term \( \IF^{-1} \nabla \zeta \) is 
due to random errors in observations
and the \emph{semiparametric bias} term \( \IF^{-1} \IFT_{\tarpv\nupv} (\hat{\nupv} - \nupvs) \)
is due to the use of the pilot \( \hat{\nupv} \) in place of the truth \( \nupvs \).

\begin{theorem}
\label{Tsemirisk}
Under conditions of Theorem~\ref{PFisemi}, it holds 
\begin{EQA}
	&& \nquad
	\Bigl| \E \bigl\{ \| \QP (\hat{\tarpv} - \tarpvs) \| \Ind_{\Omega(\xx)} \bigr\}
	- \E \Bigl\{ \bigl\| \QP \bigl\{ 
			\IF^{-1} \nabla \zeta 
			- \IF^{-1} \IFT_{\tarpv\nuiv} (\hat{\nupv} - \nupvs) 
		\bigr\} \bigr\| \Ind_{\Omega(\xx)} \Bigr\} 
	\Bigr|
	\\
	& \leq &
	\| \QP \, \IF^{-1} \DPN \| \, 
	\Bigl\{ (\dltwuns + \dmax \, \dltwunn) \, \dimA_{\HPN}
		+ (2 \dltwu_{3} + \dmax^{-1} \dltwunn/2) \, \dimA_{\DPN}
	\Bigr\}
	\, ,
	\qquad
\label{fgvnjiuet65hftdjhey6uji}
\end{EQA}
where
\begin{EQA}
	\dimA_{\DPN}
	& \eqdef &
	\E \, \| \DPN \, \IF^{-1} \scorem{\tarpv} \zeta \|^{2} \, ,
	\qquad
	\dimA_{\HPN}
	\eqdef
	\E\bigl\{ \| \HPN (\hat{\nuiv} - \nuivs) \|_{\nano}^{2} \Ind_{\Omega(\xx)} \bigr\}
	\leq 
	\rru^{2} \, .
\label{du8dufh8uwhe8723r4fcjue}
\end{EQA}
Moreover, if
\begin{EQA}
	\E \, \bigl\{ \| \HPN (\hat{\nupv} - \nupvs) \|_{\nano}^{4} \Ind_{\Omega(\xx)} \bigr\}
	& \leq &
	\CONSTi_{\HPN}^{2} \, \dimA_{\HPN}^{2} \, ,
	\qquad
	\E \, \bigl\{ \| \DPN \, \IF^{-1} \scorem{\tarpv} \zeta \|^{4} \Ind_{\Omega(\xx)} \bigr\}
	\leq 
	\CONSTi_{\DPN}^{2} \, \dimA_{\DPN}^{2} \, ,
\label{du8jw3d76uy2if523jhv}
\end{EQA}
then it holds with \( \riskt_{\QP} \eqdef 
\E \bigl\{ \| \QP \IF^{-1} \nabla \zeta - \QP \IF^{-1} \IFT_{\tarpv\nuiv} (\hat{\nupv} - \nupvs) \|^{2} \Ind_{\Omega(\xx)} \bigr\} \),
\begin{EQA}
	&& \nquad
	\biggl| \sqrt{\E \, \| \QP (\hat{\tarpv} - \tarpvs) \|^{2} \Ind_{\Omega(\xx)} }
	- \sqrt{\riskt_{\QP}} \biggr|
	\\
	& \leq &
	\| \QP \, \IF^{-1} \DPN \| \, 
	\Bigl\{ 
		(\dltwuns + \dmax \, \dltwunn) \, \CONSTi_{\HPN} \, \dimA_{\HPN}	
		+ (2 \dltwu_{3} + \dmax^{-1} \dltwunn/2) \CONSTi_{\DPN} \, \dimA_{\DPN} 
	\Bigr\} \, .
	\qquad
\label{jv7euejrftdhqttydhj}
\end{EQA}
\end{theorem}

\begin{proof}
Bound \eqref{fgvnjiuet65hftdjhey6uji} follows by \eqref{usdhyw6hikhurnsem}.
Define \( \epsv_{\QP} =	\QP \bigl\{ \hat{\tarpv} - \tarpvs 
- \IF^{-1} \IFT_{\tarpv\nuiv} (\hat{\nupv} - \nupvs)
- \IF^{-1} \, \nabla \zeta \bigr\} \).
Then again by \eqref{usdhyw6hikhurnsem}
\begin{EQA}
	&& 
	\Bigl| \sqrt{\E \, \| \QP (\hat{\tarpv} - \tarpvs) \|^{2} \Ind_{\Omega(\xx)} }
	- \sqrt{\riskt_{\QP}} \Bigr|
	\leq 
	\sqrt{ \E \, \| \epsv_{\QP} \|^{2} \Ind_{\Omega(\xx)} } 
	\leq 
	\| \QP \, \IF^{-1} \DPN \| \times
	\\
	&&
	\times \biggl\{ (\dltwuns + \dmax \, \dltwunn) \, \sqrt{\E \| \HPN (\hat{\nuiv} - \nuivs) \|_{\nano}^{4} \Ind_{\Omega(\xx)} }
	+ \Bigl( 2 \dltwu_{3} + \frac{\dltwunn}{2\dmax} \Bigr) \, 
	\sqrt{\E \| \DPN \, \IF^{-1} \scorem{\tarpv} \zeta \|^{4} \Ind_{\Omega(\xx)} }
	\biggr\}
	\, ,
\label{hc6y2jwei9eiwudhjeei}
\end{EQA}
and the last assertion follows as well.
\end{proof}

\Subsection{Remainders in the expansions}
\label{Sremsemi}

Expansions \eqref{fgvnjiuet65hftdjhey6uji} and \eqref{jv7euejrftdhqttydhj} involve some remainders.
Here we briefly sketch the conditions ensuring a small impact of these remainders.
To be more certain, fix \( \QP = \DPN \) yielding \( \| \QP \, \IF^{-1} \DPN \| \leq \dmax^{2} \).
For more transparency, we also assume \( \dmax = 1 \), \( \dltwu_{3} \asymp n^{-1/2} \) as in \nameref{LLtS3ref} 
for the ``sample size'' \( n \)
and similarly \( \dltwun \asymp n^{-1/2} \), \( \dltwunn \asymp n^{-1/2} \) yielding \( \dltwuns \asymp n^{-1/2} \).
Expansions \eqref{fgvnjiuet65hftdjhey6uji} and \eqref{jv7euejrftdhqttydhj} are meaningful if the remainder terms 
\( (\dltwuns + \dltwunn) \, \| \HPN (\hat{\nuiv} - \nuivs) \|_{\nano}^{2} \)
and \( (2 \dltwu_{3} + \dltwunn/2) \, \| \DPN \, \IF^{-1} \scorem{\tarpv} \zeta \|^{2} \)
are small compared with the magnitude of \( \| \DPN \, \IF^{-1} \scorem{\tarpv} \zeta \| \).
%
%
The corresponding relation can be formulated as
\begin{EQA}
	n^{-1/2} (\dimA_{\DPN} + \dimA_{\HPN})
	& = &
	o(\sqrt{\dimA_{\DPN}}) \, .
\label{dfgt8uy43fvrhi3f5}
\end{EQA}
It is important to stress that the value 
\( \dimA_{\HPN} = \E\{ \| \HPN (\hat{\nuiv} - \nuivs) \|_{\nano}^{2} \Ind_{\Omega(\xx)} \} \)
only involves the \( \| \cdot \|_{\nano} \)-norm of the error \( \HPN (\hat{\nuiv} - \nuivs) \) of the pilot \( \hat{\nupv} \).
For the case of sup-norm estimation, this gradually helps to relax the critical dimension condition.
Section~\ref{ScoBTL} provides some empirical evidence of this conclusion. 

\Subsection{Semiparametric adaptivity}
\label{Scritsemi}
This section discusses the issue of semiparametrically adaptive estimation. 
The question under study is the set of conditions ensuring for the plug-in estimator \( \hat{\tarpv} \)
an ``oracle'' accuracy corresponding to the case of known nuisance parameter \( \nupv \).
Expansion \eqref{usdhyw6hikhurnsem} involves the linear term 
\( \DPN \, \IF^{-1} \IFT_{\tarpv\nuiv} (\hat{\nuiv} - \nuivs) \),
it should be smaller in magnitude than the term \( \DPN \, \IF^{-1} \scorem{\tarpv} \zeta \).
In a high-dimensional situation, the squared norm \( \| \DPN \, \IF^{-1} \scorem{\tarpv} \zeta \|^{2} \)
concentrates around its expectation \( \dimA_{\DPN} \).
By \eqref{f8vuehery6gv65ftehweeG}, it holds on \( \Omega(\xx) \)
\begin{EQA}
	\| \DPN \, \IF^{-1} \IFT_{\tarpv\nuiv} (\hat{\nuiv} - \nuivs) \|
	& \leq &
	\crossnorm \| \HPN (\hat{\nuiv} - \nuivs) \|_{\nano} 
	\leq 
	\crossnorm \, \rru 
	\, ,
\label{yusgfu73uwikrf7g787re}
\end{EQA}
and it is sufficient to check that \( \crossnorm \, \rru \ll \dimA_{\DPN}^{1/2} \).
This is the most critical condition, it requires the cross-correlation between the target and the nuisance parameter
to be sufficiently small. 
The next section illustrates this finding for the sup-norm estimation.

%% file: semi_sup.tex
\Section{Sup-norm bounds for the maximum likelihood estimator}
\label{Ssupsemi}
Let \( \tilde{\upsv} = \argmax_{\upsv} \LL(\upsv) \) be the MLE for a (quasi) log-likelihood function \( \LL(\upsv) \).
The target \( \upsvs \) is defined by maximization of \( \fs(\upsv) = \E \LL(\upsv) \).
This section discusses the problem of bounding the error of estimation \( \tilde{\upsv} - \upsvs \) in sup-norm
(component-wise).
The results from 
\cite{LiRi2022} 
and \cite{GSZ2023} 
demonstrate that accuracy bounds for the MLE can be substantially improved 
compared to the estimation in the more conventional \( \ell_{2} \)-norm.
This paper offers a slightly different viewpoint on this phenomenon 
using the ideas from semiparametric estimation. 

Our results assume that \( \fs(\upsv) \) is concave and the Fisher information matrix \( \IFT = - \nabla^{2} \fs(\upsvs) \)
is well posed at the true point.
We also assume stochastic linearity, that is, \( \zeta(\upsv) = \LL(\upsv) - \E \LL(\upsv) \) is linear in \( \upsv \) 
and the gradient \( \nabla \zeta(\upsv) \) does not depend on \( \upsv \).
The proposed approach suggests treating each one entry \( \tilde{\ups}_{j} \) of \( \tilde{\upsv} \) as a target and the remaining coordinates 
as a nuisance parameter and apply the results of 
Proposition~\ref{Psemibiassup} from Section~\ref{Ssupnormco}. 
Fix a diagonal metric tensor \( \DPN = \diag\{ \DPN_{1},\ldots,\DPN_{\dimp} \} \) with \( \DPN_{j}^{2} = \IFT_{jj} \).
The main condition for applicability of the proposed approach is \( \crosssup \ll 1 \) for
\begin{EQA}[c]
	\crosssup^{2}
	= 
	\max_{j = 1,\ldots,\dimp} \frac{1}{\DPN_{j}^{2}}
	\sum_{m \neq j} \frac{\IFT_{jm}^{2}}{\DPN_{m}^{2}} \, .
\end{EQA}
Also, we need a bound on \( \DPN^{-1} \nabla \zeta \) in a sup-norm.
We keep the result as general as possible, we only assume that on a random set \( \Omega(\xx) \) of a high probability, it holds
\begin{EQA}
	\frac{\sqrt{2}}{1 - \crosssup} \| \DPN^{-1} \nabla \zeta \|_{\infty} 
	& \leq &
	\rrinf
	\, 
\label{8drywjwuycrt5sdeteyhfuwke}
\end{EQA}
for some fixed \( \rrinf \).
Local smoothness of \( \fs(\upsv) \) in a vicinity of \( \upsvs \) will be described by condition \nameref{LLpsupref}
from Section~\ref{Ssupnormco}.
The next result follows directly from Proposition~\ref{Psemibiassup}.

\begin{theorem}
\label{Tsemibiassup}
Let \( \fs(\upsv) = \E \LL(\upsv) \) be concave function. 
For the diagonal metric tensor \( \DPN^{2} = \diag(\IFT) \), assume \( \crosssup < 1 \) and let 
\eqref{8drywjwuycrt5sdeteyhfuwke} hold on \( \Omega(\xx) \) for some \( \rrinf \).
Let \( \fs(\upsv) \) follow \nameref{LLpsupref} with this \( \rrinf \) and
\( \dltwun, \dltwunn, \dltwu_{3} \) satisfying  
\begin{EQA}[c]
	\dltwbun = \dltwunn \, \rrinf \leq 1/4 ,
	\qquad
	\dltwun \, \rrinf \leq 1/4 \, ,
	\qquad	
	\dltwunss \| \DPN^{-1} \nabla \zeta \|_{\infty} 
	\leq 
	\sqrt{2} - 1
	\, ,
\label{dhfiejfowelocuyehbrfe}
\end{EQA}
where
\begin{EQA}[rcl]
	\dltwuns
	& \eqdef &
	\frac{1}{1 - \dltwbun} \, 
	\biggl( \crosssup \, \dltwunn + \frac{\dltwun}{2}
		+ \frac{3(\crosssup + \dltwbun/2)^{2} \, \dltwu_{3}}{4(1 - \dltwbun)^{2}} 
	\biggr) 
	\, ,
\label{36gfijh94ejdvtwekoisedie}
	\\
	\dltwunss
	& \eqdef &
	2 \dltwu_{3} + \frac{\dltwunn}{2}  + \frac{2(\dltwuns + \dltwunn)}{(1 - \crosssup)^{2}} 
	\, .
\end{EQA}
Then \( \tilde{\upsv} \) satisfies on \( \Omega(\xx) \)
\begin{EQA}
	\| \DPN (\tilde{\upsv} - \upsvs) \|_{\infty}
	\eqdef
	\max_{j \leq \dimp} | \DPN_{j} (\tilde{\ups}_{j} - \upss_{j}) |
	& \leq &
	\rrinf \, .
\label{7ytdufchskmls7rghnke}
\end{EQA} 
Furthermore, 
\begin{EQA}
\label{usdhyw6hikhurnetrsppe}
	\bigl\| \DPN^{-1} \bigl\{ \IFT (\tilde{\upsv} - \upsvs) - \nabla \zeta \bigr\} \bigr\|_{\infty}
	& \leq &
	\dltwunss \, \| \DPN^{-1} \nabla \zeta \|_{\infty}^{2} 
	\, ,
	\\
	\bigl\| \DPN (\tilde{\upsv} - \upsvs - \IFT^{-1} \nabla \zeta) \bigr\|_{\infty}
	& \leq &
	\frac{\dltwunss}{1 - \crosssup} \, \| \DPN^{-1} \nabla \zeta \|_{\infty}^{2} 
	\, ,
	\qquad
\label{usdhyw6hikhurnetrspp1e}
	\\
	\bigl\| \DPN (\tilde{\upsv} - \upsvs) - \DPN^{-1} \nabla \zeta) \bigr\|_{\infty}
	& \leq &
	\frac{\dltwunss}{1 - \crosssup} \, \| \DPN^{-1} \nabla \zeta \|_{\infty}^{2} + \frac{\crosssup}{1 - \crosssup} \, \| \DPN^{-1} \nabla \zeta \|_{\infty}
	\, .
\label{usdhyw6hikhurnetrspp2e}
\end{EQA}
\end{theorem}

As a corollary, we obtain a component-wise expansion: on \( \Omega(\xx) \)
\begin{EQA}
	\bigl| \DPN_{j} (\tilde{\ups}_{j} - \upss_{j}) - \DPN_{j}^{-1} \nabla_{j} \zeta) \bigr|
	& \leq &
	\frac{\dltwunss}{1 - \crosssup} \, \| \DPN^{-1} \nabla \zeta \|_{\infty}^{2} + \frac{\crosssup}{1 - \crosssup} \, \| \DPN^{-1} \nabla \zeta \|_{\infty}
	\\
	& \leq &
	\frac{\dltwunss}{2} \, \rrinf^{2} + \frac{\crosssup}{\sqrt{2}} \, \rrinf \, .
\label{6cty6v6et3re4hfdyw}
\end{EQA}

%% file: semi_opt_short.tex

\def\AUv{\Delta}
\def\cmax{\CONST}
\def\upsn{\ups^{\circ}}
\def\pert{t}

\Section{Conditional and marginal optimization}
\label{Ssemiopt}
This section describes the problem of conditional/partial and marginal optimization. 
Consider a function \( \fs(\prmtv) \) of a parameter \( \prmtv \in \R^{\dimttl} \)
which can be represented as \( \prmtv = (\targv,\nuiv) \), where \( \targv \in \R^{\dimp} \)
is the target subvector while \( \nuiv \in \R^{\dimq} \) is a nuisance variable.
Our goal is to study the solution to the optimization problem 
\( \prmtvs = (\targvs,\nuivs) = \argmax_{\prmtv} \fs(\prmtv) \) and, in particular, 
its target component \( \targvs \).
Later we consider a localized setup with a set \( \Ups \) of pairs \( \prmtv = (\targv,\nuiv) \) 
to be fixed around \( \prmtvs \).

\Subsection{Partial optimization}
For any fixed value of the nuisance variable \( \nuiv \in \Nui \), consider 
\( \fs_{\nuiv}(\targv) = \fs(\targv,\nuiv) \) as a function of \( \targv \) only.
Below we assume that \( \fs_{\nuiv}(\targv) \) is concave in \( \targv \) for any \( \nuiv \in \Nui \).
Define 
\begin{EQA}[rcl]
	\targv_{\nuiv}
	& \eqdef &
	\argmax_{\targv \colon (\targv,\nuiv) \in \Ups} \fs_{\nuiv}(\targv) \, .
\label{jfoiuy2wedfv7tr2qsdzxdfso}
\end{EQA}
Our goal is to describe variability of the partial solution \( \targv_{\nuiv} \) in \( \nuiv \) in terms of \( \targv_{\nuiv} - \targvs \)
and \( \fs(\prmtvs) - \fs_{\nuiv}(\targv_{\nuiv}) \).
Introduce 
\begin{EQ}[rcccc]
	\Av_{\nuiv}
	& \eqdef &
	\nabla \fs_{\nuiv}(\targvs)
	&=&
	\nabla_{\targv} \fs(\targvs,\nuiv) 
	\, ,
	\\
	\IFN_{\nuiv}
	& \eqdef &
	- \nabla^{2} \fs_{\nuiv}(\targvs) 
	&=&
	- \nabla_{\targv\targv}^{2} \fs(\targvs,\nuiv) 
	\, .
\label{ge8qwefygw3qytfyju8qfhbdso}
\end{EQ}
Local smoothness of each function \( \fs_{\nuiv}(\cdot) \) around \( \targv_{\nuiv} \) 
can be well described under the self-concordance property.
Let for any \( \nuiv \in \Nui \), some radius \( \rrn_{\nuiv} \) be fixed.
We also assume that a local metric on \( \R^{\dimp} \)
for the target variable \( \targv \) is defined by a matrix \( \DFN_{\nuiv} \in \Matr_{\dimp} \)
that may depend on \( \nuiv \in \Nui \).
Later we assume \nameref{LLsT3ref} 
to be fulfilled for all \( \fs_{\nuiv} \), \( \nuiv \in \Nui \).

\begin{description}

    \item[\label{LLsoT3ref} \( \bb{(\mathcal{T}_{3|\nuiv}^{*})} \)]
      \emph{For \( \nuiv \in \Nui \), 
      it holds
\begin{EQA}
	\sup_{\uv \in \R^{\dimp} \colon \| \DFN_{\nuiv} \uv \| \leq \rrn_{\nuiv}}  \,\, \sup_{\zv \in \R^{\dimp}} \,\,
	\frac{\bigl| \langle \nabla^{3} \fs_{\nuiv}(\targv_{\nuiv} + \uv), \zv^{\otimes 3} \rangle \bigr|}
		 {\| \DFN_{\nuiv} \zv \|^{3}}
	& \leq &
	\dltwu_{3} \, .
\label{d6f53ye5vry4fddfgeyd}
\end{EQA}
}

\end{description}

Our first result describes the \emph{semiparametric bias} \( \targvs - \targv_{\nuiv} \) caused by using
the value \( \nuiv \) of the nuisance variable in place of \( \nuivs \).

\begin{proposition}
\label{PbiassemiN}
Let \( \fs_{\nuiv}(\targv) \) be a strongly concave function with 
\( \fs_{\nuiv}(\targv_{\nuiv}) = \max_{\targv} \fs_{\nuiv}(\targv) \)  
and \( \IFN_{\nuiv} = - \nabla^{2} \fs_{\nuiv}(\targv_{\nuiv}) \).
Let \( \Av_{\nuiv} \) and \( \IFN_{\nuiv} \) be given by \eqref{ge8qwefygw3qytfyju8qfhbdso}.
Assume \nameref{LLsoT3ref} at \( \targv_{\nuiv} \) with \( \DFN_{\nuiv}^{2} \), \( \rrn_{\nuiv} \), and \( \dltwu_{3} \) such that 
\begin{EQA}[c]
	\DFN_{\nuiv}^{2} \leq \dmax^{2} \, \IFN_{\nuiv} ,
	\quad
	\rrn_{\nuiv} \geq \frac{3}{2} \| \DFN_{\nuiv} \, \IFN_{\nuiv}^{-1} \Av_{\nuiv} \| \, ,
	\quad
	\dmax^{2} \dltwu_{3} \| \DFN_{\nuiv} \, \IFN_{\nuiv}^{-1} \Av_{\nuiv} \| < \frac{4}{9} \, .
\label{8difiyfc54wrboer7bjfr}
\end{EQA}
Then \( \| \DFN_{\nuiv} (\targv_{\nuiv} - \targvs) \| \leq (3/2) \| \DFN_{\nuiv} \, \IFN_{\nuiv}^{-1} \Av_{\nuiv} \| \) and moreover,
\begin{EQ}[rcccl]
    \| \DFN_{\nuiv}^{-1} \IFN_{\nuiv} (\targv_{\nuiv} - \targvs + \IFN_{\nuiv}^{-1} \Av_{\nuiv}) \|
    & \leq &
    \frac{3\dltwu_{3}}{4} \| \DFN_{\nuiv} \, \IFN_{\nuiv}^{-1} \Av_{\nuiv} \|^{2} 
    & \leq &
    \frac{\dltwu_{3} \, \rrn_{\nuiv}}{2} \| \DFN_{\nuiv} \, \IFN_{\nuiv}^{-1} \Av_{\nuiv} \|
	\, .
	\qquad
\label{jhcvu7ejdytur39e9frtfw}
\end{EQ}
Moreover,
\begin{EQA}
	\bigl| 
		2 \fs_{\nuiv}(\targv_{\nuiv}) - 2 \fs_{\nuiv}(\targvs) - \| \IFN_{\nuiv}^{-1/2} \Av_{\nuiv} \|^{2} 
	\bigr|
	& \leq &
	\frac{5\dltwu_{3}}{2} \, \| \DFN_{\nuiv} \, \IFN_{\nuiv}^{-1} \Av_{\nuiv} \|^{3} \, .
	\qquad
\label{gtxddfujhyfdytet6ywerfd}
\end{EQA}
\end{proposition}

\Subsection{Conditional optimization under (semi)orthogonality}
\label{SortLaplsemi}
Here we study variability of the value \( \targv_{\nuiv} = \argmax_{\targv} \fs(\targv,\nuiv) \)
w.r.t. the nuisance para\-meter \( \nuiv \).
It appears that local quadratic approximation of the function \( \fs \) in a vicinity of the extreme point
\( \prmtvs \) yields a nearly linear dependence of \( \targv_{\nuiv} \) on \( \nuiv \).
We illustrate this fact on the case of a quadratic function \( \fs(\cdot) \).
Consider the negative Hessian \( \IFT = - \nabla^{2} \fs(\prmtvs) \) in the block form:
\begin{EQA}
	\IFT
	& \eqdef &
	- \nabla^{2} \fs(\prmtvs)
	=
	\begin{pmatrix}
	\IFT_{\targv\targv} & \IFT_{\targv\nuiv}
	\\
	\IFT_{\nuiv\targv} & \IFT_{\nuiv\nuiv}
	\end{pmatrix} 
\label{hwe78yf2diwe76tfw67etfwtbso}
\end{EQA}
with \( \IFT_{\nuiv\targv} = \IFT_{\targv\nuiv}^{\T} \).
If \( \fs(\prmtv) \) is quadratic then \( \IFT \) and its blocks do not depend on \( \prmtv \).

\begin{lemma}
\label{Lpartmaxq}
Let \( \fs(\prmtv) \) be {quadratic}, {strongly concave}, and \( \nabla \fs(\prmtvs) = 0 \).
Then
\begin{EQA}
	\targv_{\nuiv} - \targvs 
	&=& 
	- \IFT_{\targv\targv}^{-1} \, \IFT_{\targv\nuiv} \, (\nuiv - \nuivs) .
\label{0fje7fhihy84efiewkw}
\end{EQA}
\end{lemma}

\begin{proof}
The condition \( \nabla \fs(\prmtvs) = 0 \) implies 
\( \fs(\prmtv) = \fs(\prmtvs) - (\prmtv - \prmtvs)^{\T} \IFT \, (\prmtv - \prmtvs)/2 \) 
with \( \IFT = - \nabla^{2} \fs(\prmtvs) \).
For \( \nuiv \) fixed, the point \( \targv_{\nuiv} \) maximizes 
\( - (\targv - \targvs)^{\T} \IFT_{\targv\targv} \, (\targv - \targvs) /2 - (\targv - \targvs)^{\T} \IFT_{\targv\nuiv} \, (\nuiv - \nuivs) \)
and thus, \( \targv_{\nuiv} - \targvs = - \IFT_{\targv\targv}^{-1} \, \IFT_{\targv\nuiv} \, (\nuiv - \nuivs) \).
\end{proof} 

This observation \eqref{0fje7fhihy84efiewkw} is in fact discouraging 
because the bias \( \targv_{\nuiv} - \targvs \) has the same magnitude as the nuisance parameter \( \nuiv - \nuivs \).
However, the condition \( \IFT_{\targv\nuiv} = 0 \) yields \( \targv_{\nuiv} \equiv \targvs \) 
and the bias vanishes.
If \( \fs(\prmtv) \) is not quadratic, the \emph{orthogonality} condition \( \nabla_{\nuiv} \nabla_{\targv} \, \fs(\targv,\nuiv) \equiv 0 \)
for all \( (\targv,\nuiv) \in \WV \) still ensures a vanishing bias.

\begin{lemma}
\label{Lpartmax}
Let \( \fs(\targv,\nuiv) \) be continuously differentiable and 
\( \nabla_{\nuiv} \nabla_{\targv} \, \fs(\targv,\nuiv) \equiv 0 \).
Then the point 
\( \targv_{\nuiv} = \argmax_{\targv} f(\targv,\nuiv) \) does not depend on \( \nuiv \).
\end{lemma}

\begin{proof}
The condition \( \nabla_{\nuiv} \nabla_{\targv} \, \fs(\targv,\nuiv) \equiv 0 \) implies the decomposition 
\( \fs(\targv,\nuiv) = \fs_{1}(\targv) + \fs_{2}(\nuiv) \) for some functions 
\( \fs_{1} \) and \( \fs_{2} \).
This in turn yields \( \targv_{\nuiv} \equiv \targvs \).
\end{proof}

In some cases, one can check \emph{semi-orthogonality} condition 
\begin{EQA}
	\nabla_{\nuiv} \nabla_{\targv} \, \fs(\targvs,\nuiv)  
	&=&
	0,
	\qquad
	\forall \nuiv \in \Nui \, .
\label{jdgyefe74erjscgygfydhwse}
\end{EQA}
\ifapp{A typical example is given by nonlinear regression; see Section~\ref{SnonlinLA}.}{}

\begin{lemma}
\label{Lsemiorto}
Assume \eqref{jdgyefe74erjscgygfydhwse}.
Then 
\begin{EQA}[rcccl]
	\nabla_{\targv} \, \fs(\targvs,\nuiv)
	& \equiv &
	0,
	\qquad
	\nabla_{\targv\targv}^{2} \fs(\targvs,\nuiv) 
	& \equiv &
	\nabla_{\targv\targv}^{2} \fs(\targvs,\nuivs),
	\qquad
	\nuiv \in \Nui .
\label{0xyc5ftr4j43vefvuruendt}
\end{EQA}
Moreover, if \( \fs(\targv,\nuiv) \) is concave in \( \targv \) given \( \nuiv \) then
\begin{EQA}
	\targv_{\nuiv} 
	& \eqdef & 
	\argmax_{\targv} \fs(\targv,\nuiv) \equiv \targvs ,
	\qquad
	\forall \nuiv \in \Nui \,.
\label{irdtyrnjjutu45r7gjhr}
\end{EQA}
\end{lemma}

\begin{proof}
Consider the vector \( \Av_{\nuiv} = \nabla_{\targv} \, \fs(\targvs,\nuiv) \).
Obviously \( \Av_{\nuivs} = 0 \).
Moreover, \eqref{jdgyefe74erjscgygfydhwse} implies that \( \Av_{\nuiv} \) 
does not depend on \( \nuiv \) and thus, vanishes everywhere.
As \( \fs \) is concave in \( \targv \), 
this implies \( \fs(\targvs,\nuiv) = \max_{\targv} \fs(\targv,\nuiv) \)
and \( \targv_{\nuiv} = \argmax_{\targv} \fs(\targv,\nuiv) \equiv \targvs \).
Similarly, by \eqref{jdgyefe74erjscgygfydhwse}, it holds
\( \nabla_{\nuiv} \nabla_{\targv\targv} \, \fs(\targvs,\nuiv) \equiv 0 \) 
and \eqref{0xyc5ftr4j43vefvuruendt} follows.
Concavity of \( \fs(\targv,\nuiv) \) in \( \targv \) for \( \nuiv \) fixed and \eqref{0xyc5ftr4j43vefvuruendt} 
imply \eqref{irdtyrnjjutu45r7gjhr}.
\end{proof}

Orthogonality or semi-orthogonality \eqref{jdgyefe74erjscgygfydhwse} conditions are rather 
restrictive and fulfilled only in special situations.
A weaker condition of 
\emph{one-point orthogonality} means \( \nabla_{\nuiv} \nabla_{\targv} \, \fs(\targvs,\nuivs) = 0 \).
This condition is not restrictive and can always be enforced by a linear transform of the nuisance variable \( \nuiv \).

\Subsection{Semiparametric bias}
\label{Sonepoint}


This section explains how one can bound variability of 
\( \IFN_{\nuiv} = - \nabla_{\targv\targv}^{2} \fs(\targvs,\nuiv) \)
and of the norm of \( \targv_{\nuiv} - \targvs \) under  
conditions on the cross-derivatives of \( \fs(\targv,\nuiv) \) for \( \nuiv = \nuivs \)
or \( \targv = \targvs \).
Suppose that the nuisance variable \( \nuiv \) is already localized to a small vicinity \( \Nui \) of \( \nuivs \).
A typical example is given by the level set \( \Nui \) of the form
\begin{EQA}
	\Nui
	&=&
	\bigl\{ \nuiv \colon \| \HFN (\nuiv - \nuivs) \|_{\nano} \leq \rru \bigr\} ,
\label{6chbdtw2hydfuguye4223w23}
\end{EQA}
where \( \HFN \) is a metric tensor in \( \R^{\dimq} \), 
\( \| \cdot \|_{\nano} \) is a norm in \( \R^{\dimq} \), and \( \rru > 0 \).
Often \( \| \cdot \|_{\nano} \) is the usual \( \ell_{2} \)-norm.
However, in some situations, it is beneficial to use different topology for the target parameter \( \targv \)
and the nuisance parameter \( \nuiv \).
One example is given by estimation in \( \sup \)-norm with \( \| \cdot \|_{\nano} = \| \cdot \|_{\infty} \).

Assume the following condition.

\begin{description}


    \item[\label{LLpT3ref} \( \bb{(\mathcal{T}_{3,\Nui}^{*})} \)]
      \emph{
      It holds with some \( \dltwun, \dltwunn \)
}
\begin{EQA}[rccl]
	\sup_{\nuiv \in \Nui} \,\, 
	\sup_{\zv \in \R^{\dimp}, \wv \in \R^{\dimq}} \,
	&
	\frac{| \langle \nabla_{\targv\nuiv\nuiv}^{3} \fs(\targvs,\nuiv), \zv \otimes \wv^{\otimes 2} \rangle |}
		 {\| \DFN \zv \| \, \| \HFN \wv \|_{\nano}^{2}}
	& \leq &
	\dltwun \, ,
	\qquad
\label{c6ceyecc5e5etcT2}
	\\
	\sup_{\nuiv \in \Nui} \,\, 
	\sup_{\zv \in \R^{\dimp}, \wv \in \R^{\dimq}} \,
	&
	\frac{| \langle \nabla_{\targv\targv\nuiv}^{3} \fs(\targvs,\nuiv), \zv^{\otimes 2} \otimes \wv \rangle |}
		 {\| \DFN \zv \|^{2} \, \| \HFN \wv \|_{\nano}}
	& \leq &
	\dltwunn \, .
\label{c6ceyecc5e5etctwhcyegwc}
\end{EQA}
\end{description}

\begin{remark}
Later we establish some bounds on the semiparametric bias assuming 
\( (\dltwun + \dltwunn) \, \rru \ll 1 \). 
\end{remark}

\begin{remark}
Condition \nameref{LLpT3ref} only involves mixed derivatives 
\( \nabla_{\targv\nuiv\nuiv}^{3} \fs(\targvs,\nuiv) \) and \( \nabla_{\targv\targv\nuiv}^{3} \fs(\targvs,\nuiv) \) 
of \( \fs(\targvs,\nuiv) \) for the fixed value \( \targv = \targvs \),
while condition \nameref{LLsoT3ref} only concerns smoothness of \( \fs(\targv,\nuiv) \) w.r.t. \( \targv \) for \( \nuiv \) fixed.
Therefore, the combination of \nameref{LLsoT3ref} and \nameref{LLpT3ref} is much weaker than the full dimensional 
condition \nameref{LLsT3ref}.
\end{remark}

The next result provides an expansion of the semiparametric bias \( \targv_{\nuiv} - \targvs \). 
For the norm \( \| \cdot \|_{\nano} \), 
let \( \| \BBH \|_{\dual} \) be the corresponding dual norm of an operator
\( \BBH \colon \R^{\dimq} \to \R^{\dimp} \):
\begin{EQA}
	\| \BBH \|_{\dual}
	&=&
	\sup_{\zv \colon \| \zv \|_{\nano} \leq 1} \| \BBH \zv \| .
\label{ygtuesdfhiwdsoif9ruta}
\end{EQA}
If \( p=1 \) and \( \| \cdot \|_{\nano} \) is the sup-norm \( \| \cdot \|_{\infty} \) then 
\( \| \BBH \|_{\dual} = \| \BBH \|_{1} \).
Define
\begin{EQA}
\label{f8vuehery6gv65ftehwee}
	\crossnorm
	& \eqdef &
	\| \DFN^{-1} \IFT_{\targv\nuiv} \, \HFN^{-1} \|_{\dual} \, ,
	\\
	\crossgrad
	& \eqdef &
	\| \DFN^{-1} \IFT_{\targv\nuiv} \, \HFN^{-1} \|_{\dual} + \dltwun \, \rru /2 \, .
\label{f8vuehery6gv65ftehwee1}
\end{EQA}
Proposition~\ref{PbiassemiN} evaluates the semiparametric bias \( \targvs - \targv_{\nuiv} \) in terms of 
the matrix 
\( \IFN_{\nuiv} = - \nabla_{\targv\targv}^{2} \fs(\targvs,\nuiv) \) and the vector 
\( \Av_{\nuiv} = \nabla \fs_{\nuiv}(\targvs) \).
The next result explains how these quantities can be controlled under \nameref{LLpT3ref}.
To simplify the formulation, we assume that \( \DFN^{2} \leq \IFN \); see Remark~\ref{Rdmax2} for an extension.

\begin{proposition}
\label{PsemiAvex}
Assume \nameref{LLpT3ref} with \( \DFN^{2} \leq \IFN \).
Fix \( \nuiv \in \Nui \), set \( \dltwbun \eqdef \dltwunn \| \HFN (\nuiv - \nuivs) \|_{\nano} \), 
and assume \nameref{LLsoT3ref} with \( \DFN_{\nuiv} \equiv \DFN \), \( \rrn_{\nuiv} \equiv \rrn \), 
and \( \dltwu_{3} \),
such that 
\begin{EQA}[c]
	\dltwbun < 1 \, ,
	\quad
	\rr 
	\geq 
	\frac{3 \crossgrad}{2(1 - \dltwbun)} \, \| \HFN (\nuiv - \nuivs) \|_{\nano} \, ,
	\quad
	\frac{\crossgrad \, \dltwu_{3}}{1 - \dltwbun} \, \| \HFN (\nuiv - \nuivs) \|_{\nano}
	\leq 
	\frac{4}{9} \, ,
	\qquad
\label{6cuydjd5eg3jfggu8eyt4h}
\end{EQA}
for \( \crossgrad \) from \eqref{f8vuehery6gv65ftehwee1}.
Then the partial solution \( \targv_{\nuiv} \) obeys 
\begin{EQA}
	\| \DFN (\targv_{\nuiv} - \targvs) \|
	& \leq &
	\frac{3 \crossgrad}{2(1 - \dltwbun)} \, \| \HFN (\nuiv - \nuivs) \|_{\nano} \,. 
\label{rhDGtuGmusGU0a2spt}
\end{EQA}
Moreover, with \( \crossnorm \) from \eqref{f8vuehery6gv65ftehwee}, it holds
\begin{EQA}[c]
\label{36gfijh94ejdvtwekoise}
	\| \QP \{ \targv_{\nuiv} - \targvs + \IFN^{-1} \IFT_{\targv\nuiv} (\nuiv - \nuivs) \} \|
	\leq 
	\| \QP \, \IFN^{-1} \DFN \| \, \dltwuns
	\, \| \HFN (\nuiv - \nuivs) \|_{\nano}^{2}
	\, ,
	\qquad
	\\
	\dltwuns
	\eqdef 
	\frac{1}{1 - \dltwbun} \, 
	\biggl( \crossnorm \, \dltwunn + \frac{\dltwun}{2}
		+ \frac{3\crossgrad^{2} \, \dltwu_{3}}{4(1 - \dltwbun)^{2}} 
	\biggr) \, .
\label{36gfijh94ejdvtwekoised}
\end{EQA}
\end{proposition}

\begin{remark}
\label{Rdmax2}
An extension to the case \( \DFN^{2} \leq \dmax^{2} \IFN \) can be done by replacing everywhere
\( \dltwu_{3} \),
\( \dltwun \),
\( \dltwunn \),
\( \dltwuns \),
\( \crossnorm \),
\( \crossgrad \) with \( \dmax^{3} \dltwu_{3} \), \( \dmax \, \dltwun \), 
\( \dmax^{2} \dltwunn \), 
\( \dmax \, \dltwuns \), 
\( \dmax^{-1} \crossnorm \), 
\( \dmax^{-1} \crossgrad \)
respectively.
\end{remark}

\Subsection{A linear perturbation}
Let \( \fn(\targv,\nuiv) \) be a linear perturbation of \( \fs(\targv,\nuiv) \):
\begin{EQA}
	\fn(\targv,\nuiv) - \fn(\targvs,\nuivs)
	&=&
	\fs(\targv,\nuiv) - \fn(\targvs,\nuivs) + \langle \AAv,\targv - \targvs \rangle + \langle \CCv,\nuiv - \nuivs \rangle ;
\label{s6cdt76ws6656ysjwsdftc6}
\end{EQA}
cf. \eqref{4hbh8njoelvt6jwgf09}.
Let also \( \Nui \) be given by \eqref{6chbdtw2hydfuguye4223w23} and \( \nuiv \in \Nui \).
We are interested in quantifying the distance between \( \targvs \) and \( \targvn_{\nuiv} \), where 
\begin{EQA}
	\targvn_{\nuiv}
	& \eqdef &
	\argmax_{\targv} \fn(\targv,\nuiv) .
\label{8dcudxjd5t6t5mkwsyddf}
\end{EQA}
The linear perturbation \( \langle \CCv,\nuiv - \nuivs \rangle \) does not depend on \( \targv \) and it can be ignored.

\begin{proposition}
\label{Pconcsupp}
For \( \nuiv \in \Nui \),
assume the conditions of Proposition~\ref{PsemiAvex} and let 
\begin{EQA}[c]
	\dltwbun \eqdef \dltwunn \| \HFN (\nuiv - \nuivs) \|_{\nano} \leq 1/4 \, .
\end{EQA}
Then 
\begin{EQA}
	\| \DFN (\targvn_{\nuiv} - \targvs) \|
	& \leq &
	\frac{3 }{2(1 - \dltwbun)} \, 
	\bigl( \crossgrad \, \| \HFN (\nuiv - \nuivs) \|_{\nano} + \| \DFN^{-1} \AAv \| \bigr) 
\label{ysd7euwejfg653wthf873k}
\end{EQA}
and for any linear mapping \( \QP \), it holds with 
\( \dltwuns \) from \eqref{36gfijh94ejdvtwekoised} 
\begin{EQA}[l]
	\!\!\!
	\| \QP \{ \targvn_{\nuiv} - \targvs + \IFN^{-1} \IFT_{\targv\nuiv} (\nuiv - \nuivs) - \IFN^{-1} \AAv \} \|
	\\
	\leq 
	\| \QP \, \IFN^{-1} \DFN \| \, \Bigl( \dltwuns \, \| \HFN (\nuiv - \nuivs) \|_{\nano}^{2}
		+ 2 \dltwu_{3} \, \| \DFN \, \IFN^{-1} \AAv \|^{2}
		+ \frac{4\dltwunn}{3} \, \| \DFN \, \IFN^{-1} \AAv \| \, \| \HFN (\nuiv - \nuivs) \|_{\nano}
	\Bigr)	
	\\
	\leq 
	\| \QP \, \IFN^{-1} \DFN \| \, \Bigl\{ (\dltwuns + \dltwunn) \, \| \HFN (\nuiv - \nuivs) \|_{\nano}^{2}
		+ (2 \dltwu_{3} + \dltwunn/2) \, \| \DFN \, \IFN^{-1} \AAv \|^{2}
	\Bigr\}	
	\, .
	\qquad 
\label{usdhyw6hikhurnetr256}
\end{EQA}
\end{proposition}

\begin{proof}
Proposition~\ref{PbiassemiN} applied to \( \fs_{\nuiv}(\targv) \) and 
\( \fn_{\nuiv}(\targv) = \fs_{\nuiv}(\targv) - \langle \AAv, \targv \rangle \) implies  by \eqref{uciuvcdiu3eoooodffcym}
\begin{EQA}
	\| \DFN (\targvn_{\nuiv} - \targv_{\nuiv}) \|
	& \leq &
	\frac{3}{2} \| \DFN \, \IFN_{\nuiv}^{-1} \AAv \|
	\leq 
	\frac{3}{2(1 - \dltwbun)} \| \DFN^{-1} \AAv \| 
	\, .
\label{djvhehfief6yh3jfgiu2}
\end{EQA}
This and \eqref{rhDGtuGmusGU0a2spt} imply \eqref{ysd7euwejfg653wthf873k}.
Next we check \eqref{usdhyw6hikhurnetr256} using the decomposition
\begin{EQA}
	&& \nquad
	\targvn_{\nuiv} - \targvs + \IFN^{-1} \IFT_{\targv\nuiv} (\nuiv - \nuivs) - \IFN^{-1} \AAv
	\\
	&=&
	(\targvn_{\nuiv} - \targv_{\nuiv} - \IFN_{\nuiv}^{-1} \AAv) 
	+ (\IFN_{\nuiv}^{-1} \AAv - \IFN^{-1} \AAv)
	+ \{ \targv_{\nuiv} - \targvs + \IFN^{-1} \IFT_{\targv\nuiv} (\nuiv - \nuivs) \} 
	\, .
\label{yhdt56swhew454erjhfgi}
\end{EQA}
Proposition~\ref{PsemiAvex} evaluates the last term.
Lemma~\ref{L3IFNNui} helps to bound
\begin{EQA}
	\| \DFN^{-1} \IFN (\IFN_{\nuiv}^{-1} - \IFN^{-1}) \IFN_{\nuiv} \, \DFN^{-1} \|
	& = &
	\| \DFN^{-1} (\IFN_{\nuiv} - \IFN) \DFN^{-1} \|
	\leq 
	\dltwunn \| \HFN (\nuiv - \nuivs) \|_{\nano} 
	\, .
\label{cuyedgfcuedgtfbvj}
\end{EQA}
This yields 
\begin{EQA}
	\| \QP (\IFN_{\nuiv}^{-1} - \IFN^{-1}) \AAv \|
	& \leq &
	\| \QP \, \IFN^{-1} \DFN \| 
	\, \| \DFN^{-1} \IFN_{\nuiv} (\IFN_{\nuiv}^{-1} - \IFN^{-1}) \IFN \, \DFN^{-1} \|
	\, \| \DFN \, \IFN_{\nuiv}^{-1} \AAv \|
	\\
	& \leq &
	\| \QP \, \IFN^{-1} \DFN \| 
	\, \dltwunn \, \| \HFN (\nuiv - \nuivs) \|_{\nano} \,
	\, \frac{1}{1 - \dltwbun} \| \DFN \, \IFN^{-1} \AAv \|
	\, .
\label{7duewjf87he4fyjwucyfne}
\end{EQA}
Moreover, 
for \( \epsv_{\nuiv} \eqdef \targvn_{\nuiv} - \targv_{\nuiv} - \IFN_{\nuiv}^{-1} \AAv \), it holds
\begin{EQA}
	\| \DFN^{-1} \IFN_{\nuiv} \, \epsv_{\nuiv} \|
	& \leq &
	\frac{3 \dltwu_{3}}{4} \| \DFN \, \IFN_{\nuiv}^{-1} \AAv \|^{2}  
	\leq 
	\frac{3 \dltwu_{3}}{4(1 - \dltwbun)^{2}} \| \DFN \, \IFN^{-1} \AAv \|^{2}
	\, ,
\label{ud82uegkftryyde6ewhde}
\end{EQA} 
and by \eqref{uciuvcdiu3eoooodffcym} and \( \dltwbun \leq 1/4 \)
\begin{EQA}
	\| \QP \, \epsv_{\nuiv} \|
	& \leq &
	\| \QP \, \IFN^{-1} \DFN \| \,\, \| \DFN^{-1} \IFN \, \epsv_{\nuiv} \|
	\leq 
	\| \QP \, \IFN^{-1} \DFN \| \, \frac{1 - \dltwbun}{1 - 2 \dltwbun} \, \| \DFN^{-1} \IFN_{\nuiv} \, \epsv_{\nuiv} \|
	\\
	& \leq &
	\| \QP \, \IFN^{-1} \DFN \| \,\, 
	\frac{3 \dltwu_{3}}{4(1 - \dltwbun) (1 - 2 \dltwbun)} \| \DFN \, \IFN^{-1} \AAv \|^{2} 
	\leq 
	\| \QP \, \IFN^{-1} \DFN \| \, 2 \dltwu_{3} \, \| \DFN \, \IFN^{-1} \AAv \|^{2}
	\, .
\label{ysdudhwuw8ey87wejue4bf}
\end{EQA}
The obtained bounds imply \eqref{usdhyw6hikhurnetr256} in view of
\( (4/3)ab \leq a^{2} + b^{2}/2 \) for any \( a,b \).
\end{proof}

\Section{Sup-norm expansions for linearly perturbed optimization}
\label{Ssupnormco}
Let \( \fs(\upsv) \) be a concave function and \( \upsvs = \argmax_{\upsv} \fs(\upsv) \).
Define \( \IFT = - \nabla^{2} \fs(\upsvs) = (\IFT_{jm}) \).
Let also \( \fn(\upsv) \) be obtained by a linear perturbation of \( \fs(\upsv) \) as in \eqref{4hbh8njoelvt6jwgf09}
and let \( \upsvn = \argmax_{\upsv} \fn(\upsv) \).
We intend to bound the corresponding change \( \upsvn - \upsvs \) in a sup-norm.
The strategy of the study is to fix one component \( \ups_{j} \) of \( \upsv \) and treat the remaining entries as the nuisance 
parameter.
Fix a metric tensor \( \DFN \) in diagonal form:
\begin{EQA}[c]
	\DFN \eqdef \diag(\DFN_{1}, \ldots, \DFN_{\dimp}).
\label{7dycf8mwdy6ey43hf98yh}
\end{EQA}
Later we assume \( \DFN_{j}^{2} = \IFT_{jj} \).
For each \( j \leq \dimp \), we will use the representation \( \upsv = (\ups_{j},\nuiv_{j}) \), where 
\( \nuiv_{j} \in \R^{\dimp-1} \) collects all the remaining entries of \( \upsv \).
Similarly, define \( \DFN = (\DFN_{j},\HFN_{j}) \) with \( \HFN_{j} = \diag(\DFN_{m}) \) for \( m \neq j \).
Define also for some radius \( \rrinf \)
\begin{EQA}[c]
	\Upsd = \{ \upsv \colon \| \DFN (\upsv - \upsvs) \|_{\infty} \leq \rrinf \} 
	=
	\Bigl\{ \upsv \colon \max_{j \leq \dimp} \bigl| \DFN_{j} (\ups_{j} - \upss_{j}) \bigr| 
	\leq \rrinf \Bigr\}.
\label{du7duydy7cy6fc6df3}
\end{EQA}
Assume the following condition.

\begin{description}
    \item[\label{LLpsupref} \( \bb{(\mathcal{T}_{\infty}^{*})} \)]
      \emph{For each \( j \leq \dimp \), 
     the function \( \fs(\ups_{j},\nuiv_{j}) \) fulfills}
\begin{EQA}[rcl]
	\sup_{\nuiv_{j} \colon \| \HFN_{j} (\nuiv_{j} - \nuivs_{j}) \|_{\infty} \leq \rrinfi} \,\, 
	\sup_{\ups \colon \DFN_{j} |\ups - \upss_{j}| \leq 2\rrinfi} \,
	\frac{\bigl| \nabla_{\ups_{j}\ups_{j} \ups_{j}}^{(3)} \fs(\ups_{j},\nuiv_{j}) \bigr|}
		 {\DFN_{j}^{3}}
	& \leq &
	\dltwu_{3} \, ,
\label{d6f53ye5vry4fddfgsup}
\end{EQA}
\emph{and}
\begin{EQA}[lccl]
	\sup_{\upsv = (\upss_{j},\nuiv_{j}) \in \Upsd} \,
	\sup_{\zv_{j} \in \R^{\dimp-1}} \,
	&
	\frac{| \langle \nabla_{\ups_{j}\ups_{j}\nuiv_{j}}^{(3)} \fs(\upss_{j},\nuiv_{j}), \zv_{j} \rangle |}
		 {\DFN_{j}^{2} \, \| \HFN_{j} \zv_{j} \|_{\infty}}
	& \leq &
	\dltwunn 
	\, ,
\label{c6hyi8ietctwhcsup}
	\\
	\sup_{\upsv = (\upss_{j},\nuiv_{j}) \in \Upsd}  \,
	\sup_{\zv_{j} \in \R_{j}^{\dimp}} \,
	&
	\frac{| \langle \nabla_{\ups_{j}\nuiv_{j}\nuiv_{j}}^{3} \fs(\upss_{j},\nuiv_{j}), \zv_{j}^{\otimes 2} \rangle |}
		 {\DFN_{j} \, \| \HFN_{j} \zv_{j} \|_{\infty}^{2}}
	& \leq &
	\dltwun 
	\, .
	\qquad
\label{c6ceyecc5e5etcT2sup}
\end{EQA}
\end{description}

For the dual norm \( \| \cdot \|_{\dual} \) from \eqref{ygtuesdfhiwdsoif9ruta}, define
\begin{EQA}
	\crosssup
	& \eqdef &
	\max_{j = 1,\ldots,\dimp} \| \DFN_{j}^{-1} \IFT_{\ups_{j},\nuiv_{j}} \, \HFN_{j}^{-1} \|_{\dual} 
	=
	\max_{j = 1,\ldots,\dimp} \,\, \sup_{\| \zv \|_{\infty} \leq 1}
	\| \DFN_{j}^{-1} \IFT_{\ups_{j},\nuiv_{j}} \, \HFN_{j}^{-1} \zv \|
	\, .
	\qquad
\label{7tdsyf8iuwopkrtg4576}
\end{EQA}

\begin{lemma}
\label{Lcrosssup}
It holds
\begin{EQA}[c]
	\crosssup^{2}
	\leq 
	\max_{j = 1,\ldots,\dimp} \frac{1}{\DFN_{j}^{2}}
	\sum_{m \neq j} \frac{\IFT_{jm}^{2}}{\DFN_{m}^{2}} \, .
\end{EQA}
\end{lemma}

\begin{proof}
By definition,
\begin{EQA}[rcl]
	\crosssup^{2}
	& \leq &
	\max_{j = 1,\ldots,\dimp} \,\, \sup_{\| \zv \|_{\infty} \leq 1}
	\sum_{m \neq j} \frac{1}{\DFN_{j}^{2} \, \DFN_{m}^{2}} \, \IFT_{jm}^{2} z_{m}^{2} 
	\leq 
	\max_{j = 1,\ldots,\dimp} \frac{1}{\DFN_{j}^{2}}
	\sum_{m \neq j} \frac{\IFT_{jm}^{2}}{\DFN_{m}^{2}} \,  
\end{EQA}
as claimed.
\end{proof}
It is mandatory for the proposed approach that \( \crosssup < 1 \).
The next result provides an upper bound on \( \| \upsvd - \upsvs \|_{\infty} \).

\begin{proposition}
\label{Psemibiassup}
Let \( \fs(\upsv) \) be concave function and 
\( \fn(\upsv) \) be a linear perturbation of \( \fs(\upsv) \) with a vector \( \AAv \). 
For a diagonal metric tensor \( \DFN \), assume \( \crosssup < 1 \); see \eqref{7tdsyf8iuwopkrtg4576}. 
Fix 
\begin{EQA}
	\rrinf 
	&=& 
	\frac{\sqrt{2}}{1 - \crosssup} \| \DFN^{-1} \AAv \|_{\infty} 
	\, .
\label{8drywjwuycrt5sdeteyhfuwk}
\end{EQA}
Let \nameref{LLpsupref} hold with this \( \rrinf \) and
\( \dltwun, \dltwunn, \dltwu_{3} \) satisfying  
\begin{EQA}[c]
	\dltwbun = \dltwunn \, \rrinf \leq 1/4 ,
	\qquad
	\dltwun \, \rrinf \leq 1/4 \, ,
	\qquad	
	\dltwunss \| \DFN^{-1} \AAv \|_{\infty} 
	\leq 
	\sqrt{2} - 1
	\, ,
\label{dhfiejfowelocuyehbrf}
\end{EQA}
where
\begin{EQA}[rcl]
	\dltwuns
	& \eqdef &
	\frac{1}{1 - \dltwbun} \, 
	\biggl( \crosssup \, \dltwunn + \frac{\dltwun}{2}
		+ \frac{3(\crosssup + \dltwbun/2)^{2} \, \dltwu_{3}}{4(1 - \dltwbun)^{2}} 
	\biggr) 
	\, ,
\label{36gfijh94ejdvtwekoisedi}
	\\
	\dltwunss
	& \eqdef &
	2 \dltwu_{3} + \frac{\dltwunn}{2}  + \frac{2(\dltwuns + \dltwunn)}{(1 - \crosssup)^{2}} 
	\, .
\end{EQA}
Then \( \upsvd = \argmax_{\upsv} \fn(\upsv) \) satisfies
\begin{EQA}
	\| \DFN (\upsvd - \upsvs) \|_{\infty}
	\eqdef
	\max_{j \leq \dimp} | \DFN_{j} (\upsd_{j} - \upss_{j}) |
	& \leq &
	\rrinf \, .
\label{7ytdufchskmls7rghnk}
\end{EQA} 
Furthermore, 
\begin{EQA}
\label{usdhyw6hikhurnetrspp}
	\bigl\| \DFN^{-1} \bigl\{ \IFT (\upsvd - \upsvs) - \AAv \bigr\} \bigr\|_{\infty}
	& \leq &
	\dltwunss \, \| \DFN^{-1} \AAv \|_{\infty}^{2} 
	\, ,
	\\
	\bigl\| \DFN (\upsvd - \upsvs - \IFT^{-1} \AAv) \bigr\|_{\infty}
	& \leq &
	\frac{\dltwunss}{1 - \crosssup} \, \| \DFN^{-1} \AAv \|_{\infty}^{2} 
	\, ,
	\qquad
\label{usdhyw6hikhurnetrspp1}
\end{EQA}
and with \( \Delta \eqdef \Id_{\dimp} - \DFN^{-1} \IFT \, \DFN^{-1} \)
\begin{EQ}[rcl]
	\bigl\| \DFN (\upsvd - \upsvs) - \DFN^{-1} \AAv) \bigr\|_{\infty}
	& \leq &
	\frac{\dltwunss}{1 - \crosssup} \, \| \DFN^{-1} \AAv \|_{\infty}^{2} + \frac{\crosssup}{1 - \crosssup} \, \| \DFN^{-1} \AAv \|_{\infty}
	\, ,
	\\
	\bigl\| \DFN (\upsvd - \upsvs) - (\Id_{\dimp} + \Delta) \DFN^{-1} \AAv) \bigr\|_{\infty}
	& \leq &
	\frac{\dltwunss}{1 - \crosssup} \, \| \DFN^{-1} \AAv \|_{\infty}^{2} + \frac{\crosssup^{2}}{1 - \crosssup} \, \| \DFN^{-1} \AAv \|_{\infty}
	\, .
	\qquad
\label{usdhyw6hikhurnetrspp2}
\end{EQ}
\end{proposition}

\begin{remark}
To gain an intuition about the result of the proposition, consider a typical situation with 
\( \dltwun \leq \dltwu_{3} \),
\( \dltwunn \leq \dltwu_{3} \),
\( \dltwbun \leq 1/4 \),
\( 1 - \crosssup \leq 1/\sqrt{2} \).
Then \( \dltwuns \) from \eqref{36gfijh94ejdvtwekoisedi} fulfills \( \dltwuns \leq 1.37 \), and it holds
\( \dltwunss \leq 12 \dltwu_{3} \).
\end{remark}

\Section{Sup-norm expansions for a separable perturbation}
A linear perturbation is a special case of a separable perturbation of the form 
\( \pert(\upsv) = \sum_{j} \pert_{j}(\ups_{j}) \).
Another popular example is a Tikhonov regularization with \( \pert_{j}(\ups_{j}) = \lambda \ups_{j}^{2} \) or, 
more generally, a Sobolev-type penalty \( \pert_{j}(\ups_{j}) = \gp_{j}^{2} \ups_{j}^{2} \) for a given sequence 
\( \gp_{j}^{2} \).
Later we assume each component \( \pert_{j}(\ups_{j}) \) to be sufficiently smooth.
This does not allow to incorporate a sparse penalty \( \pert_{j}(\ups_{j}) = \lambda |\ups_{j}| \) or
complexity penalty \( \pert_{j}(\ups_{j}) = \lambda \Ind(\ups_{j} \neq 0) \).

Let \( \fs(\upsv) \) be concave, \( \upsvs = \argmax_{\upsv} \fs(\upsv) \), and
\( \IFT = - \nabla^{2} \fs(\upsvs) = (\IFT_{jm}) \).
Define a separable perturbation \( \fn(\upsv) = \fs(\upsv) - \pert(\upsv) \) of \( \fs(\upsv) \) by smooth functions \( \pert_{j}(\ups_{j}) \):
\begin{EQA}
	\fn(\upsv)
	&=&
	\fs(\upsv) - \pert(\upsv)
	=
	\fs(\upsv) - \sum_{j} \pert_{j}(\ups_{j}), 
\label{c7c7ctefwct5enfvyewsd}
\end{EQA} 
and let \( \upsvn = \argmax_{\upsv} \fn(\upsv) \).
We intend to bound the corresponding change \( \upsvn - \upsvs \) in a sup-norm.
Separability of the penalty is very useful, the cross derivatives of \( \fn \) are the same as for \( \fs \).
We only update the tensor \( \DFN(\upsv) \) assuming concavity of each \( \pert_{j}(\cdot) \):
\begin{EQA}
	\DFN(\upsv)
	& \eqdef &
	\diag\{ \DFN_{1}(\upsv),\ldots,\DFN_{\dimp}(\upsv) \},
	\qquad
	\DFN_{j}^{2}(\upsv)
	\eqdef
	\IF_{jj}(\upsv) - \pert''_{j}(\ups_{j}) \, .
\label{7fyc5wgdf7vrewgdhfywh}
\end{EQA}

\begin{proposition}
\label{Psemibiassupp}
Let \( \fs(\upsv) \) be concave function and 
\( \fn(\upsv) \) be a separable perturbation of \( \fs(\upsv) \); see \eqref{c7c7ctefwct5enfvyewsd}. 
For \( \upsvn = \argmax_{\upsv} \fn(\upsv) \), define the diagonal metric tensor \( \DFN = \DFN(\upsvn) \); 
see \eqref{7fyc5wgdf7vrewgdhfywh}.
Let \( \crosssup \) from \eqref{7tdsyf8iuwopkrtg4576} fulfill \( \crosssup < 1 \). 
Fix \( \Avm \eqdef \nabla \pert(\upsvs) = (\pert'_{j}(\upss_{j})) \) and 
\begin{EQA}
	\rrinf 
	&=& 
	\frac{\sqrt{2}}{1 - \crosssup} \| \DFN^{-1} \Avm \|_{\infty} 
	\, .
\label{8drywjwuycrt5sdeteyhfuwkp}
\end{EQA}
Let \nameref{LLpsupref} hold for \( \fn(\cdot) \) with this \( \rrinf \) and
\( \dltwun, \dltwunn, \dltwu_{3} \) satisfying  
\begin{EQA}[c]
	\dltwbun = \dltwunn \, \rrinf \leq 1/4 ,
	\qquad
	\dltwun \, \rrinf \leq 1/4 \, ,
	\qquad	
	\dltwunss \| \DFN^{-1} \Avm \|_{\infty} 
	\leq 
	\sqrt{2} - 1
	\, ,
\label{dhfiejfowelocuyehbrfp}
\end{EQA}
where
\begin{EQA}[rcl]
	\dltwuns
	& \eqdef &
	\frac{1}{1 - \dltwbun} \, 
	\biggl( \crosssup \, \dltwunn + \frac{\dltwun}{2}
		+ \frac{3(\crosssup + \dltwbun/2)^{2} \, \dltwu_{3}}{4(1 - \dltwbun)^{2}} 
	\biggr) 
	\, ,
\label{36gfijh94ejdvtwekoisedip}
	\\
	\dltwunss
	& \eqdef &
	2 \dltwu_{3} + \frac{\dltwunn}{2}  + \frac{2(\dltwuns + \dltwunn)}{(1 - \crosssup)^{2}} 
	\, .
\end{EQA}
Then \( \upsvd = \argmax_{\upsv} \fn(\upsv) \) satisfies
\begin{EQA}
	\| \DFN (\upsvd - \upsvs) \|_{\infty}
	\eqdef
	\max_{j \leq \dimp} | \DFN_{j} (\upsd_{j} - \upss_{j}) |
	& \leq &
	\rrinf \, .
\label{7ytdufchskmls7rghnkp}
\end{EQA} 
Furthermore, 
\begin{EQA}
\label{usdhyw6hikhurnetrsppp}
	\bigl\| \DFN^{-1} \bigl\{ \IFT (\upsvd - \upsvs) - \Avm \bigr\} \bigr\|_{\infty}
	& \leq &
	\dltwunss \, \| \DFN^{-1} \Avm \|_{\infty}^{2} 
	\, ,
	\\
	\bigl\| \DFN (\upsvd - \upsvs - \IFT^{-1} \Avm) \bigr\|_{\infty}
	& \leq &
	\frac{\dltwunss}{1 - \crosssup} \, \| \DFN^{-1} \Avm \|_{\infty}^{2} 
	\, ,
	\qquad
\label{usdhyw6hikhurnetrspp1p}
\end{EQA}
and with \( \Delta \eqdef \Id_{\dimp} - \DFN^{-1} \IFT \, \DFN^{-1} \)
\begin{EQ}[rcl]
	\bigl\| \DFN (\upsvd - \upsvs) - \DFN^{-1} \Avm) \bigr\|_{\infty}
	& \leq &
	\frac{\dltwunss}{1 - \crosssup} \, \| \DFN^{-1} \Avm \|_{\infty}^{2} + \frac{\crosssup}{1 - \crosssup} \, \| \DFN^{-1} \Avm \|_{\infty}
	\, ,
	\\
	\bigl\| \DFN (\upsvd - \upsvs) - (\Id_{\dimp} + \Delta) \DFN^{-1} \Avm) \bigr\|_{\infty}
	& \leq &
	\frac{\dltwunss}{1 - \crosssup} \, \| \DFN^{-1} \Avm \|_{\infty}^{2} + \frac{\crosssup^{2}}{1 - \crosssup} \, \| \DFN^{-1} \Avm \|_{\infty}
	\, .
	\qquad
\label{usdhyw6hikhurnetrspp2p}
\end{EQ}
\end{proposition}

\begin{remark}
\label{RLLsup}
If each \( \pert_{j}(\ups_{j}) \) is quadratic, then condition \nameref{LLpsupref} for \( \fn(\upsv) \) follows from 
the same condition for \( \fs(\upsv) \) as the third derivatives of these two functions coincide.
For general smooth perturbations \( \pert_{j} \), \eqref{d6f53ye5vry4fddfgsup} should be checked
for \( \nabla^{3}_{\ups_{j}\ups_{j} \ups_{j}} \fn(\upsv) = \nabla^{3}_{\ups_{j}\ups_{j} \ups_{j}} \fs(\upsv) + \pert'''_{j}(\ups_{j}) \).
The other conditions in \nameref{LLpsupref} can be checked for \( \fs \).
\end{remark}

%% file: BTL_short.tex

\def\alpG{p}
\def\dimGrV{n}
\Chapter{Estimation for the Bradley-Terry-Luce model}
\label{ScoBTL}

Let \( \Graph = (\GrV, \GrE) \) stand for a {comparison graph}, where the
vertex set \( \GrV = \{ 1,2,\ldots,\dimGrV \} \) represents the \( \dimGrV \) items of interest. 
The items \( j \) and \( m \) are {compared} if and only if \( (jm) \in \GrE \).
One observes independent paired comparisons \( Y^{(\ell)}_{jm} \), \( \ell = 1,\ldots,\nbin_{jm} \),
and \( Y^{(\ell)}_{jm} = 1 - Y^{(\ell)}_{mj} \).
For modeling and risk analysis, {Bradley-Terry-Luce} (BTL) model is frequently used; see \cite{BT1952}, \cite{Lu1959}.
The chance of each item winning a paired comparison is determined by the {relative scores} 
\begin{EQA}
	\P\bigl( \text{item \( j \) is preferred over item } m \bigr)
	&=&
	\P\bigl( Y_{jm}^{(\ell)} = 1 \bigr)
	=
	\frac{\ex^{\upss_{j}}}{\ex^{\upss_{j}} + \ex^{\upss_{m}}} 
	=
	\frac{1}{1 + \ex^{\upss_{m} - \upss_{j}}} \, .
\label{Pjoieujeui}
\end{EQA}
%
%
%
%
%
%
The goal is to recover the {score vector} 
\( \upsv = (\ups_{1},\ldots,\ups_{\dimGrV})^{\T} \) and top--\( k \) items.
%
%
%
%
%
%
%
%
%
%
%
%
Most of the theoretical results for a BTL model have been established under the following assumptions:
\begin{myitem}
	\item \( \Graph \) is a random Erd\H{o}s-R\'{e}nyi graph with the edge probability \( \alpG \); 
		\( \alpG \geq \CONST \dimGrV^{-1} \log \dimGrV \) ensures with overwhelming probability a connected graph;
	\item the values \( \nbin_{jm} \) are all the same; \( \nbin_{jm} \equiv L \), 
	for all \( (j,m) \in \GrE \);
	\item for the ordered sequence \( \upss_{(1)} \geq \upss_{(2)} \geq \ldots \geq \upss_{(\dimGrV)} \), it holds \( \upss_{(k)} - \upss_{(k+1)} > \Delta \);
	\item \( \upss_{(1)} - \upss_{(\dimGrV)} \leq \Range \);
\end{myitem}
see e.g. \cite{chen2020partial,GSZ2023} for a related discussion.
Under such assumptions, 
\cite{CFMW2019} and \cite{chen2020partial} showed that the conditions
\begin{EQA}
	\Delta^{2}
	& \geq &
	\CONST \frac{\log \dimGrV}{\dimGrV \, \alpG \, L}
\end{EQA}
enables to identify the top-\( k \) set with a high probability.
Both {regularized MLE} and a {spectral method} are rate optimal. 
We refer to \cite{GSZ2023} for an extensive overview and recent results for the BTL model including a non-asymptotic MLE expansion.

Unfortunately, some of the mentioned assumptions could be very restrictive in practical applications.
This especially concerns the graph structure and design of the set of comparisons.
An assumption of a bounded dynamic range \( \Range \) is very useful for the theoretical study because it allows to bound the success probability of each game away from zero and one.
However, it seems to be questionable for many real-life applications. 
Our aim is to demonstrate that the general approach of the paper enables us to get accurate results applying under
\begin{myitem}
	\item arbitrary configuration of the graph \( \Graph \); 
	\item {heterogeneous} numbers \( \nbin_{jm} \) of comparisons par edge;
	\item {unbounded} range \( \upss_{(1)} - \upss_{(\dimGrV)} \).

\end{myitem}
We still assume that \( \Graph \) is connected; otherwise, 
each connected component should be considered separately.

For \( j < m \),
denote \( S_{jm} = \sum_{\ell=1}^{\nbin_{jm}} Y_{jm}^{(\ell)} \) and \( S_{jm} = 0 \) if \( \nbin_{jm} = 0 \).
With \( \cdens(\ups) = \log(1 + \ex^{\ups}) \), the log-likelihood for the parameter vector \( \upsvs \) reads as follows:
\begin{EQA}
	L(\upsv)
	&=&
	\sum_{m=1}^{\dimGrV} \sum_{j =1}^{m-1} 
	\bigl\{ (\ups_{j} - \ups_{m}) S_{jm} - \nbin_{jm} \cdens(\ups_{j} - \ups_{m}) \bigr\} ,
	\qquad
\label{LusijEuiuj}
\end{EQA}
leading to the MLE
\begin{EQA}
	\tilde{\upsv}
	&=&
	\argmax_{\upsv} L(\upsv).
\end{EQA}
The function \( \cdens(\ups) = \log(1 + \ex^{\ups}) \) is convex, hence, \( L(\upsv) \) is {concave}.
However, representation \eqref{LusijEuiuj} reveals \emph{lack-of-identifiability} problem: 
\( \tilde{\upsv} \) is {not unique}, any shift 
\( \upsv \to \upsv + a \ev \)
does not change \( L(\upsv) \), \( \ev = \dimGrV^{-1/2} (1,\ldots,1)^{\T} \in \R^{\dimGrV} \).
Therefore, the Fisher information matrix \( \IF(\upsv) = -\nabla^{2} L(\upsv) \) is not positive definite and
\( L(\upsv) \) is not strongly concave, thus, \( \tilde{\upsv} \) is not uniquely defined.
For a connected graph \( \Graph \), assumed later on,
this issue can be resolved by fixing one component of \( \upsv \), e.g. \( \ups_{1} = 0 \), or by the condition 
\( \sum_{j} \ups_{j} = 0 \).
In general, we need {one condition} per {connected component} of the graph \( \Graph \).
Alternatively, one can use penalization
with a {quadratic} penalty \( \| \GP \upsv \|^{2}/2 \).
A ``non-informative'' choice is \( \GP^{2} = \gp^{2} \Id_{\dimGrV} \); cf. \cite{CFMW2019}.
Another option is to replace the constraint \( \sum_{j} \ups_{j} = 0 \) by the penalty 
\( \| \GP \upsv \|^{2} = \gp^{2} \langle \upsv,\ev \rangle^{2} \)
with \( \GP^{2} = \gp^{2} \ev \, \ev^{\T} \).
It is obvious to see that this penalization replaces the condition \( \sum_{j} \ups_{j} = 0 \)
because it vanishes for any such \( \upsv \).
If the true skill vector \( \upsvs \) is defined under this condition then 
\( \upsvs_{\GP} = \argmax_{\upsv} \E \LGP(\upsv) = \upsvs \).
Hence, the penalty \( \| \GP \upsv \|^{2}/2 \) does not yield any bias of estimation.
The choice of the constant \( \gp^{2} \) is not essential, we just assume it 
sufficiently large.
One more benefit of using a penalty \( \| \GP \upsv \|^{2} \) is that it ensures 
the desired stochastically dominant structure of the Fisher information matrix:
each diagonal element is larger than the sum of the off-diagonal for this row. 
This property is important for proving a concentration of the penalized MLE 
\begin{EQA}
	\tilde{\upsv}_{\GP}
	&=&
	\argmax_{\upsv} \LGP(\upsv) 
	=
	\argmax_{\upsv} \bigl( L(\upsv) - \| \GP \upsv \|^{2}/2 \bigr) .
\label{ftbc43e3e3rfdf6ddfd}
\end{EQA}
%
%
The next lemma reveals some important features of the Fisher information matrix \( \IF(\upsv) \).

\begin{lemma}
\label{LIFBTL}
With \( \cdens''(\ups) = \frac{\ex^{\ups}}{(1 + \ex^{\ups})^{2}} \), 
the entries \( \IF_{jm}(\upsv) \) of \( \IF(\upsv) = -\nabla^{2} \E L(\upsv) \) 
satisfy
\begin{EQA}
\label{sfhce8fy8efhidveuuedf8}
	\IF_{jj}(\upsv)
	&=&
	\sum_{m \neq j} \nbin_{jm} \cdens''(\ups_{j} - \ups_{m}),
	\\
	\IF_{jm}(\upsv)
	&=&
	- \nbin_{jm} \cdens''(\ups_{j} - \ups_{m}) ,
	\quad
	j \neq m \, .
\end{EQA}
This yields \( \IF_{jj}(\upsv) = - \sum_{m \neq j} \IF_{jm}(\upsv) \).
Moreover, each eigenvalue of \( \IF(\upsv) \) 
belongs to \( [0,2\IF_{jj}(\upsv)] \) for some \( j \leq \dimGrV \).
\end{lemma}

\begin{proof}
The structure of \( \IF(\upsv) \) follows directly from \eqref{LusijEuiuj}.
Gershgorin's theorem implies the final statement. 
\end{proof}

The rigorous theoretical study including finite sample expansions 
\begin{EQ}[rcl]
\label{usdhyw6hikhurnetrspF}
	\| \DPN (\tilde{\upsv}_{\GP} - \upsvs - \IF_{\GP}^{-1} \nabla \zeta) \|_{\infty}
	& \leq &
	\frac{\dltwunss}{1 - \crosssup} \, \| \DPN^{-1} \nabla \zeta \|_{\infty}^{2} \, ,
	\\
	\| \DPN (\tilde{\upsv}_{\GP} - \upsvs - \DPN^{-2} \nabla \zeta) \|_{\infty}
	& \leq &
	\frac{\dltwunss}{1 - \crosssup} \| \DPN^{-1} \nabla \zeta \|_{\infty}^{2}
	+ \frac{\crosssup}{1 - \crosssup} \, \| \DPN^{-1} \nabla \zeta \|_{\infty}
	\, .
	\qquad
\label{usdhyw6hikhurnetrsppm}
\end{EQ}
and a check of all the related conditions 
is deferred to the forthcoming paper \cite{PankSp2025}.
Here we only present some numerical results illustrating the established theoretical guarantees.

We use the following experimental setup from \cite{GSZ2023}:
an Erd\"{o}s–R\'{e}nyi comparison graph with \( n \) vertices, edge probability \( \alpG \),
the number \( L \) of comparisons between pairs of items is homogeneous,
the values \( \upss_{i} \) are i.i.d  samples from  \( U[\ups_{\min}, \ups_{\max}] \), \( \Range = \ups_{\max} - \ups_{\min}  \).
For experiments, they used \( L = 1 \), \( \alpG = (\log(n)^3 / n) \), \( \ups_{\min} = 0 \), \( \ups_{\max} = 2 \), and 
\( n \in \{100, 200, 500, 1000, 2000\} \). 
The study presents the important values \( \crosssup \) from \eqref{7tdsyf8iuwopkrtg4576} and 
checks the quality of the expansion \eqref{usdhyw6hikhurnetrsppm}.

\begin{figure}[H]
    \centering
    \includegraphics[width=0.5\textwidth,height = 0.18\textheight]{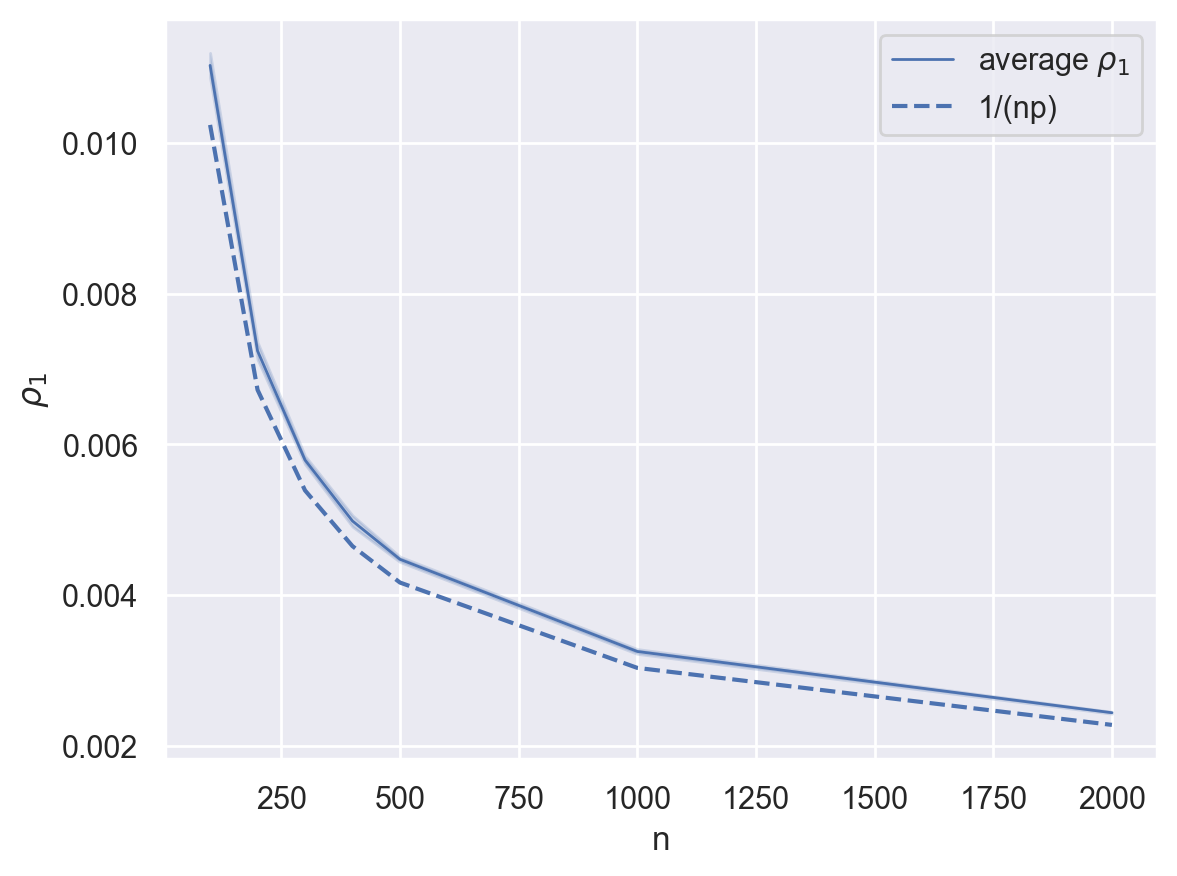}
    \caption{Average \( \crosssup \) value for different values of \( n \).}
    \label{fig:crosssup}
\end{figure}

Figure~\ref{fig:crosssup} depicts the value of \( \crosssup \) for each \( n \) 
after averaging over few realizations of an Erd\"{o}s–R\'{e}nyi comparison graph with \( L = 1 \), 
\( \alpG = \left(\frac{\log^{3}(n)}{n}\right) \), and \( \upss_{j} \) as i.i.d. samples from \( U[0, 2] \). 
The error bars represent three times the standard deviation.
The results support the hypothesis that \( \crosssup \) is small for Erd\"{o}s–R\'{e}nyi graphs and close to 
\( \frac{1}{n\alpG} \ll 1 \).

%
%

\begin{figure}[H]
    \centering
    \includegraphics[width=0.8\textwidth,height = 0.19\textheight]{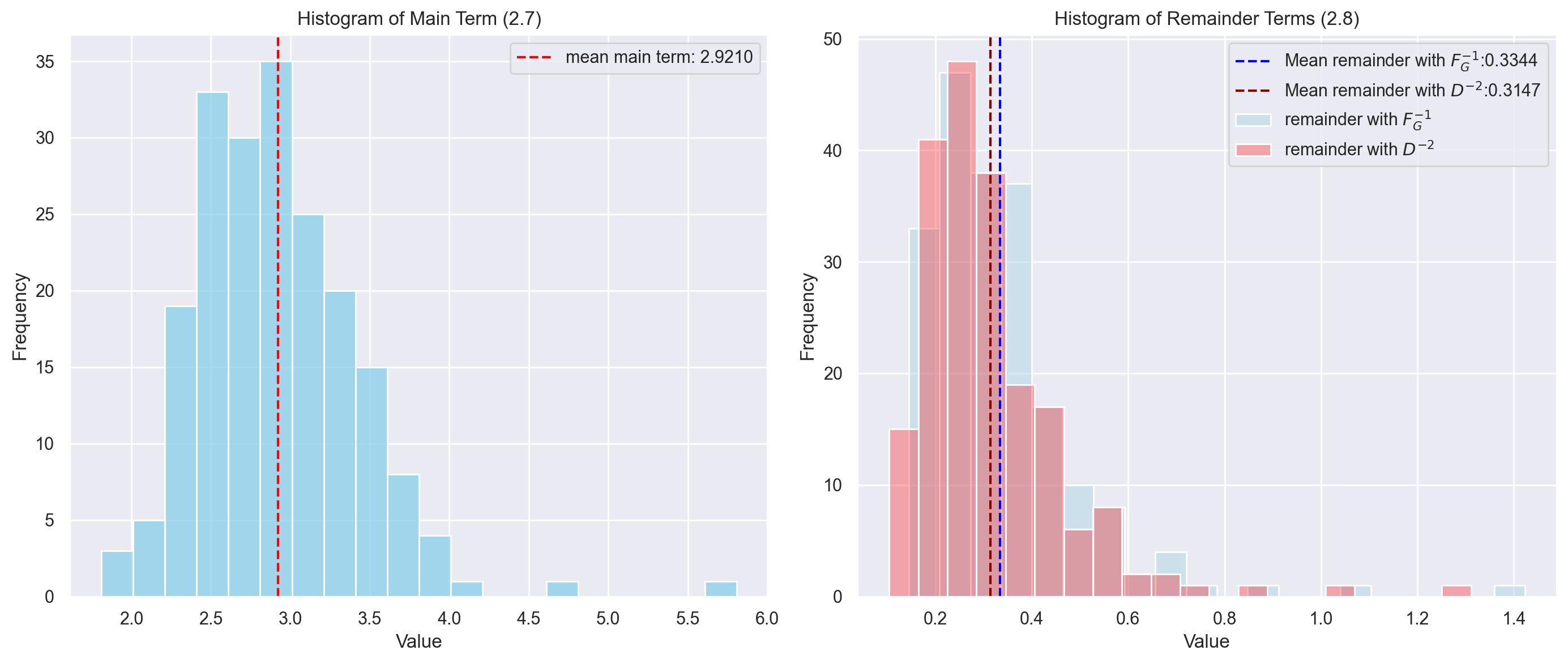}
    \caption{Distribution of the leading term and the remainder for \( n=100 \).}
    \label{fig:BTLapprox}
\end{figure}

Figure~\ref{fig:BTLapprox} left presents the norm of the leading term 
\( \| \IF_{\GP}^{-1} \nabla \zeta \|_{\infty} \) and \( \| \DPN^{-2} \nabla \zeta \|_{\infty} \) 
while the right plot shows the errors \( \| \tilde{\upsv}_{\GP} - \upsvs - \IF_{\GP}^{-1} \nabla \zeta \|_{\infty} \)
and \( \| \tilde{\upsv}_{\GP} - \upsvs - \DPN^{-2} \nabla \zeta \|_{\infty} \) for \( n=100 \).
One can see that the remainder is much smaller in the magnitude that the leading term, thus indicating 
a very good quality of the provided expansion.

\begin{figure}[H]
    \centering
    \includegraphics[width=0.8\textwidth,height = 0.2\textheight]{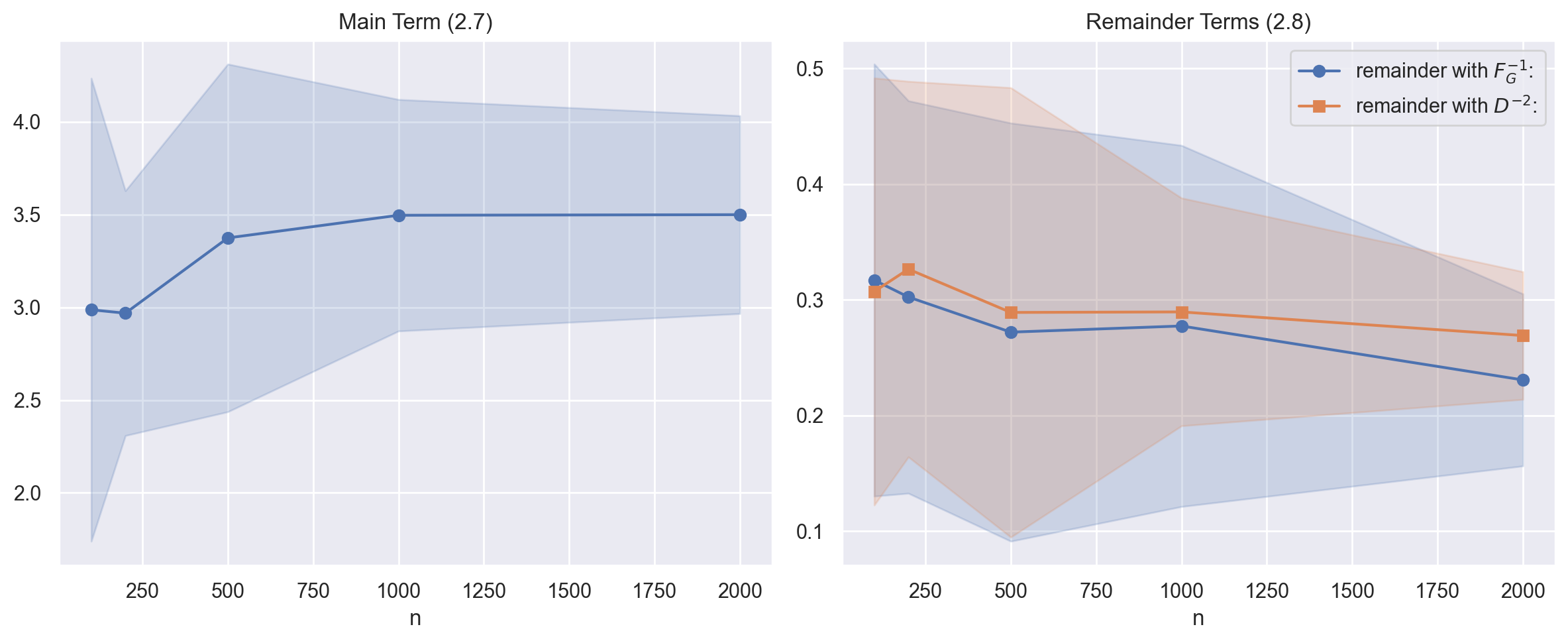}
    \caption{Comparison of the leading term and the remainder for different \( n \).}
    \label{fig:BTLapprox2}
\end{figure}

Similar results for \( n \in \{100, 200, 500, 1000, 2000\} \) are collected in Figure~\ref{fig:BTLapprox2}.

%% file: genMLE_short.tex

\Chapter{Properties of the MLE \( \tilde{\upsv} \) for SLS models}
\label{SgenBounds}
This \chname overviews general results about concentration and expansion of the MLE in the SLS setup from \cite{Sp2024}
which substantially improve the bounds from \cite{SP2013_rough} and \cite{SpLaplace2022}.
%
We assume to be given a random function \( L(\upsv) \), \( \upsv \in \Ups \subseteq \R^{\dimp} \),
\( \dimp < \infty \).
This function can be viewed as log-likelihood or negative loss.
Consider in parallel two optimization problems defining 
the MLE \( \tilde{\upsv} \) and its population counterpart (the background truth) \( \upsvs \):
\begin{EQA}[rcl]
	\tilde{\upsv} 
	&=& 
	\argmax_{\upsv} L(\upsv),
	\qquad
	\upsvs 
	=
	\argmax_{\upsv} \E L(\upsv),
	\qquad
\label{tuGauLGususGE}
\end{EQA}
Define the Fisher information matrix \( \IF(\upsv) \eqdef - \nabla^{2} \E L(\upsv) \) 
and denote \( \IF = \IF(\upsvs) \). 

\Section{Basic conditions}
\label{Scondgeneric}
Now we present our major conditions.
The most important one is about linearity of the stochastic component 
\( \zeta(\upsv) = L(\upsv) - \E L(\upsv) = L(\upsv) - \E L(\upsv) \).

\medskip
\begin{description}
    \item[\label{Eref} \( \bb{(\zeta)} \)]
      \textit{The stochastic component \( \zeta(\upsv) = L(\upsv) - \E L(\upsv) \) of the process \( L(\upsv) \) is linear in 
      \( \upsv \). 
      We denote by \( \nabla \zeta \equiv \nabla \zeta(\upsv) \in \R^{\dimp} \) its gradient
      }.
\end{description}

Below we assume some concentration properties of the stochastic vector \( \nabla \zeta \).
More precisely, we require that \( \nabla \zeta \) obeys the following condition.

\begin{description}
\item[\label{EU2ref}\( \bb{(\nabla \zeta)} \)]
	\textit{There exists \( \VP^{2} \geq \Var(\nabla \zeta) \) such that  
	for all considered  \( \BBH \in \Matr_{\dimp} \) and \( \xx > 0 \)
	}
\begin{EQA}
	\P\bigl( \| \BBH^{1/2} \VP^{-1} \nabla \zeta \| \geq \zq(\BBH,\xx) \bigr)
	& \leq &
	3 \ex^{-\xx} ,
\label{2emxGPm12nz122}
	\\
	\zq^{2}(\BBH,\xx)
	& \eqdef &
	\tr \BBH + 2 \sqrt{\xx \, \tr \BBH^{2}} + 2 \xx \| \BBH \| \,  .
\label{34rtyghuioiuyhgvftid}
\end{EQA}
\end{description}

This condition can be effectively checked if the errors in the data exhibit sub-gaussian or sub-exponential behavior; see 
\ifsqnorm{Section~\ref{SdevboundnonGauss}.}{\cite{Sp2023c}, \cite{Sp2023d}.}
The important special case corresponds to \( \BBH = \IF^{-1/2} \VP^{2} \IF^{-1/2} \) 
and \( \xx \approx \log n \) leading to the bound
\begin{EQA}
	\P\bigl( \| \IF^{-1/2} \nabla \zeta \| > \zq(\BBH,\xx) \bigr)
	& \leq &
	3/n .
\label{udyvfeyejff6777dj23}
\end{EQA}
The  value \( \dimG = \tr(\IF^{-1} \VP^{2}) \) can be called the \emph{effective dimension}; see \cite{SP2013_rough}.

We also assume that the log-likelihood \( L(\upsv) \) or, equivalently, its deterministic part 
\( \E L(\upsv) \) is a concave function.

\medskip
\begin{description}
    \item[\label{LLref} \( \bb{(\mathcal{C})} \)]
      \textit{The function \( \fs(\upsv) \eqdef \E L(\upsv) \) is concave on open and convex set \( \Ups \subset \R^{\dimp} \).      
      }
\end{description}
\medskip

Later we will also need some smoothness conditions on the function \( f(\upsv) = \E L(\upsv) \)
within a local vicinity of the point \( \upsvs \).
The notion of locality is given in terms of a metric tensor \( \DPN \in \Matr_{\dimp} \).
In most of the results later on, one can use \( \DPN = \IF^{1/2} \).
In general, we only assume \( \DPN^{2} \leq \dmax^{2} \IF \) for some \( \dmax > 0 \).
Let \( \fs(\upsv) \) be three or sometimes even four times Gateaux differentiable in \( \upsv \in \Ups \).
Introduce the following conditions.

\begin{description}
    \item[\label{LLsT3ref} \( \bb{(\mathcal{T}_{3}^{*})} \)]
    \emph{\( \fs(\upsv) \) is three times differentiable and 
	}
\begin{EQA}
    \sup_{\uv \colon \| \DFN(\upsv) \uv \| \leq \rr} \,\, \sup_{\zv \in \R^{\dimp}} \,\, 
    \frac{\bigl| \langle \nabla^{3} \fs(\upsv + \uv), \zv^{\otimes 3} \rangle \bigr|}
		 {\| \DFN(\upsv) \zv \|^{3}} 
	& \leq &
	\dltwu_{3} \, .
\label{jcxhydtferyu9j3d6vhew}
\end{EQA}

    \item[\label{LLsT4ref} \( \bb{(\mathcal{T}_{4}^{*})} \)]
    \emph{\( \fs(\upsv) \) is four times differentiable and 
	}
\begin{EQA}
    \sup_{\uv \colon \| \DFN(\upsv) \uv \| \leq \rr} \,\, \sup_{\zv \in \R^{\dimp}} \,\, 
    \frac{\bigl| \langle \nabla^{4} \fs(\upsv + \uv), \zv^{\otimes 4} \rangle \bigr|}
		 {\| \DFN(\upsv) \zv \|^{4}} 
	& \leq &
	\dltwu_{4} \, .
\label{jcxhydtferyu9j3d6vhew4}
\end{EQA}

\end{description}

%
\noindent
By Banach's characterization \cite{Banach1938}, \nameref{LLsT3ref} implies
\begin{EQA}
	\bigl| \langle \nabla^{3} \fs(\upsv + \uv), \zv_{1} \otimes \zv_{2} \otimes \zv_{3} \rangle \bigr|
	& \leq &	 
	\dltwu_{3} \| \DFN(\upsv) \zv_{1} \| \, \| \DFN(\upsv) \zv_{2} \| \, \| \DFN(\upsv) \zv_{3} \| \, 
\label{jbuyfg773jgion94euyyfg}
\end{EQA}
for any \( \uv \) with \( \| \DFN(\upsv) \uv \| \leq \rr \) and all \( \zv_{1} , \zv_{2}, \zv_{3} \in \R^{\dimp} \).
Similarly under \nameref{LLsT4ref}
\begin{EQA}
	\bigl| \langle \nabla^{4} \fs(\upsv + \uv), \zv_{1} \otimes \zv_{2} \otimes \zv_{3} \otimes \zv_{4} \rangle \bigr|
	& \leq &	 
	\dltwu_{4} \prod_{k=1}^{4} \| \DFN(\upsv) \zv_{k} \| \, ,
	\quad 
	\zv_{1} , \zv_{2}, \zv_{3}, \zv_{4} \in \R^{\dimp} \, .
	\qquad
\label{jbuyfg773jgion94euyyfg4}
\end{EQA}

The values \( \dltwu_{3} \) and \( \dltwu_{4} \) are usually very small.
Some quantitative bounds are given later in this section
under the assumption that the function \( \fs(\upsv) \) can be written in the form \( - \fs(\upsv) = n \hL(\upsv) \) 
for a fixed smooth function \( h(\upsv) \) with the Hessian \( \nabla^{2} \hL(\upsv) \). 
The factor \( n \) has meaning of the sample size%
\ifapp{; see \Chname \ref{ScritdimMLE} or \Chname \ref{SGBvM}.}{.}

\begin{description}
    \item[\label{LLtS3ref} \( \bb{(\mathcal{S}_{3}^{*})} \)]
      \emph{ \( - \fs(\upsv) = n \hL(\upsv) \) for \( \hL(\upsv) \) three times differentiable and
\begin{EQA}
	\sup_{\uv \colon \| \HL(\upsv) \uv \| \leq \rr/\sqrt{n}} 
	\frac{\bigl| \langle \nabla^{3} \hL(\upsv + \uv), \uv^{\otimes 3} \rangle \bigr|}{\| \HL(\upsv) \uv \|^{3}}
	& \leq &
	\hmax_{3} \, .
\end{EQA}
}
    \item[\label{LLtS4ref} \( \bb{(\mathcal{S}_{4}^{*})} \)]
      \emph{ the function \( \hL(\cdot) \) satisfies \nameref{LLtS3ref} and  
\begin{EQA}
	\sup_{\uv \colon \| \HL(\upsv) \uv \| \leq \rr/\sqrt{n}}
	\frac{\bigl| \langle \nabla^{4} \hL(\upsv + \uv), \uv^{\otimes 4} \rangle \bigr|}{\| \HL(\upsv) \uv \|^{4}}
	& \leq &
	\hmax_{4} \, .
\end{EQA}
}
\end{description}

\noindent
\nameref{LLtS3ref} and \nameref{LLtS4ref}
are local versions of the so-called self-concordance condition; see \cite{Ne1988} and \cite{OsBa2021}.
In fact, they require that each univariate function \( \hL(\upsv + t \uv) \), \( t \in \R \),
is self-concordant with some universal constants \( \hmax_{3} \) and \( \hmax_{4} \).

The class of models satisfying the conditions \nameref{Eref}, 
\nameref{EU2ref}, and \nameref{LLref}
with a smooth function \( f(\upsv) = \E L(\upsv) \) will be referred to as \emph{stochastically linear smooth} (SLS). 
This class includes linear regression, generalized linear models (GLM), and log-density models; 
see \cite{SpPa2019}, \cite{OsBa2021}\ifapp{ or Section~\ref{SGBvM} later.}{ or \cite{SpLaplace2022}.}
However, this class is much larger.
For instance, nonlinear regression can be adapted to the SLS framework 
by an extension of the parameter space; see 
\ifNL{Section~\ref{Snoninverse}.}{\cite{Sp2024}.}

\Section{Expansions and risk bounds under third-order smoothness}
\label{SFiWiexs3}

Define
\begin{EQA}[c]
	\BBH_{\DPN} = \DPN \IF^{-1} \VP^{2} \IF^{-1} \DPN,
	\\
	\dimD
	\eqdef
	\tr \BBH_{\DPN} \, ,
	\quad
	\rrD \eqdef \zq(\BBH_{\DPN},\xx)
	\leq 
	\sqrt{\tr \BBH_{\DPN}} + \sqrt{2\xx \, \| \BBH_{\DPN} \|} \, ;
\label{y7s7d7d77dfdy7fuegue3j}
\end{EQA}
cf. \eqref{34rtyghuioiuyhgvftid}.
By \nameref{EU2ref}, it holds 
\( \P(\| \DPN \, \IF^{-1} \nabla \zeta \| > \rrD) \leq 3\ex^{-\xx} \).
The result follows by limiting to the set \( \Omega(\xx) \) on which \( \| \DPN \, \IF^{-1} \nabla \zeta \| \leq \rrD \)
and by applying Proposition~\ref{PFiWigeneric2}.

\begin{theorem}
\label{TFiWititG3}
Assume 
\nameref{Eref},
\nameref{EU2ref},
and 
\nameref{LLref}.
Let \nameref{LLsT3ref} hold at \( \upsvs \) with a metric tensor 
\( \DPN \) and values \( \rr \) and \( \dltwu_{3} \) 
satisfying
\begin{EQA}[c]
	\DPN^{2} \leq \dmax^{2} \hspm \IF ,
	\quad
	\rr \geq \frac{3}{2} \, \rrD \, ,
	\quad
	\dltwu_{3} \, \dmax^{2} \hspm \rrD < \frac{4}{9} \, ,
\label{8difiyfc54wrboeMLE}
\end{EQA}
where \( \rrD \) is from \eqref{y7s7d7d77dfdy7fuegue3j}.
With \( \Omega(\xx) = \{ \| \DPN \, \IF^{-1} \nabla \zeta \| \leq \rrD \} \), 
it holds \( \P(\Omega(\xx)) \geq 1 - 3 \ex^{-\xx} \) and the following statements hold:

\begin{description}
	\item[1. Concentration:] On \( \Omega(\xx) \),
\begin{EQA}
	\| \IF^{1/2} (\tilde{\upsv} - \upsvs) \|  
	& \leq &
	\frac{4}{3} \rrB \, ,
	\qquad
	 \| \DPN (\tilde{\upsv} - \upsvs) \|  
	\leq 
	\frac{4\dmax}{3} \, \rrB \, ,
\label{rhDGtuGmusGU0a2MLE}
\end{EQA}

	\item[2. Fisher expansion:] on \( \Omega(\xx) \)
\begin{EQA}[rcl]
    \| \DPN^{-1} \IF (\tilde{\upsv} - \upsvs - \IF^{-1} \nabla \zeta) \|
    & \leq &
    \frac{3\dltwu_{3}}{4} \| \DPN \, \IF^{-1} \nabla \zeta \|^{2} 
	\, .
\label{DGttGtsGDGm13rG22}
\end{EQA}

	\item[3. Wilks expansion:] on \( \Omega(\xx) \)
\begin{EQ}[rcl]
    \Bigl| 2 L(\tilde{\upsv}) - 2 L(\upsvs) - \| \IF^{-1/2} \nabla \zeta \|^{2} \Bigr|
    & \leq &
    \frac{\dltwu_{3}}{2} \, \| \DPN \IF^{-1} \nabla \zeta \|^{3} 
    \, .
\label{3d3Af12DGttG2}
\end{EQ}

	\item[4. PAC loss bounds:] on \( \Omega(\xx) \), for any linear mapping \( \QP \colon \R^{\dimp} \to \R^{\dimq} \), 
\begin{EQA}
	\| \QP (\tilde{\upsv} - \upsvs - \IF^{-1} \nabla \zeta) \|
	& \leq &
	\| \QP \IF^{-1} \DPN \| \, \frac{3\dltwu_{3}}{4} \, \| \DPN \IF^{-1} \nabla \zeta \|^{2}
	\, .
	\qquad
	\quad
\label{g25re9fjfregdndg}
\end{EQA}

	\item[5. \( \ell_{2} \) risk bounds]
Introduce 
\begin{EQA}[c]
	\riskt_{\QP} \eqdef \E \{ \| \QP \IF^{-1} \nabla \zeta \|^{2} \Ind_{\Omega(\xx)} \} 
	\leq 
	\dimQ  \, 
\label{7djhed8cjfct534etgdhdyQP}
\end{EQA}
with \( \dimQ \eqdef \E \| \QP \IF^{-1} \nabla \zeta \|^{2} = \tr \Var(\QP \IF^{-1} \nabla \zeta) \).
Then
\begin{EQA}
	\E \bigl\{ \| \QP (\tilde{\upsv} - \upsvs) \| \Ind_{\Omega(\xx)} \bigr\} 
	& \leq &
	\riskt_{\QP}^{1/2} 
	+ \| \QP \IF^{-1} \DPN \| \, \frac{3\dltwu_{3}}{4} \, \dimD \, .
\label{EtuGus11md3GQP}
\end{EQA}
	\item[6. Quadratic risk bounds]
Let \( \E \bigl\{ \| \DPN \IF^{-1} \nabla \zeta \|^{4} \Ind_{\Omega(\xx)} \bigr\} \leq \CONSTi_{4}^{2} \, \dimD^{2} \). 
Define
\begin{EQ}[rcl]
	\alp_{\QP}
	& \eqdef & 
	\frac{\| \QP \IF^{-1} \DPN \| \, (3/4) \dltwu_{3} \, \CONSTi_{4} \, \dimD } {\sqrt{\riskt_{\QP}}}
	\, . 
\label{6dhx6whcuydsds655srew}
\end{EQ}
If \( \alp_{\QP} < 1 \) then
\begin{EQA}
	(1 - \alp_{\QP})^{2} \riskt_{\QP} 
	\leq 
	\E \bigl\{ \| \QP \, (\tilde{\upsv} - \upsvs) \|^{2} \Ind_{\Omega(\xx)} \bigr\}
	& \leq &
	(1 + \alp_{\QP})^{2} \riskt_{\QP} \, .
\label{EQtuGmstrVEQtGQP}
\end{EQA}

\end{description}
\end{theorem}

\Section{Effective and critical dimension in ML estimation}
\label{ScritdimMLE}
This section discusses the important question of the critical parameter dimension 
still ensuring the validity of the presented results.
To be more specific, we only consider the 3S-results of Theorem~\ref{TFiWititG3}.
Also, assume \( \dmax \equiv 1 \).
The important constant \( \dltwu_{3} \) is identified by \nameref{LLtS3ref}: \( \dltwu_{3} = \hmax_{3} /\sqrt{n} \),
where the scaling factor \( n \) means the sample size.
It can be defined as the smallest eigenvalue of the Fisher operator \( \IF \).

First, we discuss the case \( \QP = \DPN = \IF^{1/2} \).
It appears that in this full dimensional situation, all the obtained results apply and are meaningful under the condition \( \dimA \ll n \),
where \( \dimA = \tr(\BBH) \) for \( \BBH = \IF^{-1/2} \VP^{2} \IF^{-1/2} \) is the \emph{effective dimension} of the problem.
Indeed,  \( \rrD^{2} = \rrB^{2} \approx \tr(\BBH) = \dimA \), and
condition \eqref{8difiyfc54wrboeMLE} requires \( \dltwu_{3} \, \rrD \ll 1 \) which can be spelled out as \( \dimA \ll n \).
Expansion \eqref{DGttGtsGDGm13rG22} means
\begin{EQA}
	\| \IF^{1/2} (\tilde{\upsv} - \upsvs) \|
    & \leq &
    \| \IF^{-1/2} \nabla \zeta \| + \frac{3\dltwu_{3}}{4} \| \IF^{-1/2} \nabla \zeta \|^{2} \, ,
\label{Y6DF66FCC6dseeEGHBE}
\end{EQA}
and the second term on the right-hand side of this bound is smaller 
than the first one under the same condition \( \dltwu_{3} \, \rrD \ll 1 \).
Similar observations apply to bound \eqref{EQtuGmstrVEQtGQP} of Theorem~\ref{TFiWititG3} 
which is meaningful only if \( \alp_{\QP} \) in \eqref{6dhx6whcuydsds655srew} is small.
As \( \riskt_{\QP} \approx \dimQ = \dimA \), the condition \( \dltwu_{3} \, \rrD \ll 1 \) implies \( \alp_{\QP} \ll 1 \)
and hence, the bound \eqref{EQtuGmstrVEQtGQP} is sharp.
We conclude that the main properties of the MLE \( \tilde{\upsv} \) 
are valid under the condition \( \dimA \ll n \) meaning sufficiently many observations 
per effective number of parameters.

\ifNL{The situation changes drastically if \( \QP \) is not full-dimensional as e.g. in semiparametric estimation,
when \( \QP \) projects onto a low-dimensional target component.
We will see in \Chname \ref{SsemiMLE} that in this case,
\( \alp_{\QP} \ll 1 \) requires \( \dimA^{2} \ll n \).}{}
%

%

%% file: localbounds_short.tex

\def\AFN{\mathbbmsl{U}}
\def\Avm{\bb{M}}

\Chapter{Optimization after linear perturbation. A basic lemma}
\label{Squadnquad}
This section presents a basic lemma from \cite{Sp2024}.
Let \( \fs(\upsv) \) be a smooth concave function, 
\begin{EQA}
	\upsvs
	&=&
	\argmax_{\upsv} \fs(\upsv),
\label{fg5hg3gf98tkj3dciryt}
\end{EQA}
and \( \IFN = - \nabla^{2} \fs(\upsvs) \).
Later we study the question of how the point of maximum and the value of maximum of \( \fs \) change if we add a linear or quadratic 
component to \( \fs \).
More precisely, let another function \( \fn(\upsv) \) satisfy for some vector \( \Av \)
\begin{EQA}
	\fn(\upsv) - \fn(\upsvs) 
	&=&
	\bigl\langle \upsv - \upsvs, \Av \bigr\rangle + \fs(\upsv) - \fs(\upsvs) .
\label{4hbh8njoelvt6jwgf09}
\end{EQA}
A typical example corresponds to \( \fs(\upsv) = \E L(\upsv) \) and \( \fn(\upsv) = L(\upsv) \) 
for a random function \( L(\upsv) \) with a linear stochastic component \( \zeta(\upsv) = L(\upsv) - \E L(\upsv) \)%
\ifapp{; see \nameref{Eref}.}{.}
Then \eqref{4hbh8njoelvt6jwgf09} is satisfied with \( \Av =	\nabla \zeta \).
Define
\begin{EQA}
	\upsvn
	& \eqdef &
	\argmax_{\upsv} \fn(\upsv),
	\qquad
	\fn(\upsvn)
	=
	\max_{\upsv} \fn(\upsv) .
\label{6yc63yhudf7fdy6edgehy} 
\end{EQA}
The aim of the analysis is to evaluate the quantities \( \upsvn - \upsvs \) and
\( \fn(\upsvn) - \fn(\upsvs) \).

\begin{proposition}
\label{PFiWigeneric2}
Let \( \fs(\upsv) \) be a strongly concave function with \( \fs(\upsvs) = \max_{\upsv} \fs(\upsv) \)  
and \( \IFN = - \nabla^{2} \fs(\upsvs) \).
Assume \nameref{LLsT3ref} at \( \upsvs \) with \( \DFN^{2} \), \( \rrn \), and \( \dltwu_{3} \) such that 
\begin{EQA}[c]
	\DFN^{2} \leq \dmax^{2} \, \IFN ,
	\quad
	\rrn \geq \frac{3}{2} \| \DFN \, \IFN^{-1} \Av \| \, ,
	\quad
	\dmax^{2} \dltwu_{3} \| \DFN \, \IFN^{-1} \Av \| < \frac{4}{9} \, .
\label{8difiyfc54wrboer7bjfr}
\end{EQA}
Then \( \| \DFN (\upsvn - \upsvs) \| \leq (3/2) \| \DFN \, \IFN^{-1} \Av \| \) and moreover,
\begin{EQA}[rcl]
    \| \DFN^{-1} \IFN (\upsvn - \upsvs - \IFN^{-1} \Av) \|
    & \leq &
    \frac{3\dltwu_{3}}{4} \| \DFN \, \IFN^{-1} \Av \|^{2} 
	\, .
\label{DGttGtsGDGm13rGa2}
\end{EQA}
Moreover, 
\begin{EQA}
    \Bigl| 2 \fn(\upsvn) - 2 \fn(\upsvs) - \| \IFN^{-1/2} \Av \|^{2} \Bigr|
    & \leq &
    \frac{\dltwu_{3}}{2} \, \| \DFN \, \IFN^{-1} \Av \|^{3} \, .
    \qquad
\label{3d3Af12DGttGa2}
\end{EQA}
\end{proposition}

\begin{remark}
\label{Rbiasgeneric}
The roles of \( \fs \) and \( \fn \) can be exchanged.
In particular, \eqref{DGttGtsGDGm13rGa2} applies with \( \IFN = \IFN(\upsvn) \) provided that 
\nameref{LLsT3ref} is also fulfilled at \( \upsvn \).
\end{remark}

%% file: semi_opt_tools.tex

\begin{proof}[Proof of Proposition~\ref{PbiassemiN}]
Define \( \fn_{\nuiv}(\targv) = \fs_{\nuiv}(\targv) - \langle \targv,\Av_{\nuiv} \rangle \).
Then \( \fn_{\nuiv} \) is concave and \( \nabla \fn_{\nuiv}(\targvs) = 0 \) yielding 
\( \targvs = \argmax_{\targv} \fn_{\nuiv}(\targv) \).
Now Proposition~\ref{PFiWigeneric2} yields \eqref{jhcvu7ejdytur39e9frtfw}.
Also by \eqref{3d3Af12DGttGa2}
\begin{EQA}
	\Bigl| 2 \fn_{\nuiv}(\targvs) - 2 \fn_{\nuiv}(\targv_{\nuiv}) - \| \IFN_{\nuiv}^{-1/2} \Av_{\nuiv} \|^{2} \Bigr|
	& \leq &
	\frac{\dltwu_{3}}{2} \, \| \DFN_{\nuiv} \, \IFN_{\nuiv}^{-1} \Av_{\nuiv} \|^{3} 
\label{ytduw237832jfg656335hdi}
\end{EQA}
yielding
\begin{EQA}
	\bigl| 
		2 \fs_{\nuiv}(\targvs) - 2 \fs_{\nuiv}(\targv_{\nuiv}) - 2 \langle \targvs - \targv_{\nuiv} , \Av_{\nuiv} \rangle 
		- \| \IFN_{\nuiv}^{-1/2} \Av_{\nuiv} \|^{2} 
	\bigr|
	& \leq &
	\frac{\dltwu_{3}}{2} \, \| \DFN_{\nuiv} \, \IFN_{\nuiv}^{-1} \Av_{\nuiv} \|^{3} \, .
	\qquad
\label{gtxdbcbcgdfswerfdsdew}
\end{EQA}
By \eqref{jhcvu7ejdytur39e9frtfw}
\begin{EQA}
	&& \nquad
	\bigl| \langle \targvs - \targv_{\nuiv} , \Av_{\nuiv} \rangle + \| \IFN_{\nuiv}^{-1/2} \Av_{\nuiv} \|^{2} \bigr|
	=
	\bigl| \langle \IFN_{\nuiv}^{1/2} (\targvs - \targv_{\nuiv}) , \IFN_{\nuiv}^{-1/2} \Av_{\nuiv} \rangle 
	+ \langle \IFN_{\nuiv}^{-1/2} \Av_{\nuiv},\IFN_{\nuiv}^{-1/2} \Av_{\nuiv} \rangle \bigr|
	\\
	&=&
	\bigl| \langle \IFN_{\nuiv}^{1/2} (\targvs - \targv_{\nuiv} + \IFN_{\nuiv}^{-1} \Av_{\nuiv}) , 
		\IFN_{\nuiv}^{-1/2} \Av_{\nuiv} \rangle 
	\bigr|
	\\
	& \leq &
	\| \DFN_{\nuiv}^{-1} \, \IFN_{\nuiv} (\targvs - \targv_{\nuiv} + \IFN_{\nuiv}^{-1} \Av_{\nuiv}) \| \, \,
		\| \DFN_{\nuiv} \, \IFN_{\nuiv}^{-1} \Av_{\nuiv} \|
	\leq 
	\frac{3\dltwu_{3}}{4} \, \| \DFN_{\nuiv} \, \IFN_{\nuiv}^{-1} \Av_{\nuiv} \|^{3}
\label{duy7cyf6re6eytedhge2}
\end{EQA}
This and \eqref{gtxdbcbcgdfswerfdsdew} imply \eqref{gtxddfujhyfdytet6ywerfd}.
\end{proof}

\begin{proof}[Proof of Proposition~\ref{PsemiAvex}]
First we bound the variability of \( \IFN_{\nuiv} \) over \( \Nui \). 
Let \( \IFN = \IFN_{\nuivs} \).

\begin{lemma}
\label{L3IFNNui}
Assume \nameref{LLpT3ref}.
Then for any \( \nuiv \in \Nui \)
\begin{EQA}
	\| \DFN^{-1} \, (\IFN_{\nuiv} - \IFN) \, \DFN^{-1} \|
	& \leq &
	\dltwunn \, \| \HFN (\nuiv - \nuivs) \|_{\nano} 
	\, .
\label{8c7c67cc6c63kdldlvvudw}
\end{EQA}
\end{lemma}

\begin{proof}
Let \( \nuiv \in \Nui \).
By \eqref{c6ceyecc5e5etctwhcyegwc}, for any \( \zv \in \R^{\dimp} \), it holds
\begin{EQA}
	&& \!\!\!\!
	\bigl| \bigl\langle \DFN^{-1} \, (\IFN_{\nuiv} - \IFN) \, \DFN^{-1}, \zv^{\otimes 2} \bigr\rangle \bigr|
	=
	\bigl| \bigl\langle \IFN_{\nuiv} - \IFN, (\DFN^{-1} \zv)^{\otimes 2} \bigr\rangle \bigr|
	\\
	& \leq &
	\sup_{t \in [0,1]}
	\bigl| \langle \nabla_{\targv\targv\nuiv}^{3} \fs(\targvs,\nuivs + t (\nuiv - \nuivs)), (\DFN^{-1} \zv)^{\otimes 2} \otimes (\nuiv - \nuivs) \rangle \bigr|
	\leq 
	\dltwunn \, \| \zv \|^{2}  \, \| \HFN (\nuiv - \nuivs) \|_{\nano} \, .
\label{9487654r5tghasdfgsup}
\end{EQA}
This yields \eqref{8c7c67cc6c63kdldlvvudw}.
\end{proof}

The next result describes some corollaries of \eqref{8c7c67cc6c63kdldlvvudw}. 

\begin{lemma}
\label{LIFvari}
Assume \( \DFN^{2} \leq \dmax^{2} \IFN \) and let some other matrix \( \IFN_{1} \in \Matr_{\dimp} \) satisfy
\begin{EQA}
	\| \DFN^{-1} \, (\IFN_{1} - \IFN) \, \DFN^{-1} \|
	& \leq &
	\dmax^{-2} \dltwbun 
\label{icfuifcu7dy7e3gw3ft6e}
\end{EQA}
with \( \dltwbun < 1 \).
Then
\begin{EQA}[ccc]
\label{8c7c67cc6c63kdldlvvudw2}
	\| \IFN^{-1/2} \, (\IFN_{1} - \IFN) \, \IFN^{-1/2} \|
	& \leq &
	\dltwbun 
	\, ,
	\\
	\| \IFN^{1/2} \, (\IFN_{1}^{-1} - \IFN^{-1}) \, \IFN^{1/2} \|
	& \leq &
	\frac{\dltwbun}{1 - \dltwbun} 
	\, ,
\label{g98g75re3gf76r4egu4e}
\end{EQA}
and
\begin{EQA}
	\frac{1}{1 + \dltwbun} \, \| \DFN \, \IFN^{-1} \, \DFN \|
	\leq 
	\| \DFN \, \IFN_{1}^{-1} \DFN \|
	& \leq &
	\frac{1}{1 - \dltwbun} \, \| \DFN \, \IFN^{-1} \, \DFN \|
	\, .
\label{g98g75re3gf76r4egu4ee}
\end{EQA}
Furthermore, for any vector \( \uv \)
\begin{EQA}[rcccl]
\label{uciuvcdiu3eoooodffcy}
	(1 - \dltwbun) \| \DFN^{-1} \IFN \uv \|
	& \leq &
	\| \DFN^{-1} \IFN_{1} \uv \|
	& \leq &
	(1 + \dltwbun) \| \DFN^{-1} \IFN \uv \| \, ,
	\qquad
	\\
	\frac{1 - 2 \dltwbun}{1 - \dltwbun} \| \DFN \, \IFN^{-1} \uv \|
	& \leq &
	\| \DFN \, \IFN_{1}^{-1} \uv \|
	& \leq &
	\frac{1}{1 - \dltwbun} \| \DFN \, \IFN^{-1} \uv \| \, .
	\qquad
\label{uciuvcdiu3eoooodffcym}
\end{EQA}
\end{lemma}

\begin{proof}
Statement \eqref{8c7c67cc6c63kdldlvvudw2} follows from \eqref{icfuifcu7dy7e3gw3ft6e}
because of \( \IFN^{-1} \leq \dmax^{2} \, \DFN^{-2} \).
Define now \( \Uv \eqdef \IFN^{-1/2} \, (\IFN_{1} - \IFN) \, \IFN^{-1/2} \). 
Then \( \| \Uv \| \leq \dltwbun \) and
\begin{EQA}
	\| \IFN^{1/2} \, (\IFN_{1}^{-1} - \IFN^{-1}) \, \IFN^{1/2} \|
	&=&
	\| (\Id + \Uv)^{-1} - \Id \|
	\leq 
	\frac{1}{1 - \dltwbun} \| \Uv \|
\label{uefu83yr83y38frkfbg7e}
\end{EQA}
yielding \eqref{g98g75re3gf76r4egu4e}.
Further,
\begin{EQA}
	\| \DFN \, (\IFN_{1}^{-1} - \IFN^{-1}) \, \DFN \|
	&=&
	\| \DFN \, \IFN_{1}^{-1} \IFN_{1} (\IFN_{1}^{-1} - \IFN^{-1}) \, \IFN \, \IFN^{-1} \DFN \|
	\\
	&=&
	\| \DFN \, \IFN_{1}^{-1} \DFN \, \DFN^{-1} (\IFN_{1} - \IFN) \DFN^{-1} \, \DFN \, \IFN^{-1} \DFN \|
	\\
	& \leq &
	\| \DFN \, \IFN_{1}^{-1} \DFN \| \, \| \DFN \, \IFN^{-1} \DFN \| \, \| \DFN^{-1} (\IFN_{1} - \IFN) \DFN^{-1} \| 
	\leq 
	\dltwbun \| \DFN \, \IFN_{1}^{-1} \DFN \| \, .
\label{odfu7eujgvy55r44rf}
\end{EQA}
This implies \eqref{g98g75re3gf76r4egu4ee}.
Also, by \( \DFN^{2} \leq \dmax^{2} \IFN \)
\begin{EQA}
	\| \DFN^{-1} \IFN_{1} \uv \|
	& \leq &
	\| \DFN^{-1} \IFN \uv \| + \| \DFN^{-1} (\IFN_{1} - \IFN) \DFN^{-1} \DFN \uv \|
	\leq 
	\| \DFN^{-1} \IFN \uv \| + \dltwbun \, \| \DFN \uv \|
	\\
	& \leq & 
	\| \DFN^{-1} \IFN \uv \| + \dltwbun \| \DFN^{-1} \IFN \uv \| 
	\leq 
	(1 + \dltwbun) \| \DFN^{-1} \IFN \uv \| \, ,
	\qquad
\label{uciuvcdiu3ery7httrwbgfcy}
\end{EQA}
and \eqref{uciuvcdiu3eoooodffcy} follows.
Similarly
\begin{EQA}
	\| \DFN \, (\IFN_{1}^{-1} - \IFN^{-1}) \uv \| 
	& = &
	\| \DFN \, \IFN_{1}^{-1} (\IFN_{1} - \IFN) \IFN^{-1} \uv \|
	=
	\| \DFN \, \IFN_{1}^{-1} \DFN \,\, \DFN^{-1} (\IFN_{1} - \IFN) \DFN^{-1} \, \DFN \, \IFN^{-1} \uv \|
	\\
	& \leq &
	\| \DFN^{-1} \, (\IFN_{1} - \IFN) \, \DFN^{-1} \| \,\,
	\| \DFN \, \IFN_{1}^{-1} \, \DFN \| \,\, \| \DFN \, \IFN^{-1} \uv \|
	\leq 
	\frac{\dltwbun}{1 - \dltwbun} \| \DFN \, \IFN^{-1} \uv \|
\label{rtiuru8bttg3hfy7vrgf}
\end{EQA}
and \eqref{uciuvcdiu3eoooodffcym} follows as well.
\end{proof}

Definition \eqref{ygtuesdfhiwdsoif9ruta} implies the following bound.
\begin{lemma}
\label{Lsemidual}
For any \( \nuiv \in \Nui \) and any linear mapping \( \QP \), it holds
\begin{EQA}
	\| \QP \IFT_{\targv\nuiv} (\nuiv - \nuivs) \|
	& \leq &
	\| \QP \IFT_{\targv\nuiv} \, \HFN^{-1} \|_{\dual} \, \| \HFN (\nuiv - \nuivs) \|_{\nano} \, .
\label{difkjv78frhjrjhfh6wekd}
\end{EQA}
\end{lemma}

The next lemma shows that \( \Av_{\nuiv} = \nabla \fs_{\nuiv}(\targvs) \) is nearly linear in \( \nuiv \).

\begin{lemma}
\label{LvarDVa}
Assume \eqref{c6ceyecc5e5etcT2}.
Then
\begin{EQA}
	\| \DFN^{-1} \{ \Av_{\nuiv} + \IFT_{\targv\nuiv} (\nuiv - \nuivs) \} \|
	& \leq &
	\frac{\dltwun}{2} \, \| \HFN (\nuiv - \nuivs) \|_{\nano}^{2} \, .
\label{4cbijjm09j9hhrjfifkrs3}
\end{EQA}
and with \( \crossgrad \) from \eqref{f8vuehery6gv65ftehwee1}, it holds
\begin{EQA}
	\| \DFN^{-1} \Av_{\nuiv} \|
	& \leq &
	\Bigl( \| \DFN^{-1} \IFT_{\targv\nuiv} \, \HFN^{-1} \|_{\dual} 
	+ \frac{\dltwun}{2} \, \| \HFN (\nuiv - \nuivs) \|_{\nano} \Bigr) \, \| \HFN (\nuiv - \nuivs) \|_{\nano}
	\\
	& \leq &
	\crossgrad \, \| \HFN (\nuiv - \nuivs) \|_{\nano} \, .
\label{4cbijjm09j9hhrjfifkrs}
\end{EQA}
\end{lemma}

\begin{proof}
Fix \( \nuiv \in \Nui \) and define \( \av_{\nuiv}(t) \eqdef \Av_{\nuivs + t (\nuiv - \nuivs)} \).
Then \( \av_{\nuiv}(0) = \Av_{\nuivs} = 0 \), \( \av_{\nuiv}(1) = \Av_{\nuiv} \), 
and
\begin{EQA}
	\av_{\nuiv}(1) - \av_{\nuiv}(0) 
	&=&
	\int_{0}^{1} \av'_{\nuiv}(t) \, dt \, ,
\label{dsf7hgduw3heiuwed53ww3}
\end{EQA}
where \( \av'_{\nuiv}(t) = \frac{d}{dt} \av_{\nuiv}(t) \) for \( t \in [0,1] \).
Similarly, by \( \av'_{\nuiv}(0) = - \IFT_{\targv\nuiv} (\nuiv - \nuivs) \), we derive 
\begin{EQA}
	\av_{\nuiv}(1) - \av_{\nuiv}(0) - \av'_{\nuiv}(0)
	&=&
	\int_{0}^{1} (\av'_{\nuiv}(t) - \av'_{\nuiv}(0)) \, dt
	=
	\int_{0}^{1} (1 - t) \, \av''_{\nuiv}(t) \, dt \, ,
\label{g8eh35fgtg7t76868jdy}
\end{EQA}
where 
\( \av''_{\nuiv}(t) = \frac{d^{2}}{dt^{2}} \av_{\nuiv}(t) \).
By condition \nameref{LLpT3ref} 
\begin{EQA}
	\bigl| \langle \av''_{\nuiv}(t), \zv \rangle \bigr|
	&=&
	\bigl| \bigl\langle \nabla_{\targv\nuiv\nuiv}^{3} \fs(\targvs,\nuivs + t (\nuiv - \nuivs)), \zv \otimes (\nuiv - \nuivs)^{\otimes 2} \bigr\rangle \bigr|
	\leq 
	\dltwun \, \| \DFN \zv \| \, \| \HFN (\nuiv - \nuivs) \|_{\nano}^{2} \, 
\label{5cbgcfyt6webhwefuy6}
\end{EQA}
and hence,
\begin{EQA}
	\| \DFN^{-1} \av''_{\nuiv}(t) \|
	&=&
	\sup_{\zv \colon \| \zv \| \leq 1} \bigl| \langle \DFN^{-1} \av''_{\nuiv}(t), \zv \rangle \bigr|
	=
	\sup_{\zv \colon \| \zv \| \leq 1}  
		\bigl| \langle \av''_{\nuiv}(t), \DFN^{-1} \zv \rangle \bigr|
	\\
	& \leq &
	\dltwun \sup_{\zv \colon \| \zv \| \leq 1}
	\| \zv \| \, \| \HFN (\nuiv - \nuivs) \|_{\nano}^{2} \, 
	=
	\dltwun  \| \HFN (\nuiv - \nuivs) \|_{\nano}^{2} \, .
\label{hgvytew3hw36fdfnhrg}
\end{EQA}
This yields
\begin{EQA}
	\| \DFN^{-1} \{ \Av_{\nuiv} + \IFT_{\targv\nuiv} (\nuiv - \nuivs) \} \|
	& \leq &
	\dltwun \, \| \HFN (\nuiv - \nuivs) \|_{\nano}^{2} \, \int_{0}^{1} (1 - t) \, dt
	\leq 
	\frac{\dltwun \, \| \HFN (\nuiv - \nuivs) \|_{\nano}^{2}}{2} \, 
\label{4cbijjm09j9hhrjfifkr}
\end{EQA}
as claimed in \eqref{4cbijjm09j9hhrjfifkrs3}. 
Lemma~\ref{Lsemidual} implies \eqref{4cbijjm09j9hhrjfifkrs}.
\end{proof}

\noindent
Now Proposition~\ref{PbiassemiN} helps to show \eqref{rhDGtuGmusGU0a2spt} and to bound
\( \targv_{\nuiv} - \targvs + \IFN_{\nuiv}^{-1} \Av_{\nuiv} \).

\begin{lemma}
If \( \dltwbun = \dltwunn \| \HFN (\nuiv - \nuivs) \|_{\nano} < 1 \) then
\begin{EQA}
	\| \DFN (\targv_{\nuiv} - \targvs) \|
	& \leq &
	\frac{3}{2} \| \DFN \, \IFN_{\nuiv}^{-1} \Av_{\nuiv} \| 
	\leq 
	\frac{3 \crossgrad}{2(1 - \dltwbun)} \, \| \HFN (\nuiv - \nuivs) \|_{\nano}\,. 
\label{rhDGtuGmusGU0a2sp}
\end{EQA}
Moreover, 
\begin{EQA}
	\| \DFN^{-1} \IFN_{\nuiv} (\targv_{\nuiv} - \targvs + \IFN_{\nuiv}^{-1} \Av_{\nuiv}) \|
	& \leq &
	\frac{3 \crossgrad^{2} \, \dltwu_{3}}{4(1 - \dltwbun)^{2}} \, 
	\| \HFN (\nuiv - \nuivs) \|_{\nano}^{2}
	\, .
	\qquad
\label{d73jejg7tre3fdy3fujell}
\end{EQA}
\end{lemma}

\begin{proof}
By \( \DFN^{2} \leq \IFN \), \eqref{8c7c67cc6c63kdldlvvudw},
\eqref{uciuvcdiu3eoooodffcym} of Lemma~\ref{LIFvari}, and \eqref{4cbijjm09j9hhrjfifkrs} 
\begin{EQA}
	\| \DFN \, \IFN_{\nuiv}^{-1} \Av_{\nuiv} \|
	& \leq &
	\| \DFN \, \IFN_{\nuiv}^{-1} \DFN \| \,\, \| \DFN^{-1} \Av_{\nuiv} \|
	\leq 
	\frac{1}{1 - \dltwbun} \, \| \DFN^{-1} \Av_{\nuiv} \|
	\leq 
	\frac{\crossgrad}{1 - \dltwbun} \, \| \HFN (\nuiv - \nuivs) \|_{\nano} \, .
\label{6dvgycuw76yeficf7wnfo}
\end{EQA}
This and the conditions of the proposition enable us to apply Proposition~\ref{PbiassemiN} which implies \eqref{rhDGtuGmusGU0a2sp} and
\begin{EQA}
	&& \nquad
	\| \DFN^{-1} \IFN_{\nuiv} (\targv_{\nuiv} - \targvs + \IFN_{\nuiv}^{-1} \Av_{\nuiv}) \|
	\leq 
	\frac{3 \dltwu_{3}}{4} \, \| \DFN \, \IFN_{\nuiv}^{-1} \Av_{\nuiv} \|^{2}
	\leq 
	\frac{3 \crossgrad^{2} \, \dltwu_{3}}{4(1 - \dltwbun)^{2}} \, \, \| \HFN (\nuiv - \nuivs) \|_{\nano}^{2}
	\, .
\label{suhyd78wyedw378yedh2}
\end{EQA}
This completes the proof.
\end{proof}

Now we can finalize the proof of the proposition.
With \( \AUv_{\nuiv} \eqdef - \IFT_{\targv\nuiv} (\nuiv - \nuivs) \), it holds
\begin{EQA}
	&& \nquad
	\targv_{\nuiv} - \targvs + \IFN^{-1} \AUv_{\nuiv}
	=
	\targv_{\nuiv} - \targvs + \IFN_{\nuiv}^{-1} \Av_{\nuiv}
	- \IFN_{\nuiv}^{-1} \Av_{\nuiv}
	+ \IFN^{-1} \AUv_{\nuiv} 
	\\
	& = &
	\targv_{\nuiv} - \targvs + \IFN_{\nuiv}^{-1} \Av_{\nuiv}
	- \IFN_{\nuiv}^{-1} (\Av_{\nuiv} - \AUv_{\nuiv}) 
	- (\IFN_{\nuiv}^{-1} - \IFN^{-1}) \AUv_{\nuiv}
	\, .
\label{36gfijh94ejdvtwekoDg}
\end{EQA}
The use of \eqref{8c7c67cc6c63kdldlvvudw} and \eqref{f8vuehery6gv65ftehwee} yields
\begin{EQA}
	&& \! \! \! \! \! \! \! \! 
	\| \DFN^{-1} \IFN_{\nuiv} (\IFN^{-1} - \IFN_{\nuiv}^{-1}) \AUv_{\nuiv} \|
	= 
	\| \DFN^{-1} (\IFN - \IFN_{\nuiv}) \IFN^{-1} \AUv_{\nuiv} \|
	\\
	& \leq &
	\| \DFN^{-1} (\IFN - \IFN_{\nuiv}) \DFN^{-1} \| \,\, 
	\| \DFN \, \IFN^{-1} \DFN \| \,\, \| \DFN^{-1} \IFT_{\targv\nuiv} \HFN^{-1} \|_{\dual} \,\, \| \HFN (\nuiv - \nuivs) \|_{\nano}
	\\
	& \leq &
	\dltwunn \, \| \HFN (\nuiv - \nuivs) \|_{\nano} \, \,
	\crossnorm \, \| \HFN (\nuiv - \nuivs) \|_{\nano}.
	\qquad
\label{jufuf7yfy7tg643yrh}
\end{EQA}
This together with \eqref{4cbijjm09j9hhrjfifkrs} and \eqref{d73jejg7tre3fdy3fujell} implies
\begin{EQA}[c]
	\| \DFN^{-1} \IFN_{\nuiv} (\targv_{\nuiv} - \targvs + \IFN^{-1} \AUv_{\nuiv}) \|
	\leq 
	\Bigl\{ \frac{3 \crossgrad^{2} \, \dltwu_{3}}{4(1 - \dltwbun)^{2}} \, 
	+ \crossnorm \, \dltwunn + \frac{\dltwun}{2}
	\Bigr\} \| \HFN (\nuiv - \nuivs) \|_{\nano}^{2}
	\, .
\label{7tfce7sdfy7w7y73wrytk3f7}
\end{EQA}
The use of \eqref{uciuvcdiu3eoooodffcy} allows to bound
\begin{EQA}[c]
	\| \DFN^{-1} \IFN (\targv_{\nuiv} - \targvs + \IFN^{-1} \AUv_{\nuiv}) \|
	\leq
	\frac{1}{1 - \dltwbun} \, 
	\| \DFN^{-1} \IFN_{\nuiv} (\targv_{\nuiv} - \targvs + \IFN^{-1} \AUv_{\nuiv}) \| 
	\, ,
\label{9jjjjh3d4dd32x32sg3thj}
\end{EQA}
and \eqref{36gfijh94ejdvtwekoise} follows.
\end{proof}

\begin{proof}[Proof of Proposition~\ref{Psemibiassup}]
Let us fix any \( j \leq \dimp \), e.g. \( j=1 \).
Represent any \( \upsv \in \Upsd \) as \( \upsv = (\ups_{1},\nuiv_{1}) \), where \( \nuiv_{1} = (\ups_{2},\ldots,\ups_{\dimp})^{\T} \). 
By \eqref{7tdsyf8iuwopkrtg4576}
\begin{EQA}
	&& \nquad
	\bigl| \DFN_{1} \IFT_{11}^{-1} \IFT_{\ups_{1}\nuiv_{1}}(\nuiv_{1} - \nuivs_{1}) \bigr|
	\leq 
	\bigl| \DFN_{1}^{-1} \IFT_{\ups_{1}\nuiv_{1}}(\nuiv_{1} - \nuivs_{1}) \bigr|
	\leq 
	\crosssup \, \| \HFN_{1} (\nuiv_{1} - \nuivs_{1}) \|_{\infty}
	\, .
\label{dekweivytre34wfgbc}
\end{EQA}
For any \( \nupv_{1} \), define
\begin{EQA}
	\upsd_{1}(\nupv_{1})
	& \eqdef &
	\argmax_{\ups_{1}} \fn(\ups_{1},\nupv_{1}) .
\label{fifwi9jhrio3jiofcwgyt7we}
\end{EQA}
Now we apply Proposition~\ref{Pconcsupp} 
with \( \QP = \DFN_{1} \), \( \IFN = \IFT_{11} \), \( \HFN = \HFN_{1} = \diag(\DFN_{2},\ldots,\DFN_{\dimp}) \),
and \( \| \HFN_{1} (\nuiv_{1} - \nuivs_{1}) \|_{\infty} \leq \rrinf \).
As \( \DFN_{1} \IFT_{11}^{-1} \leq \DFN_{1}^{-1} \), bound \eqref{usdhyw6hikhurnetr256} yields
\begin{EQA}
	&& \nquad
	\bigl| \DFN_{1} \{ \upsd_{1}(\nuiv_{1}) - \upss_{1} - \IFT_{11}^{-1} \AAv_{1} + \IFT_{11}^{-1} \IFT_{\ups_{1}\nuiv_{1}}(\nuiv_{1} - \nuivs_{1}) \} \bigr|
	\\
	& \leq &
	(\dltwuns + \dltwunn) \, \| \HFN_{1} (\nuiv_{1} - \nuivs_{1}) \|_{\nano}^{2}
		+ \bigl(2 \dltwu_{3} + \frac{\dltwunn}{2} \bigr) \, \DFN_{1}^{-2} \AAv_{1}^{2}
	\\
	& \leq &
	(\dltwuns + \dltwunn) \rrinf^{2} + \bigl(2 \dltwu_{3} + \frac{\dltwunn}{2} \bigr) \, \| \DFN^{-1} \AAv \|_{\infty}^{2} 
	=
	\dltwunss \, \| \DFN^{-1} \AAv \|_{\infty}^{2} 
	\, .
	\qquad
\label{usdhyw6hikhurnetr256sup}
\end{EQA}
This, \eqref{dekweivytre34wfgbc}, and \eqref{dhfiejfowelocuyehbrf} imply by \( \rrinf = \cmax \| \DFN^{-1} \AAv \|_{\infty} \)
and \( \cmax = \sqrt{2}/(1 - \crosssup) \)
\begin{EQA}
	\bigl| \DFN_{1} \{ \upsd_{1}(\nupv) - \upss_{1} \} \bigr|
	& \leq &
	\| \DFN^{-1} \AAv \|_{\infty} + \crosssup \, \rrinf
	+ \dltwunss \, \| \DFN^{-1} \AAv \|_{\infty}^{2}
	\\
	& \leq &
	\bigl( 1 + \cmax \, \crosssup + \sqrt{2} - 1 \bigr) \| \DFN^{-1} \AAv \|_{\infty}
	=
	\cmax \| \DFN^{-1} \AAv \|_{\infty}
	= 
	\rrinf \, .
\label{7deje8vg7g76gvgv5f43w}
\end{EQA}
The same bounds apply to each \( j \leq \dimp \).
Therefore, when started from any point \( \upsv \in \Upsd \), the sequential optimization procedure 
which maximizes the objective function w.r.t. one coordinate while keeping all the remaining coordinates
has its solution within \( \Upsd \).
As the function \( \fs \) is strongly concave, and its value improves at every step, the procedure converges to a unique solution \( \upsvd \in \Upsd \).
This implies \eqref{7ytdufchskmls7rghnk}.
Further, with \( \upsv = \upsvd \), 
\begin{EQA}
	\IFT_{11} \{ \ups_{1}(\nuiv_{1}) - \upss_{1} \} + \IFT_{\ups_{1}\nuiv_{1}}(\nuiv_{1} - \nuivs_{1})
	&=&
	\IFT (\upsvd - \upsvs) 
	\, .
\label{tdtdtd6tf78gf98ggt}
\end{EQA}
This yields \eqref{usdhyw6hikhurnetrspp}.
Moreover, with \( \uv \eqdef \DFN (\upsvd - \upsvs) \) and \( B \eqdef \DFN^{-1} \IFT \, \DFN^{-1} \), it holds
\begin{EQA}
	\DFN^{-1} \bigl\{ \IFT (\upsvd - \upsvs) - \AAv \bigr\} 
	&=&
	B (\uv - \DFN \, \IFT^{-1} \AAv) \, ,
\label{fi8fu7v7y6vyv44v3eged}
\end{EQA}
and by \eqref{usdhyw6hikhurnetrspp} and Lemma~\ref{LinvGersh} 
\begin{EQA}
	&& \nquad
	\bigl\| \DFN (\upsvd - \upsvs - \IFT^{-1} \AAv) \bigr\|_{\infty}
	=
	\bigl\| \uv - \DFN \, \IFT^{-1} \AAv \bigr\|_{\infty}
	\\
	&=&
	\bigl\| B^{-1} \DFN^{-1} \{ \IFT (\upsvd - \upsvs) - \AAv \} \bigr\|_{\infty}
	\leq 
	\frac{\dltwunss}{1 - \crosssup} \, \| \DFN^{-1} \AAv \|_{\infty}^{2} 
	\, .
\label{6f6f67vfg3we5dgwytquw}
\end{EQA}
Finally, by \eqref{Thd7w3jhdujrrfrrfgdyw3} of Lemma~\ref{LinvGersh} and by \( \DFN \, \IFT^{-1} \AAv = \BBN^{-1} \DFN^{-1} \AAv  \)
\begin{EQA}[rcccl]
	\| \DFN \, \IFT^{-1} \AAv - \DFN^{-1} \AAv \|_{\infty}
	& = &
	\bigl\| (\BBN^{-1} - \Id_{\dimp}) \DFN^{-1} \AAv \bigr\|_{\infty}
	& \leq &
	\frac{\crosssup}{1 - \crosssup} \| \DFN^{-1} \AAv \bigr\|_{\infty} \, ,
	\\
	\| \DFN \, \IFT^{-1} \AAv - (\Id_{\dimp} + \Delta) \DFN^{-1} \AAv \|_{\infty}
	& = &
	\bigl\| (\BBN^{-1} - \Id_{\dimp} - \Delta) \DFN^{-1} \AAv \bigr\|_{\infty}
	& \leq &
	\frac{\crosssup^{2}}{1 - \crosssup} \| \DFN^{-1} \AAv \bigr\|_{\infty} \, .
\label{6c6tyc6c6cf6de5tw5}
\end{EQA}
This and \eqref{6f6f67vfg3we5dgwytquw} imply the final inequalities of the proposition.
\end{proof}

\begin{lemma}
\label{LinvGersh}
Let \( \BBN = (\BBN_{ij}) \in \Matr_{\dimp} \) with \( \BBN_{ii} = 1 \) and 
\begin{EQA}[c]
	\sup_{\uv \colon \| \uv \|_{\infty} \leq 1} \, \| (\BBN - \Id_{\dimp}) \uv \|_{\infty} \leq \crosssup < 1 .
\label{8f98dej237yt656543n}
\end{EQA}
Then \( \| \BBN \uv \|_{\infty} \leq (1 - \crosssup)^{-1} \| \uv \|_{\infty} \) for any \( \uv \in \R^{\dimp} \).
Similarly,
\begin{EQA}[c]
	\| (\BBN^{-1} - \Id_{\dimp}) \uv \|_{\infty} \leq \frac{\crosssup}{1 - \crosssup} \, \| \uv \|_{\infty}
	\, ,
	\quad
	\| (\BBN^{-1} - 2 \Id_{\dimp} + \BBN) \uv \|_{\infty} \leq \frac{\crosssup^{2}}{1 - \crosssup} \, \| \uv \|_{\infty}
	\, .
	\qquad
\label{Thd7w3jhdujrrfrrfgdyw3}
\end{EQA}
\end{lemma}

\begin{proof}
Represent \( \BBN = \Id_{\dimp} - \Delta \).
Then \eqref{8f98dej237yt656543n} implies \( \| \Delta \uv \|_{\infty} \leq \crosssup \| \uv \|_{\infty} \).
The use of \( \BBN^{-1} = \Id_{\dimp} + \Delta + \Delta^{2} + \ldots \) yields for any \( \uv \in \R^{\dimp} \)
\begin{EQA}[c]
	\| \BBN^{-1} \uv \|_{\infty}
	\leq 
	\sum_{m=0}^{\infty} \| \Delta^{m} \uv \|_{\infty} \, .
\end{EQA}
Further, with \( \uv_{m} \eqdef \Delta^{m} \uv \), it holds
\begin{EQA}
	\| \Delta^{m+1} \uv \|_{\infty}
	&=&
	\| \Delta \uv_{m} \|_{\infty}
	\leq 
	\crosssup \| \uv_{m} \|_{\infty} \, .
\label{cd7dysrtd4ftgewhged4r}
\end{EQA}
By induction, this yields \( \| \Delta^{m} \uv \|_{\infty} \leq \crosssup^{m} \| \uv \|_{\infty} \) and thus,
\begin{EQA}[c]
	\| \BBN^{-1} \uv \|_{\infty}
	\leq 
	\sum_{m=0}^{\infty} \crosssup^{m} \| \uv \|_{\infty}
	= 
	\frac{1}{1 - \crosssup} \, \| \uv \|_{\infty}
\end{EQA}
as claimed. 
The proof of \eqref{Thd7w3jhdujrrfrrfgdyw3} is similar in view of \( \Id_{\dimp} + \Delta = 2 \Id_{\dimp} - \BBN \).
\end{proof}